\documentclass[11pt]{amsart}	
\usepackage[utf8]{inputenc}
\usepackage[T1]{fontenc}
\usepackage{amsmath,amsfonts,amssymb,dsfont}
\usepackage{graphicx} 
\usepackage[english]{babel}
\usepackage{listings}  
\usepackage[colorinlistoftodos]{todonotes}
\usepackage{pgf,pgfkeys,pgfmath}
\usepackage{upgreek} 
\usepackage{fullpage}

\usepackage[backref=page]{hyperref}
\hypersetup{
	plainpages=false,
	unicode=false,          
	pdftoolbar=true,        
	pdfmenubar=true,        
	pdffitwindow=true,      
	pdftitle={My title},    
	pdfauthor={Author},     
	pdfsubject={Subject},   
	pdfcreator={Creator},   
	pdfproducer={Producer}, 
	pdfkeywords={keywords}, 
	pdfnewwindow=true,      
	colorlinks=true,       
	linkcolor=blue,          
	citecolor=blue,        
	filecolor=blue,      
	urlcolor=magenta,          
	urlbordercolor=0 1 1
}

\def\dist{\operatorname{dist}}

\newtheorem{thm}{Theorem}[section]

\def\ackname{Acknowledgements}%
 \def\acknowledgement{\par\addvspace{17pt}\small\rmfamily
\trivlist\if!\ackname!\item[]\else
\item[\hskip\labelsep
{\bfseries\ackname}]\fi}

\usepackage{todonotes}

\newcommand\numeq[1]%
  {\stackrel{\scriptscriptstyle(\mkern-1.5mu#1\mkern-1.5mu)}{=}}

\newcommand{\myqed}{\thinspace\null\nobreak\hfill\hbox{\vbox{\kern-.2pt\hrule
height.2pt depth.2pt\kern-.2pt\kern-.2pt \hbox to2.5mm{\kern-.2pt\vrule
width.4pt \kern-.2pt\raise2.5mm\vbox to.2pt{}\lower0pt\vtop to.2pt{}\hfil
\kern-.2pt \vrule width.4pt\kern-.2pt}\kern-.2pt\kern-.2pt\hrule
height.2pt depth.2pt \kern-.2pt}}\par\medbreak}

\definecolor{lbcolor}{rgb}{0.95,0.95,0.95}
\definecolor{cblue}{rgb}{0.,0.0,0.6}

\definecolor{Lightgray}{rgb}{0.85, 0.85, 0.85}

\def\R{\mathbb{R}}
\def\ds{\displaystyle}

\def\dist{\operatorname{dist}}
\def\ds{\displaystyle}
\def\H{{\mathcal{H}}}
\def\haus{\H^{d-1}}
\def\N{{\mathds{N}}}
\def\e{{\varepsilon}}

\newtheorem{rem}[thm]{Remark}
\newtheorem{prop}[thm]{Proposition}

\def\R{\mathbb{R}}

\usepackage{palatino}

 \title{Approximation of multiphase mean curvature flows\\ with arbitrary nonnegative mobilities}
 \author{Eric Bonnetier}
 \address{Institut Fourier, Universit\'e Grenoble-Alpes, CS 40700, 38058 Grenoble Cedex 09, France}
 \email{eric.bonnetier@univ-grenoble-alpes.fr}
 \author{Elie Bretin}
 \address{Univ Lyon, INSA de Lyon, CNRS UMR 5208, Institut Camille Jordan, 20 avenue Albert Einstein, F-69621 Villeurbanne, France}
 \email{elie.bretin@insa-lyon.fr}
 \author{Simon Masnou}
 \address{Univ Lyon, Universit\'e Claude Bernard Lyon 1, CNRS UMR 5208, Institut Camille Jordan, 43 boulevard du 11 novembre
 1918, F-69622 Villeurbanne, France}
 \email{masnou@math.univ-lyon1.fr}

\makeatletter
\@namedef{subjclassname@2020}{\textup{2020} Mathematics Subject Classification}
\makeatother


\begin{document}
\maketitle
\begin{abstract}
This paper is devoted to the robust approximation with a variational phase field approach of multiphase mean curvature flows with possibly highly contrasted mobilities. The case of harmonically additive mobilities has been addressed recently using a suitable metric to define the gradient flow of the phase field approximate energy.
We generalize this approach to arbitrary nonnegative mobilities using a decomposition as sums of harmonically additive mobilities.
We establish the consistency of the resulting method by analyzing the 
sharp interface limit of the flow: a formal expansion of the phase field
shows that the method is of second order. We propose a simple numerical scheme to approximate the solutions to our new model.
Finally, we present some numerical experiments in dimensions $2$ and $3$
that illustrate the interest and effectiveness of our approach, in particular for approximating flows in which the mobility of some phases is zero.
\end{abstract}

\section{Introduction }

Motion by mean curvature is the driving mechanism  
of many physical systems, in which interfaces are moving due to the 
thermodynamics of phase changes.
Such situations are encountered in the modeling of epitaxial
growth of thin films ~\cite{MR1898515}, in the fabrication of nano-wire by
vapor-liquid-solid growth~\cite{BRETIN2018324,PhysRevMaterials.2.033402}, in the modeling of wetting or 
de-wetting of substrates by  crystalline materials~\cite{doi:10.1021/la500312q,Cahn_Allen_Cahn_angle}, or in
the evolution of grain boundaries in polycrystalline materials~\cite{Mullins1999}.
 
A time-dependent collection of smooth domains $t \rightarrow \Omega(t) \subset \R^d$ is a motion by
mean curvature if, for every $t$,  the normal velocity $V(x,t)$ at each point $x \in \partial \Omega(t)$ is proportional to the mean curvature $H(x,t)$ of $\partial \Omega(t)$ at $x$. 
Up to a time rescaling, the equation of evolution takes the form
\begin{eqnarray*}
V(x,t) &=& H(x,t), \quad x \in \partial \Omega(t),
\end{eqnarray*}
and can be viewed as the $L^2$-gradient flow of the perimeter of $\Omega(t)$
\begin{eqnarray*}
P(\Omega(t)) &=& \H^{d-1}(\partial \Omega(t)),
\end{eqnarray*}
where $\H^{d-1}$ denotes the $(d-1)$-dimensional Hausdorff measure. 
The seminal work of Modica and Mortola~\cite{Modica1977}
has shown that the perimeter can be approximated (in the sense of $\Gamma$-convergence) by
the smooth Van der Waals-Cahn-Hilliard functional
\begin{eqnarray} \label{P_eps}
P_\e(u) &=& \ds\int_Q \left( \ds\frac{\e}{2} |\nabla u|^2 + \ds\frac{1}{\e} W(u) \right) \,dx.
\end{eqnarray}
defined for smooth functions $u$, with  $Q \subset \R^2$ a fixed bounded domain {that contains strictly the convex envelope of $\Omega(0)$ (so that $\partial \Omega(t)$ stays at positive distance from $\partial Q$)}, 
$\e>0$ a small parameter, and $W$ a smooth double-well potential,
typically
\begin{eqnarray*}
W(s) &=& \ds\frac{1}{2}s^2(1-s)^2.
\end{eqnarray*}
It follows from Modica-Mortola's $\Gamma$-convergence result~\cite{Modica1977} that, when $\Omega$ is a set of finite perimeter,
its characteristic function $1_{\Omega}$ can be approximated in $L^1$ by sequences of functions of the form $u_\e = q(\dist(x,\Omega)/\e)$,
such that $P_\e(u_\e) \rightarrow c_W P(\Omega)$, with
$c_W = \int_0^1 \sqrt{2 W(s)} ds$. Here, $\dist(x,\Omega)$ denotes the signed distance to the set $\Omega$ (negative inside, positive outside), and $q$ is a so-called {\it optimal profile} that depends on the potential $W$ and is defined by
$$
q =  \mathop{\rm argmin} \limits_{p } 
\left\{ \int_{\R} \sqrt{ {W}(p(s))} \, |p'(s)| ds,\;  
p(-\infty) = 1,\, p(0) = 1/2,\, p(+\infty) = 0 
\right\},
$$  
where $p : \R \to \R$ ranges over the set of Lipschitz continuous functions.
A simple derivation of the Euler equation associated with this minimization problem shows that
\begin{equation} \label{def_q}
 q'(s) = -\sqrt{2 {W}(q(s))} \quad\text{and}\quad q''(s) = W'(q(s)), \quad  
 \text{ for all } s \in \R,
\end{equation}
which implies that $q(s) = (1- \tanh(s))/2$ in the case of the standard double-well potential 
${W}(s) = \frac1 2 s^2 (1-s)^2 $ considered above.

The $L^2$-gradient flow of the Van der Waals--Cahn--Hilliard energy $P_\varepsilon$, 
gives the Allen-Cahn equation \cite{Allen1979}.
Up to a time rescaling, it takes the form
\begin{eqnarray}\label{eq_AllenCahn}  
u_t &=&  \Delta u - \frac{1}{\varepsilon^2}W'(u).
\end{eqnarray}
Given smooth initial and boundary conditions, this nonlinear parabolic equation has a unique solution in short time which satisfies a comparison
principle~ \cite{Ambrosio2000}. Furthermore, a smooth motion by mean curvature $t\mapsto \Omega(t)$ can be approximated by
$$ 
\Omega^\varepsilon(t) 
= \left\{x \in \R^d,\; u^\varepsilon(x,t) \geq \frac{1}{2} \right\},
$$
where $u^\varepsilon$ solves~\eqref{eq_AllenCahn} with initial condition
$$u^{\varepsilon}(x,0) = q \left( \frac{\dist(x,\Omega(0))}{\varepsilon}\right).
$$
A formal asymptotic expansion of $u^{\varepsilon}$ near the boundary $\partial \Omega^\varepsilon(t)$ 
shows~\cite{BellettiniBook2013}  that $u^{\varepsilon}$ is quadratically close to the optimal profile, i.e.
$$ 
u^{\varepsilon}(x,t) = q \left( \frac{\dist(x,\Omega^{\varepsilon}(t))}{\varepsilon}\right) + O(\varepsilon^2),
$$
and the normal velocity $V^{\varepsilon}$ along $\partial \Omega^\varepsilon(t)$  satisfies
$$ V^{\varepsilon} = H_{\partial \Omega^\varepsilon(t)} + O(\varepsilon^2).$$
Convergence of $\partial \Omega_{\varepsilon}(t)$ to $\partial \Omega(t)$ has been rigorously
proved for smooth flows with a quasi-optimal 
convergence order $O(\varepsilon^2 | \log \varepsilon |^2)$~\cite{Chen1992, Mottoni1995, Bellettini1995} 
The fact that $u^{\varepsilon}$ is quadratically close to the optimal profile
has inspired the development of very effective numerical methods
~\cite{Modica1977,Chen1992,BenceMerrimanOsher,Ruuth_efficient,BrasselBretin2011}.


\subsection{Multiphase flows}

In the presence of several phases, the motion of interfaces obeys a relation of the form
\begin{eqnarray*}
\ds V_{ij} &=& m_{ij}\sigma_{ij} H_{ij},
\end{eqnarray*}
where $V_{ij}, H_{ij}$, and $\sigma_{ij}$ denote, respectively, the normal velocity,
the mean curvature and the surface tension along the interface $\Gamma_{ij}$
that separates the phases $i$ and $j$.
The mobilities $m_{ij}$ describe how fast adatoms from one phase may be adsorbed 
in another phase as the front advances.
These parameters are associated with the kinetics of the moving front,
not with the equilibrium shape of the crystal, contrarily to the
surface tensions $\sigma_{ij}$.
\medskip
 
Assuming that the material phases partition an open region $Q \subset \R^n, n = 2,3$ 
into closed sets $\Omega_i$ occupied by the phase $i$,
the perimeter functional takes the form
\begin{eqnarray*}
   P(\Omega_1, \Omega_2, \dots,\Omega_N) &=& 
   \frac{1}{2} \sum_{1\leq i < j\leq N}  \sigma_{ij}  \haus(\Gamma_{ij}\cap Q).
\end{eqnarray*}
with $\Gamma_{ij}=\partial \Omega_i\cap \partial \Omega_j$. We assume throughout this work that the surface tensions are additive, i.e., that there
exist $\sigma_k\geq 0$, $1 \leq k \leq N$, such that
\begin{eqnarray*}
\sigma_{ij} &=& \sigma_i + \sigma_j, \quad 1 \leq i < j \leq N.
\end{eqnarray*}
The additivity property is always satisfied when $N \leq 3$ and when the set of coefficients $\{\sigma_{ij}\}$
satisfy the triangle inequality.
In particular, this is the case of the evolution of a single chemical species
in its liquid, vapor and solid phases.
The perimeter functional can then be rewritten in the form
\begin{eqnarray*}
   P(\Omega_1, \Omega_2, \dots,\Omega_N) &=&  
   \sum_{i}^{N}  \sigma_{i}  \haus(\partial \Omega_i\cap Q),
\end{eqnarray*}
and therefore lends itself to approximation by the multiphase Cahn-Hilliard energy defined for 
${\boldsymbol u} = (u_1, u_2, \dots, u_N)$ by
$$P_{\varepsilon}({\boldsymbol u}) =
\begin{cases}
\frac{1}{2}\ds \sum_{i=1}^{N} \int_{Q} \sigma_i \left( \varepsilon \frac{|\nabla u_i|^2}{2} + \frac{1}{\varepsilon} W(u_i) \right) dx, &\text{if~} \sum_{i=1}^{N} u_i = 1 \text{,} \\
+ \infty & \text{otherwise.}
\end{cases}
$$
Modica-Mortola's scalar $\Gamma$-convergence result was generalized to multiphase in~\cite{oudet2011} when $\sigma_i = 1, 1 \leq i \leq N$.
For more general $\Gamma-$convergence results, we refer to~\cite{baldo1990,Bretin_Masnou_multiphase} 
for inhomogeneous surface tensions, and to~\cite{Garcke_amulti,Garcke199887} for anisotropic
surface tensions. 
\medskip

The $L^2$-gradient flow of $P_\varepsilon$ yields the following system of Allen-Cahn equations:
\begin{equation} \label{eqn:allencahn_simple}
 \partial_t u^{\varepsilon}_k = \sigma_k \left[ \Delta u^{\varepsilon}_k - \frac{1}{\varepsilon^2} W'(u^{\varepsilon}_k) \right] + \lambda^{\varepsilon}, \quad   \forall k=1, \dots, N \text{,}
\end{equation}
where the Lagrange multiplier $\lambda^{\varepsilon}$ accounts for the constraint 
$\sum_{k=1}^{N} u^{\varepsilon}_k = 1$.
In practice, however, the numerical schemes derived from~(\ref{eqn:allencahn_simple})
do not prove as accurate as in the single-phase case.
To improve the convergence, one may localize the Lagrange multiplier $\lambda$
near the diffuse interface, as was proposed in~\cite{MR3961089},
and consider instead of~\eqref{eqn:allencahn_simple} the modified system
\begin{equation} \label{eqn:allencahn_modif}
  \partial_t u^{\varepsilon}_k = \sigma_k 
  \left[ \Delta u^{\varepsilon}_k - \frac{1}{\varepsilon^2} W'(u^{\varepsilon}_k) \right] + \lambda^{\varepsilon} \sqrt{2W(u_k)} \quad   \forall k=1, \dots, N,
\end{equation}
where the effect of $\lambda$ is essentially felt in the vicinity of the interfaces.
A rigorous proof of convergence of this modified Allen-Cahn system to
multiphase Brakke's mean curvature flow is established in~\cite{Takasaopreprint}.


\subsection{Incorporating mobilities} 

As mentioned above, mobilities are kinetic parameters that model how fast
adatoms get attached to an evolving front. 
In~\cite{Garcke199887,Garcke_amulti}, mobilities are included in the
definition of the surface potential $f({\boldsymbol u},\nabla {\boldsymbol u})$ 
and of the multi-well potential ${\bf W}({\boldsymbol u})$
that define the Allen-Cahn approximate energy
${\boldsymbol u} = (u_1, u_2, \dots, u_N)$ as
$$P_{\varepsilon}({\boldsymbol u}) =
\begin{cases}
\int_{Q} \varepsilon f({\boldsymbol u},\nabla {\boldsymbol u}) 
+ \frac{1}{\varepsilon} {\bf W}({\boldsymbol u}) dx, &\text{if~} \sum_{i=1}^{N} u_i = 1 \text{,} \\
+ \infty & \text{otherwise.}
\end{cases}
$$

Examples of surface potential $f$ and multiple well potential ${\bf W}$ 
that have been considered are
\[ \left\{ \begin{array}{ccl}
f({\boldsymbol u},\nabla {\boldsymbol u}) &=& 
\sum_{i<j} m_{ij} \sigma_{ij} \left| u_i \nabla u_j - u_j \nabla u_i \right|^2,
\\ [6pt]
{{\bf W}}({\boldsymbol u}) &=& 
\sum_{i<j} \frac{\sigma_{ij}}{m_{ij}} u_i^2 u_j^2  + \sum_{i<j<k} \sigma_{ijk} u_i^2 u_j^2 u_k^2.
\end{array} \right.
\]
In these models, both surface tensions and mobilities appear in the Cahn-Hilliard energy.
It is shown in~\cite{MR1740846,Garcke_amulti} that taking the sharp interface limit imposes
constraints on the limiting values of the surface tensions and mobilities, in particular
in the anisotropic case. From a numerical perspective
it follows that the mobilities are likely to impact the size of the diffuse interfaces, 
as they appear in the energy, especially in situations where the contrast of mobilities is large.
\medskip

In this work, we assume that the flux of adatoms is a linear function of the normal velocity 
of the interface $\Gamma_{ij}$, with a proportionality constant equal to $m_{ij}$. 
From the modeling point of view, this amounts to considering 
the surface tensions as geometric parameters which govern the equilibrium,
and the mobilities as parameters related to the evolution of the system from 
an out-of-equilibrium configuration, 
which only affect the metric used for the gradient flow.

It is proposed in~\cite{BRETIN2018324} to take the mobilities into account through the metric used to define the gradient flow.
The mobilities that are considered in~\cite{BRETIN2018324} mimic the properties of 
additive surface tensions, i.e., it is assumed that the $m_{ij}$'s, for $1\leq i<j\leq N$, can be
decomposed as
\begin{equation} \label{def:mob_harmo}
\frac{1}{m_{ij}} = \frac{1}{m_i} + \frac{1}{m_j},
\end{equation}
for a suitable collection of coefficients $m_k > 0, 1 \leq k \leq N$.  We extend the definition of $m_{ij}$ to all $1\leq i,j\leq N$, $i\not=j$, by the natural symmetrization $m_{ji}=m_{ij}$. 

The Allen-Cahn system associated to a set of mobilities with such a decomposability property takes the form
\begin{equation} \label{eqn:allencahn_mob_additive}
\partial_t u^{\varepsilon}_k =   
m_k \left[ \sigma_k \left( \Delta u^{\varepsilon}_k - \frac{1}{\varepsilon^2} W'(u^{\varepsilon}_k) \right) 
+  \lambda^{\varepsilon} \sqrt{2 W(u_k)} \right], \quad \forall k \in \{1,\dots,N\},
\end{equation}
where the Lagrange multiplier $\lambda^{\varepsilon}$ is again associated 
to the constraint  $\sum u_i^{\varepsilon} = 1$ and given by
$$ 
\lambda^{\varepsilon} = - \frac{ \sum_k m_k  \sigma_k \left( \Delta u^{\varepsilon}_k - \frac{1}{\varepsilon^2} W'(u^{\varepsilon}_k)  \right)}{\sum_k m_k \sqrt{2 W(u_k)}}.
$$

This model has the following advantages~\cite{BRETIN2018324}: 
\begin{itemize}
\item 
It is quantitative in the sense that the coefficients $\sigma_i$ and $m_i$ can be identified from the mobilities and surface tensions $m_{ij}$ and $\sigma_{ij}$,
\item 
Numerical tests indicate an accuracy of order two in $\varepsilon$, and suggest  that the size of the diffuse interface does not depend on the ${m_{ij}}$'s,
\item 
A simple and effective numerical scheme can be derived to approximate the
solutions to \eqref{eqn:allencahn_mob_additive}.
\end{itemize}
Positive mobilities that satisfy \eqref{def:mob_harmo} are called {\it harmonically additive}.
For convenience, we extend the definition to nonnegative mobilities using the convention $\frac 1{0^+}=+\infty$ and $\frac 1{+\infty}=0^+$. The Allen-Cahn equation associated with a null coefficient $m_i=0$ reduces to $\partial_t u^{\varepsilon}_i = 0$.


\subsection{General mobilities} 

The main motivation of the paper is to introduce a phase field model 
similar to \eqref{eqn:allencahn_mob_additive}, but not limited to 
harmonically additive mobilities.
For example, in the case of a 3-phase system ($N=3$), the triplet of mobility coefficients
$(m_{12},m_{13},m_{23}) = (1,0,0)$ is indeed harmonically additive as one can choose
$m_1 = m_2 = 2$ and $m_3 = 0$. 
However, this is far from general, and there seems to be no physical (even practical)
reason that justifies this hypothesis. 
The situation studied in~\cite{BRETIN2018324},  
that models the vapor-liquid-solid (VLS) growth of nanowires, 
is an illustration of this remark. 
Indeed, VLS growth can be viewed as a system of three phases
with mobilities $m_{LS} = m_{LV} = 1, m_{SV} = 0$. In such a system, 
the vapor-solid interface remains fixed, as growth only takes place
along the liquid-solid interface.
It is easy to check that a triplet of mobilities of the form 
$(m_{12},m_{13},m_{23}) = (1,1,0)$ fails to be harmonically additive
(or more generally, any triplet $(1,1,\beta)$ as soon as $\beta<1/2$).
\medskip

To derive a numerical scheme adapted to general nonnegative mobilities and ensuring that the width of the diffuse interface does
not depend on the possible degeneracy of the mobilities, 
we decompose each mobility as a sum of  harmonically additive mobilities.
In other words, for each $m_{ij}$, $i\not=j$, we consider $P \in \N$ and 
{nonnegative} coefficients $\{m^{p}_{ij}\}$ and $\{m^p_{i}\}$ such that  
\begin{eqnarray}\label{m-additive}
 m_{ij} = \sum_{p=1}^{P} m_{ij}^p\quad
 \text{ and } \quad
\frac{1}{m^p_{ij}} = \frac{1}{m^p_i} + \frac{1}{m^p_j},
\end{eqnarray}
{with the convention that $\frac{1}{0^+} = + \infty$ and $\frac{1}{+\infty} = 0^+$.}\\

It is easy to check that one can {\it always} find such a decomposition, 
provided  all the $m_{ij}$'s are nonnegative.
For instance, a canonical choice is  
\begin{eqnarray}\label{decomp_m}
m_{ij} &=& \sum_{1 \leq k < \ell \leq N} m_{ij}^{k\ell},
\end{eqnarray}
with {$m_{ij}^{k\ell} =m_{ij} \delta_k(i) \delta_\ell(j)$} satisfying
\begin{eqnarray*}
\ds\frac{1}{m_{ij}^{k\ell}} &=& \ds\frac{1}{m_i^{k\ell}} + \ds\frac{1}{m_j^{k\ell}},
\end{eqnarray*}
where, for every $1 \leq \alpha \leq N$, $m_\alpha^{k\ell} = 2m_{k\ell}(\delta_k(\alpha) + \delta_\ell(\alpha))$.

We associate to this decomposition a phase field model of the form
\begin{equation} \label{eqn:allencahn_mob_additive2}
 \partial_t u^{\varepsilon}_k =   m^{*}_k \left[ \sigma_k \left( \Delta u^{\varepsilon}_k - \frac{1}{\varepsilon^2} W'(u^{\varepsilon}_k) \right) +  \lambda^{\varepsilon}_k \sqrt{2 W(u_k)} \right], \quad \forall k \in \{1,\dots,N\}
\end{equation}
where we define
\begin{itemize}
 \item the coefficients $m^{*}_k$ as 
 $$ 
 m^{*}_k = \sum_{p=1}^{P} m^p_{k}.
 $$
 \item the Lagrange multipliers $\lambda^{\varepsilon}_k$ as
 $$ 
 \lambda^{\varepsilon}_k = \frac{1}{m^{*}_k} \sum_{p=1}^{P} m_{k}^{p}~\lambda^{p,\varepsilon}, \quad \text{ with } \quad  \lambda^{p,\varepsilon} = - \left( \frac{ \sum_{k=1}^{N} m_k^p  \sigma_k \left( \Delta u^{\varepsilon}_k - \frac{1}{\varepsilon^2} W'(u^{\varepsilon}_k)  \right)}{\sum_{k=1}^{N} m_k^p \sqrt{2 W(u_k)}} \right).
 $$
\end{itemize}

\begin{rem} 
The difference between the two models \eqref{eqn:allencahn_mob_additive} and  \eqref{eqn:allencahn_mob_additive2}  
lies in the  definition of the Lagrange multipliers $\lambda^{\varepsilon}_k$. 
In the first case, the components $\lambda^{\varepsilon}_k$ are identical 
and do not differentiate interfaces according to the mobilities
for the satisfaction of the constraint $\sum \partial_t u_k = 0$. 
In the second model, the $\lambda^{\varepsilon}_k$'s are weighted in terms
of the $m^p_k$'s.
\end{rem}
\begin{rem}
There is, in general, no unique way of decomposing a given set of nonnegative mobilities 
$(m_{ij})_{1 \leq i<j \leq N}$ as a sum of harmonically additive mobilities.
In view of the tests we performed, it seems that the particular choice of 
decomposition does not have a strong influence on the numerical results.
\end{rem}

Proving the consistency of our new phase field model \eqref{eqn:allencahn_mob_additive2}  is the
main theoretical result of the present work.
More precisely, we show that smooth solutions to the above system are
close up to order 2 in $\varepsilon$ to a sharp interface motion.

\begin{prop} \label{claim:allencahn_mob_additive2} 

Assume that ${\boldsymbol u}^{\varepsilon}$ is a smooth solution 
to~\eqref{eqn:allencahn_mob_additive2} and define  the set 
$$E^{\varepsilon}_i(t) = \{ x \in \Omega,\quad u_i^\varepsilon(x,t) \geq 1/2 \}$$
and the interface
$$\Gamma_{ij}^{\varepsilon}(t) =  \partial E^{\varepsilon}_i(t) \cap \left\{  x \in \Omega,\quad u_j^\varepsilon(x,t) \geq  u_k^\varepsilon(x,t) \quad  \forall k \neq i \right\}.$$
Then, in a neighborhood of  $\Gamma_{ij}^{\varepsilon}$,  ${\boldsymbol u}^{\varepsilon}$ satisfies
$$
\begin{cases}
u_i^{\varepsilon} &=
q \left( \frac{ d_i^{\varepsilon}(x,t)}{\varepsilon} \right) + O(\varepsilon^2),
\\
u_j^{\varepsilon} &=  1 - q \left( \frac{d_i^{\varepsilon}(x,t)}{\varepsilon} \right) + O(\varepsilon^2),\\
u_k^{\varepsilon} &=  O(\varepsilon^2), \text{ for k } \in \{1, \dots, N\}\setminus\{i,j\},
\end{cases}
$$
where    $d^{\varepsilon}_i(x,t)$ denotes the signed distance to $E^{\varepsilon}_i(t)$, with 
$d^{\varepsilon}_i(x,t) < 0$ if $x \in E^{\varepsilon}_i(t)$. Define further $V_{ij}^\varepsilon(x,t) = \partial_t d^{\varepsilon}_i(x,t)$
for $x \in \Gamma_{ij}$.
Then the following estimate holds:
\begin{eqnarray*}
 V_{ij}^{\varepsilon} &=& m_{ij} \sigma_{ij} H_{ij} + O(\varepsilon).
\end{eqnarray*}
\end{prop}
\medskip

The paper is organized as follows: Proposition~\ref{claim:allencahn_mob_additive2} 
is proven formally in Section 2, using the method of matched asymptotic expansions (the formal proof is given for general nonnegative mobilities, thus including of course the more restrictive case of harmonically additive mobilities considered in~\cite{BRETIN2018324}).
In Section $3$, we propose a numerical scheme based on the phase-field
system~\eqref{eqn:allencahn_mob_additive2}. 
To illustrate its simplicity, we give an explicit {\bf Matlab} 
implementation of the scheme in dimension $2$ that requires less than $50$ lines. 
In the last section, we provide examples of simulations of multiphase flows
in dimensions $2$ and $3$ that illustrate the consistency and effectiveness of the method, and 
the influence of mobilities on the flow.


\section{Asymptotic expansion of solutions to the Allen-Cahn system} 

This section is devoted to the formal identification of sharp interface limits of solutions
${\boldsymbol u}^\varepsilon=(u_1^\varepsilon,\dots,u^{\varepsilon}_{N})$ to the Allen-Cahn 
system \eqref{eqn:allencahn_mob_additive2}.
To this aim, we use the method of matched asymptotic expansions
proposed in  \cite{fife,pego,belpaoqo,Loreti_march,MR3961089,BRETIN2018324}, 
which we apply around each interface $\Gamma_{ij}$. 
Henceforth,  we fix $i,j \in \{1,\dots, N \}$ and we assume that ${\boldsymbol u}^\varepsilon$ is a solution 
to~\eqref{eqn:allencahn_mob_additive2} that is smooth near  the 
interface $\Gamma^{\varepsilon}_{ij}$.

\subsection{Preliminaries}

\paragraph{Outer expansion far from $\Gamma_{ij}^{\varepsilon}$} 

We assume that the {\it outer expansion} of $u^{\varepsilon}_{k}$, {i.e.} the expansion far from the front $\Gamma^{\varepsilon}_{ij}$, has the form:
$$u^{\varepsilon}_k(x,t) = u^0_k(x,t) + \varepsilon u^1_k(x,t) +  O(\varepsilon^2), \text{ for all } k \in \{1, \dots,N\}.$$
In particular, and analogously to \cite{Loreti_march}, it is not difficult to see that 
if $E^\varepsilon_i(t) = \{ x \in \Omega, u^{\varepsilon}_i \geq \frac{1}{2}\},$ then
$$u^0_i(x,t) =
\begin{cases}
 1 & \text{ if } x \in E^\varepsilon_i(t) \\
 0 & \text{otherwise}
\end{cases}, \quad  u^0_j(x,t) =
\begin{cases}
 0 & \text{ if } x \in E^\varepsilon_i(t) \\
 1 & \text{otherwise}
\end{cases}
$$
and $ u^1_i = u^1_j = 0,$ $u^0_k = u^1_k = 0$ for all $k \in \{1,\dots, N\}\setminus \{i,j\}$. \\

\paragraph{Inner expansions around $\Gamma_{ij}^{\varepsilon}$}
In a small neighborhood of $\Gamma_{ij}^{\varepsilon}$, 
we define the stretched normal distance to the front as
$ z = \frac{1}{\varepsilon}{d^\varepsilon_i(x,t)},$ where 
$d^\varepsilon_i(x,t)$ denotes the signed distance to $E^\varepsilon_i(t)$ 
such that $d^\varepsilon_i(x,t)<0$ in $E^{\varepsilon}_i(t)$.
The {\it inner expansions} of $u^{\varepsilon}_k(x,t)$ and $\lambda^{p,\varepsilon}(x,t)$, {i.e.} 
expansions close to the front, are assumed of the form
$$ u^{\varepsilon}_k(x,t) = U^{\varepsilon}_{k}(z,x,t) = U^{0}_{k}(z,x,t) 
+ \varepsilon U^{1}_k(z,x,t) + O(\varepsilon^2),  
\text{ for all } k \in \{1, \dots,N\},$$
and
$$  \lambda^{p,\varepsilon}(x,t) =  \Lambda^{p,\varepsilon}(z,x,t)  
=  \varepsilon^{-2} \Lambda^{p,-2}(z,x,t) +  \varepsilon^{-1} \Lambda^{p,-1}(z,x,t) + O(1).$$
Moreover, if $n$ denotes the unit normal to $\Gamma_{ij}$ and $V^{\varepsilon}_{ij}$ 
the normal velocity to the front (pointing to the inside of $E^{\varepsilon}_i$) for $x\in\Gamma_{ij}$
$$ V^{\varepsilon}_{ij} = \partial_t d^\varepsilon_i(x,t) 
= V^{0}_{ij} + \varepsilon V^1_{ij} + O(\varepsilon^2) , 
\quad n = \nabla d^\varepsilon_i(x,t).$$
where $\nabla$ refers to the spatial derivative only.
Following~\cite{pego,Loreti_march} we assume that
$U^{\varepsilon}_k(z,x,t)$ does not change when $x$ varies normal to $\Gamma_{ij}$ with $z$ held fixed, or equivalently
$( \nabla U^{\varepsilon}_k )_{z={\rm const.}}\cdot n = 0$. This amounts to requiring that the blow-up with respect to the parameter $\varepsilon$ is consistent with the flow.

Following~\cite{pego,Loreti_march}, it is easily seen that
\begin{equation}\label{eq_scaling}
\begin{cases}
 \nabla u^\varepsilon_k = \nabla_x U^\varepsilon_k + \varepsilon^{-1} n \partial_z U^\varepsilon_k, \\
 \Delta u^\varepsilon_k =\Delta_x U^\varepsilon_k 
 + \varepsilon^{-1} \Delta d_i ~\partial_z U^\varepsilon_k 
 + \varepsilon^{-2} \partial^2_{zz} U^\varepsilon_k, \\
 \partial_t u^\varepsilon_k = \partial_t U^\varepsilon_k 
 {+} \varepsilon^{-1} V_{ij}^{\varepsilon} \partial_z U^\varepsilon_{k}.
\end{cases}
\end{equation}
Recall also that in a sufficiently small neighborhood of $\Gamma_{ij}$, according to Lemma 14.17 in~\cite{GT}, we have
$$ \Delta d_i (x,t) = \sum_{k=1}^{d-1} \frac{\kappa_k(\pi(x))}{1 + \kappa_k(\pi(x)) d_i(x,t)} = \sum_{k=1}^{d-1} \frac{\kappa_k(\pi(x))}{1 + \kappa_k(\pi(x)) \varepsilon z },$$
where $\pi(x)$ is the projection of $x$ on $\Gamma_{ij}$ and $\kappa_k$ are the principal curvatures on $\Gamma_{ij}$.
In particular this implies that
$$ \Delta d^\varepsilon_i (x,t) = H_{ij} - \varepsilon z \|A_{ij}\|^2 + O(\varepsilon^2), $$
where $H_{ij}$ and $\|A_{ij}\|^2$ denote, respectively, the mean curvature and the squared $2$-norm of the second fundamental form of $\Gamma_{ij}$ at $\pi(x)$.

\paragraph{{Matching conditions between outer and inner expansions}} 
The matching conditions (see \cite{Loreti_march} for more details) can be written as:
$$   \lim_{z \to + \infty} U^0_i(z,x,t) = 0, \lim_{z \to - \infty} U^0_i(z,x,t) = 1, \quad \lim_{z \to \pm \infty} U^1_i(z,x,t) = 0,$$
$$   \lim_{z \to + \infty} U^0_j(z,x,t) = 1, \lim_{z \to - \infty} U^0_j(z,x,t) = 0, \quad \lim_{z \to \pm \infty} U^1_j(z,x,t) = 0,$$
and
$$    \lim_{z \to \pm \infty} U^0_k(z,x,t) =  \lim_{z \to \pm \infty} U^1_k(z,x,t) = 0, \text{for all } k \in \{1,\dots,N\}\setminus \{i,j\}.$$
{Moreover, recall that the definition of $z$ implies that $U^0_i(0,x,t) = \frac{1}{2}$ and $U^1_i(0,x,t) = U^2_i(0,x,t) = 0$.}

\subsection{Analysis of the Allen-Cahn system }

We insert the form~\eqref{eq_scaling} in~\eqref{eqn:allencahn_mob_additive2} and
match the terms according to their powers of $\varepsilon$.

\paragraph{{Order $\varepsilon^{-2}$}} Identifying the terms of order $\varepsilon^{-2}$ gives for all  $k \in \{1,\dots,N\}$:
$$\sigma_k \left( \partial^2_{zz} U^0_k - W'(U^0_k) \right)  + \frac{1}{m_k^*} \sum_{p} m_k^{p} \Lambda^{p,-2} \sqrt{2 W(U^0_k)} = 0,$$
and
$$\left[\sum_{k=1}^{N} m^p_k \sqrt{2 W(U^0_k)} \right] \Lambda^{p,-2}  = - \sum_{k=1}^{N} m^p_k \sigma_k \left(  \partial^2_{zz} U^0_k - W'(U^0_k)  \right).$$
{
Assuming  $\Lambda^{p,-2} = 0$ shows that  $ \partial^2_{zz} U^0_k - W'(U^0_k)  = 0$ for all $k \in \{1, \dots, N\} $.
Now, using boundary conditions, we deduce that  $U^{0}_k(z,x,t) = 0$  for all $k \in \{1, \dots, N\} \setminus \{i,j\}$ as 
$  \lim_{z \to \pm \infty} U^0_k(z,x,t) = 0$. About the case $k=i$, recall that the phase field profile
$q$,  defined as the solution of   $ q''(z) =  W'(q) $ with $\lim_{z \to + \infty} q(z) = 0$, $\lim_{z \to - \infty} q(z) = 1$ and 
$q(z) = 1/2$, satisfies $q(z) = (1 - \text{tanh}(z))/2$. Now by remarking that  $U^0_i(0,x,t) = \frac{1}{2}$, 
$\lim_{z \to + \infty} U^0_i(z,x,t) = 0$ and $\lim_{z \to - \infty} U^0_i(z,x,t) = 1$, we show that $U^0_i(z,x,t)  = q(z)$. The function 
$U^{0}_{j}$ can then be  identified to  $U^{0}_{j} = 1-q(z) = q(-z)$ thanks to the partition constraint $\sum_{k=1}^{N} U^0_k(z,x,t) = 1$.
  }

\paragraph{{Order $\varepsilon^{-1}$}} 
Matching the next order terms shows that for $k \neq \{1,\dots, N\}\setminus {i,j}$,
$$  \frac{1}{m^*_k} V^0_{ij} \partial_z U^0_{k} =  \sigma_k \left[ \partial^2_{zz} U^{1}_k - W''(U^0_k)U^1_k  +   H_{ij}  \partial_z U^0_{k} \right]  +  \frac{1}{m_k^*} \sum_{p} m_k^{p} \Lambda^{p,-1} \sqrt{2 W(U^0_k)}  $$
and
$$ \left[\sum_{k=1}^{N}  m^p_k  \sqrt{2 W(U^0_k)} \right] \Lambda^{p,-1}  = - \sum_{k=1}^{N}  m_k^p \sigma_k \left[  \partial^2_{zz} U^1_k - W''(U^0_k) U^1_k + H_{ij} \partial_z U^0_{k} \right].$$
Then, for all $k \in \{1,\dots, N \}\setminus\{i,j\}$,  {as  $U^0_k = 0$, we deduce that 
$ \left( \partial^2_{zz} U^{1}_k - W''(0) U^1_k \right)=0$  which, in view of the matching boundary conditions $  \lim_{z \to \pm \infty} U^1_k(z,x,t) = 0$,  yields  $U^{1}_k = 0.$ }\\

Moreover, recalling from \eqref{def_q} that  $\sqrt{2 W(q(z))} = -q'(z) $, the  equations for $U_i^1$  and $\Lambda^{p,-1}$ become
\begin{equation*}
 V^0_{ij} q'(z) =  \sigma_i m_i^* \left( \partial^2_{zz} U^{1}_i - W''(q(z)) U^1_i \right) +  \sigma_i q'(z) H_{ij}  -   \sum_{p} m_i^{p} \Lambda^{p,-1}(z,x,t) q'(z), 
\end{equation*}
and
\begin{equation*}
\begin{split}
 \left( m^p_i +  m^p_j \right) q'(z) \Lambda^{p,-1}(z,x,t) &=  m^p_i \sigma_i  \left( \partial^2_{zz} U^{1}_i - W''(q(z)) U^1_i \right) +   m^p_i \sigma_i q'(z) H_{ij}+  \\
                                                     &+  m^p_j \sigma_j  \left( \partial^2_{zz} U^{1}_j - W''(q(z)) U^1_j \right) -   m^p_j \sigma_j q'(z) H_{ij},
\end{split}
                                                     \label{lambda_AC_additive_mobility}
\end{equation*}
{where the minus sign before the term  $ m^p_j \sigma_j q'(z) H_{ij}$ comes from  $ \partial_z U^0_j = \partial_z (1-q(z)) = -q'$}.

In particular, {multiplying the last equation by $ \frac{m_i^p}{ m_i^p + m_j^p}$ and summing over $p$, we find that}

\begin{equation*}
\begin{split}
 \sum_{p} \left[   m^p_i \Lambda^{p,-1}(z,x,t)  q'(z)  \right] &=  \left(\sum_p \frac{m_i^p}{m_j^p} m_{ij}^p \right) \left[ \sigma_i  \left( \partial^2_{zz} U^{1}_i - W''(q(z)) U^1_i \right) +    \sigma_i q'(z) H_{ij} \right] \\
                                                     &+   \left(\sum_p m_{ij}^p\right)  \left[ \sigma_j  \left( \partial^2_{zz} U^{1}_j - W''(q(z)) U^1_j \right) -   \sigma_j q'(z) H_{ij} \right].
\end{split}                                                 
\end{equation*}
{where we have used the equality $ m_{ij}^p =  (\frac{1}{m_i^p} + \frac{1}{m_j^p})^{-1} = \frac{m^p_i m^p_j}{m^p_i  + m^p_j } $}.\\
{Then, by injecting this expression into the first equation, we see that}
\begin{eqnarray*} 
V^0_{ij} q'(z) &=&  \left[ \left( m_i^* -  \sum_p \frac{m_i^p}{m_j^p} m_{ij}^p  \right) \sigma_i +  \left( \sum_p m_{ij}^p \right)  \sigma_j \right] q'(z) H_{ij}  \\
&~&  \quad \quad  +  \left( m_i^* -  \sum_p \frac{m_i^p}{m_j^p} m_{ij}^p \right) \sigma_i \left( \partial^2_{zz} U^{1}_i - W''(q(z)) U^1_i \right) \\
&~& \quad  \quad  - \left( \sum_{p} m^p_{i,j} \right)\sigma_j \left( \partial^2_{zz} U^{1}_j - W''(q(z)) U^1_j \right).
\end{eqnarray*} 
Moreover, remarking that 
$${ \left( m_i^* -  \sum_p \frac{m_i^p}{m_j^p} m_{ij}^p  \right) =   \sum_p m_i^p -  \sum_p\frac{(m_i^p)^2}{m_i^p + m_j^p}  =    \sum_p   \frac{m^p_i m^p_j}{m^p_i  + m^p_j } =  \sum_p m_{ij}^p = m_{ij}} \quad{\text{ and }\quad \sigma_i + \sigma_j = \sigma_{i,j}},$$
we deduce that  $U^{1}_i$, $U^{1}_j$, and $V_{ij}^0$ satisfy  
\begin{eqnarray*} 
V_{ij}^0 q'(z) =  \sigma_{i,j} m_{ij} q'(z) H_{ij}  +  m_{ij} \sigma_i \left( \partial^2_{zz} U^{1}_i - W''(q(z)) U^1_i \right)  - m_{ij} \sigma_j \left( \partial^2_{zz} U^{1}_j - W''(q(z)) U^1_j \right).
\end{eqnarray*}

Multiplying this equation by $q'$ and integrating over $\R$ leads to the interface evolution
$$     V^0_{ij} =   m_{ij} \sigma_{i,j}  H_{ij},$$
as 
$${\int_{\R} \left( \partial_{zz}^2 U(z) - W''(z)U(z) \right) q'(z) dz =  \int_{\R} U(z) \left( q''(z) - W'(q(z)) \right)' dz = 0.}$$

{Moreover, as $U = \sigma_i U^{1}_i - \sigma_j U^{1}_j$ satisfies the equation $ \partial^2_{zz} U - W''(q) U = 0 $ 
and the boundary conditions $ \lim_{z \to \pm \infty} U = 0$, we deduce  that $ \sigma_i U^{1}_i - \sigma_j U^{1}_j = 0$. 
It follows from the partition constraint $ \sum_{k=1}^{N} U^{1}_k = U^{1}_i +  U^{1}_j = 0$ that $ U^{1}_i = U^1_j = 0$.} \\

Moreover, we have
\begin{equation*}
 \left( m^p_i +  m^p_j \right) q'(z) \Lambda^{p,-1}(z,x,t) =    m^p_i \sigma_i  q'(z) H_{ij}(x) - m^p_j \sigma_j  q'(z) H_{ij}(x),
\end{equation*}
which shows that 
$$\Lambda^{p,-1} = \left( \frac{m^p_{i,j}}{m^p_i} \sigma_i - \frac{m^p_{i,j}}{m^p_j} \sigma_j \right) H_{i,j},$$
and Proposition~\ref{claim:allencahn_mob_additive2} ensues.

\section{Numerical scheme and implementation}

In this section we introduce a Fourier spectral splitting scheme~\cite{Chen_fourier}
to approximate the solutions to the Allen-Cahn system
\begin{equation*} 
 \partial_t u^{\varepsilon}_k =   m^{*}_k \left[ \sigma_k \left( \Delta u^{\varepsilon}_k - \frac{1}{\varepsilon^2} W'(u^{\varepsilon}_k) \right) +  \lambda^{\varepsilon}_k \sqrt{2 W(u_k)} \right],\quad k \in \{1,\dots,N\},
\end{equation*}
where  $m^{*}_k = \sum_{p=1}^{P} m^p_{k}$ and   
 $$ \lambda^{\varepsilon}_k = \frac{1}{m^{*}_k} \sum_{p=1}^{P} m_{k}^{p}~\lambda^{p,\varepsilon}, \quad \text{ with } \quad  \lambda^{p,\varepsilon} = - \left( \frac{ \sum_{k=1}^{N} m_k^p  \sigma_k \left( \Delta u^{\varepsilon}_k - \frac{1}{\varepsilon^2} W'(u^{\varepsilon}_k)  \right)}{\sum_{k=1}^{N} m_k^p \sqrt{2 W(u_k)}} \right).$$
The solutions to the system are approximated numerically on a square box $Q = [0,L_1]\times \cdots \times [0,L_d]$ 
with periodic boundary conditions. \\

We recall that the Fourier $K$-approximation of a function $u$ defined in a box 
$Q$ is given by
$$u^{\boldsymbol K}(x) = \sum_{{\boldsymbol k}\in K_d  } c_{\boldsymbol k} e^{2i\pi{\boldsymbol \xi}_k\cdot x},$$
where  $K_d =  [ -\frac{K_1}{2},\frac{K_1}{2}-1 ]\times [ -\frac{K_2}{2},\frac{K_2}{2}-1] \cdots \times   [ -\frac{K_d}{2},\frac{K_d}{2}-1] $,   ${\boldsymbol k} = (k_1,\dots,k_d)$, and ${\boldsymbol \xi_k} = (k_1/L_1,\dots,k_d/L_d)$. In this formula, the $c_{\boldsymbol k}$'s denote the $K^d$ first discrete Fourier coefficients of $u$. 
The inverse discrete Fourier transform leads to 
$u^{K}_{\boldsymbol k} =   \textrm{IFFT}[c_{\boldsymbol k}]$ 
where $u^{K}_{\boldsymbol k}$ denotes the value of $u$ at the points 
$x_{\boldsymbol k} = (k_1 h_1, \dots, k_d h_d)$ and where  $h_{i} = L_{i}/N_{i}$ for $i\in\{1,\dots,d\}$. Conversely,
$c_{\boldsymbol k}$ can be computed as the discrete Fourier transform of $u^K_{\boldsymbol k},$ {i.e.}, $c_{\boldsymbol k} =
\textrm{FFT}[u^K_{\boldsymbol k}].$

\subsection{Definition of the scheme}\label{sec:def-scheme}
 
Given a time discretisation parameter $\delta_t > 0$,
we construct a sequence $({\boldsymbol u}^n)_{n \geq 0}$ 
of approximations of ${\boldsymbol u}$ at the times $n \delta_t$,
{by adapting the splitting discretization schemes proposed in the previous works} \cite{MR3961089,BRETIN2018324}. 
More precisely, we iteratively
 
\begin{itemize}
 \item minimize the Cahn-Hilliard energy without the constraint
 $\sum_{k=1}^N u_k^n = 1$. 
 \item compute the contribution of the Lagrange multipliers $\lambda_k^{\varepsilon}$
 and update the values of $u_k^n$. 
\end{itemize}
This approach provides a simple scheme, and our numerical experiments (see Section 4) together with Proposition~\ref{prop_scheme}
indicate that it is effective, stable, and that it conserves the partition constraint 
in the sense that 
$$\sum_{k=1}^{N} u_k^{n+1} = \sum_{k=1}^{N} u_k^{0}, \quad \forall n \in \N.$$
Let us now give more details about our scheme.
 
\begin{itemize}
\item[\textbf{Step $1$:}] 
{\bf Solving the decoupled Allen-Cahn system (i.e., without the partition constraint): } \\
~\\
Let ${\boldsymbol u}^{n+1/2}$ denote 
an approximation of ${\boldsymbol v}(\delta_t)$, 
where ${\boldsymbol v} = (v_1,\dots, v_N)$ is the solution with periodic boundary conditions 
on $\partial Q$ to:
$$
\begin{cases}
   \partial_t v_k(x,t) &= m_k^{*} \sigma_k \left[\Delta v_k (x,t) - \frac{1}{\varepsilon^2} W'(v_k(x,t))\right], \quad (x,t) \in Q \times [0,\delta_t] \text{,}\\
  {\boldsymbol v}(x,0) &= {\boldsymbol u}^{n}(x), \quad x \in Q.
\end{cases}
$$

Here, our motivation is to introduce a stable scheme in the sense  
that the associated Cahn-Hilliard energy decreases with the iterations.  
A totally implicit scheme would require the resolution of a nonlinear system
at each iteration, which in practice would prove costly and not very accurate.
Rather, we opt for a semi-implicit scheme in which the non linear term $W'(v_k)$
is integrated explicitly.
More precisely, we consider the  scheme 
$$ 
\left( I_d - \delta_t m_k^{*} \sigma_k \left( \Delta  -  \alpha/\varepsilon^2 I_d \right)   \right) u_k^{n+1/2} =     u_k^{n}  - \frac{\delta_t m_k^{*} \sigma_k}{\varepsilon^2} 
\left( W'(u_k^n) - \alpha u_k^{n}\right),
$$
where $\alpha$ is a positive stabilization parameter, chosen sufficiently large 
to ensure the stability of the scheme. 
Indeed, it is known that the Cahn-Hilliard energy decreases unconditionally 
as soon as the explicit part, {i.e.} $s \to W'(s) - \alpha s$, is the derivative 
of a concave function~\cite{EDTK_Eyre,ShenWWW12}. 
This is the case for the potential $W(s) = \frac{1}{2} s^2 (1-s)^2$, as soon as $\alpha>2$.
We also note that even when $\alpha = 0$, the semi-implicit
scheme is stable under the classical condition 
$\delta_t \leq \frac{C}{\varepsilon^2}$, where $C = \sum_{s \in [0,1]} |W^{\prime\prime}(s)|.$
Further, as the fields $u_k^n$ are required to satisfy periodic boundary conditions
on $\partial Q$, the action of the inverse operator
$\left( I_d - m_k \delta_k\delta_t \left( \Delta  -  \alpha/\varepsilon^2 I_d \right) \right)^{-1}$ 
is easily computed in Fourier space \cite{Chen_fourier} using the Fast Fourier Transform.
Remark that this strategy can also be generalized to anisotropic flows~\cite{MR2929124}.

\item[\textbf{Step $2$:}] {\bf Explicit projection onto the partition constraint $\sum u_k = 1$ }. \\
~\\
The advantage of an implicit treatment of the Lagrange multiplier $\lambda_k^{\varepsilon}$ is not 
significant enough considering the complexity and the algorithmic cost of this approach. 
We rather prefer an explicit approach for which we will prove that the processing is exact in the sense that
$\sum_{k=1}^{N} u_k^{n+1} =    \sum_{k=1}^{N} u_k^{n},  \quad \forall n \in \N$.
More precisely,  for all $k\in \{ 1,\dots, N \}$, we define $u^{n+1}_k$ by
$$ u_k^{n+1} =  u^{n+1/2}_k + \delta_t  m_k^{*}\lambda^{n+1/2}_{k} \sqrt{2 W(u^{n+1/2}_k)}$$
where
$$\lambda^{n+1/2}_{k} =  \frac{1}{m_k^{*}} \sum_{p=1}^{P} m_k^p \lambda^{p,n+1/2}, \quad \text{ and }\quad {\lambda^{p,n+1/2} =  -\frac{\sum_{i=1}^{N} m_i^p \alpha^{n+1/2}_i }{\sum_{i=1}^{N} m_i^p \sqrt{2W(u_i^{n+1/2})}}}
.$$ 
                     
Here  $\alpha^{n+1/2}_i$ is a semi-implicit approximation of
$\sigma_i \left[\Delta u_i (x,t) - \frac{1}{\varepsilon^2} W'(u_i(x,t))\right]$ at time $t_{n+1/2}$  defined by 
$$\alpha_i^{n+1/2} = \frac{ u_i^{n+1/2} - u_i^{n} }{ \delta_t m_i^*}.$$
\end{itemize}

\begin{rem} {Notice that we can always assume that $\sum_{i} m_i^p > 0$,  as otherwise $m_{i,j}^{p} = 0$ for all $i,j$ and there is no
contribution of the $p$-th mobility. Moreover, }the above definitions of $\alpha_i^{n+1/2}$ and $\lambda^{n+1/2}_{k}$  
only make sense when the $m_k^{*}$'s  or the sum $\sum_{i=1}^{N} m_i^p \sqrt{2W(u_i^{n+1/2})}$ 
do not vanish.  In practice (see the code in Table \ref{fig:matlab_code}), to overcome this difficulty and avoid any division by zero, {it is more convenient to work with
$ \tilde{\lambda}^{n+1/2}_{k} = m_k^{*} \lambda^{n+1/2}_{k}$} and 
to use the following regularized version of the scheme : 

$$ u_k^{n+1} =  u^{n+1/2}_k + \delta_t  \tilde{\lambda}^{n+1/2}_{k}\left( \sqrt{2 W(u^{n+1/2}_k)} + \beta \right), \quad \tilde{\lambda}^{n+1/2}_{k} =  \sum_{p=1}^{P} m_k^p \tilde{\lambda}^{p,n+1/2},$$
and

$$
\tilde{\lambda}^{p,n+1/2} =  -\frac{\sum_{i=1}^{N} m_i^p {\tilde{\alpha}}^{n+1/2}_i }{\sum_{i=1}^{N} m_i^p (\sqrt{2W(u_i^{n+1/2})} + \beta ) }, \text{ and } \tilde{\alpha}_i^{n+1/2} = \frac{ u_i^{n+1/2} - u_i^{n} }{ \delta_t \text{max}\{m_i^*,\beta\}},
$$

where $\beta \simeq 2.22 10^{-16}$  is the machine precision.  
\end{rem}

{The next proposition shows that our scheme conserves the partition constraint and conserves  each phase whose mobilities at each of its interfaces are zero.}
\begin{prop} \label{prop_scheme} With the above notations:
\begin{itemize}
 \item[(1)] Assume that $m_i^{*}> \beta$ for all $i \in \{1, \dots, N\}$ and  $\sum_{i} m_i^{p}>0$ for all $p \in \{1, \dots, P \}$, 
 then the previous scheme preserves the partition constraint, i.e. 
 $$\sum_{k=1}^{N} u_k^{n+1}  = \sum_{k=1}^{N} u_k^{n}.$$
 \item[(2)] Let $i\in  \{1,\dots, N\}$. Assume that $m_{ij} = 0$ for all $j \in \{1,\dots, N\}$, $j\not=i$, then   
 $$u_i^{n+1}= u_i^{n}.$$
\end{itemize}
\end{prop} 

{\bf Proof of $(1)$ : } 
As 
$$ u_k^{n+1} =  u^{n+1/2}_k -   \sum_{p=1}^{P} m_k^p  \left[ \frac{\sum_{i=1}^{N} m_i^p  \left(u_i^{n+1/2} - u_i^{n} \right)/ \text{max}\{ m_i^*, \beta \}}{\sum_{i=1}^{N} m_i^p (\sqrt{2W(u_i^{n+1/2})} + \beta)} \right]   (\sqrt{2 W(u^{n+1/2}_k)} + \beta)$$
it follows that 
\begin{eqnarray*}
 \sum_{k=1}^{N} u_k^{n+1} &=& \sum_{k=1}^{N} u^{n+1/2}_k  -  \sum_{p=1}^{P} \left[ \sum_{i=1}^{N} m_i^p  \left(u_i^{n+1/2} - u_i^{n} \right)/m_i^* \right] \frac{ \sum_{k=1}^{N} m_k^p  (\sqrt{2 W(u^{n+1/2}_k)} + \beta)  }{ \sum_{i=1}^{N} m_i^p  ( \sqrt{2 W(u^{n+1/2}_i)} + \beta) } \\
                          &=& \sum_{k=1}^{N} u^{n+1/2}_k  -  \sum_{i=1}^{N} (u_i^{n+1/2} - u_i^{n}) \sum_{p=1}^{P}  m_i^p/m_i^*  = \sum_{k=1}^{N} u^{n}_k.
\end{eqnarray*}
{ as  $\sum_{p=1}^{P}  m_i^p = m_i^*$.} 

{\bf Proof of $(2)$ : }  {By definition, $\sum_{p} m_{ij}^{p} = m_{ij}=0$. All $m_{ij}^p$'s being nonnegative, we deduce that $m_{ij}^p = 0$, $p \in \{1,\dots, P\}$, $j \in \{1,\dots, N\}\setminus \{i\}$. Moreover, as $m_{ij}^{p}$ is harmonically additive, i.e. $1/m_{ij}^{p} = (1/m_i^p + 1/m_j^p)$ (using the convention $\frac 1{0^+}=+\infty$ and $\frac 1{+\infty}=0^+$), if follows  that, for every $q\in\{1,\dots,P\}$, either $m_i^{q}=0$, or $m_i^{q}>0$ and $m_j^{q}=0$ for all  $j \in \{1,\dots, N\}\setminus \{i\}$. In the latter case, $m_{kj}^{q}=0$ for every $k\not=j$, so the ${q}$-th term is useless in the decomposition of all mobility coefficients $m_{kj}$, $k\not=j$, and can be removed. Using the same argument for every $q\in\{1,\dots,P\}$, and discarding the trivial situation where {\it all} mobility coefficients $m_{kj}$ are zero, we finally obtain that necessarily  $m_i^{q}=0$, $\forall q\in\{1,\dots,P\}$ (using for simplicity the same notation $P$ for the new number of elements in the decomposition). Now,  the first step of our scheme yields $u_i^{n+1/2} = u_i^{n}$ as $m_i^{*} =  \sum_{p} m_i^{p}=0$}, and the second step gives 
$ \tilde{\lambda}_{i}^{n+1/2} = \sum_{p=1}^{P} m_{i}^{p} \tilde{\lambda}_{i}^{p,{n+1/2}} = 0$ and $u_i^{n+1} = u_i^{n+1/2} = u_i^{n}$,
as $\tilde{\lambda}^{p,n+1/2}$ is bounded and $m_{i}^{p} = 0$ for all $p \in \{1,\dots, P \}$.

%
\begin{rem}
 {The above argument brings the following natural question: given a collection of coefficients $\{m_{k\ell}^p, 1\leq k<\ell<N\}$, $p\in\{1,\dots,P\}$, how can we compute the coefficients $m_i^p$?
 In the case of positive coefficients, the following formula can be used, whose proof is straightfoward:
 $$\frac 1{m_i^p} = \frac 1 2 \left(\frac 1{m_{ij}^{p}} +  \frac 1{m_{ik}^{p}}  -  \frac 1{m_{kj}^{p}} \right).$$
In practice, in particular in a numerical code, this formula can be extended to general nonnegative coefficients by simply replacing each null coefficient with the machine precision.}
\end{rem}

\subsection{{\bf Matlab} code } 
We present in Table~\ref{fig:matlab_code} an example of Matlab code with less than $50$ lines which implements our splitting scheme in the case of $N=3$ phases in dimension $d=2$. Here is a short description of a few lines of the code:
\begin{itemize}
 \item Lines $5$ to $8$ correspond to the initialization of the phase $u_k$ for $k=1, 2, 3$.
 \item  Lines $11$ to $18$ implement  a canonical decomposition of the mobility coefficients $m_{ij}$: 
 \begin{eqnarray*}
  (m_{12},m_{13},m_{23}) &=& (m^1_{12},m^1_{13},m^1_{23}) +  (m^2_{12},m^2_{13},m^2_{23}) + (m^3_{12},m^3_{13},m^3_{23}) \\
                         &=&  (m_{12},0,0) +  (0,m_{13},0) + (0,0,m_{23})
 \end{eqnarray*}
 and provide the associated coefficients $m_i^p$.
 \item Lines  $20$ to $24$ provide the operators necessary in Step 1 of the scheme for the numerical resolution of the equation 
 $$ \left( I_d - \delta_t m_k^{*} \sigma_k \left( \Delta  -  \alpha/\varepsilon^2 I_d \right)   \right) u_k^{n+1/2} =     u_k^{n}  - \frac{\delta_t m_k^{*} \sigma_k}{\varepsilon^2} \left( W'(u_k^n) - \alpha u_k^{n}\right),$$
Two operators are introduced:
$$ OP(u,m) = u -  \frac{\delta_t m \sigma}{\varepsilon^2} \left( W'(u) - \alpha u\right)\quad \text{ and }\quad  OL(u,m) =  \left( I_d - \delta_t m \sigma \left( \Delta  -  \alpha/\varepsilon^2 I_d \right)   \right)^{-1} u$$
 \item Lines $33$-$39$ correspond to the computation of each Lagrange multiplier $\lambda^{p,n+1/2}$.
\end{itemize}
 
\begin{table}[htbp]
\centering
\centering
\begin{lstlisting}
%%%%%%%%% Resolution parameters 
N = 2^8; epsilon = 1/N; dt = 10/N^2; L = 1; T = 1;

%%%%%%%%%% initial condition %%%%%%%%%%%%%%
x = linspace(-L/2,L/2,N); [Y,X] = meshgrid(x,x); R = 0.3;
d1 = max(sqrt((X).^2 + (Y + 0.1).^2) - R,Y - 0.05*cos(12*pi*X)); U1 = 1/2 - 1/2*(tanh(d1/epsilon/2));
d2 = max(sqrt((X).^2 + (Y - 0.1).^2) - R, - Y + 0.05*cos(12*pi*X)); U2 = (1/2 - 1/2*(tanh(d2/epsilon/2)));
U3 = 1 - (U1 + U2);

%%%%%%%%%%%%% surface and mobility coefficients %%%%%%%%%%
sigma12 = 1; sigma13 = 1; sigma23 = 1; m12 = 1; m13 = 0; m23 =0;

%%%%%%%%%%%%% coefficients m_{i,j}^p and  m_{i}^p %%%%%%%%%%%% 
sigma1 = (sigma12 + sigma13 - sigma23)/2; sigma2 = (sigma12 + sigma23 - sigma13)/2; sigma3 = (sigma23 + sigma13 - sigma12)/2;
m12_1 = m12; m13_1 = 0; m23_1 = 0; m1_1 =  2*m12; m2_1 =  2*m12; m3_1 = 0;
m12_2 = 0;   m13_2 = m13; m23_2 = 0; m1_2 =  2*m13; m2_2 =  0; m3_2 = 2*m13;
m12_3 = 0;   m13_3 = 0; m23_3 = m23; m1_3 =  0 ; m2_3 =  2*m23; m3_3 = 2*m23;
m1 = m1_1 + m1_2 + m1_3; m2 = m2_1 + m2_2 + m2_3; m3 = m3_1 + m3_2 + m3_3;
%%%%%%%%%%%%%%%%%%%%%%%% Diffusion and reaction operators 
k = [0:N/2,-N/2+1:-1]; [K1,K2] = meshgrid(k,k); Delta = (K1.^2 + K2.^2);
sqrtWU  =  @(U) abs(U.*(1-U)); potentiel_prim = @(U) U.*(1-U).*(1 - 2*U);
alpha = 2;
OP = @(U,dt,epsilon,sigma,m,alpha) U - dt/epsilon^2*sigma*m*(potentiel_prim(U) - alpha*U); 
OL = @(U,dt,epsilon,sigma,m,alpha) ifft2((1./(1 + dt*m*sigma*(4*pi^2*Delta + alpha/epsilon^2))).*fft2(U)); 

%%%%%%%%%%%%%%%%%%%%%%%%% Computation of the solution %%%%%%%%%%%%%%%%
for n=1:T/dt, 
 %%%%%%%%%  Step 1 % Cahn Hilliard flow
 U1_p  = OL(OP(U1,dt,epsilon,sigma1,m1,alpha),dt,epsilon,sigma1,m1,alpha);
 U2_p  = OL(OP(U2,dt,epsilon,sigma2,m2,alpha),dt,epsilon,sigma2,m2,alpha);
 U3_p  = OL(OP(U3,dt,epsilon,sigma3,m3,alpha),dt,epsilon,sigma3,m3,alpha);
 %%%%%%%%  Step 2 % Lagrange multiplier Lambda 
alpha1 = (U1_p - U1)/(dt*max(m1,eps)); alpha2 = (U2_p - U2)/(dt*max(m2,eps));  alpha3 = (U3_p - U3)/(dt*max(m3,eps));
if (m1_1 +  m2_1 + m3_1>0), lambda_p1 = - ( m1_1*alpha1 + m2_1*alpha2 + m3_1*alpha3)./(m1_1*(sqrtWU(U1_p)+eps) + m2_1*(sqrtWU(U2_p)+eps) + m3_1*(sqrtWU(U3_p)+eps));
else lambda_p1 = 0;    end
if ( m1_2 + m2_2 + m3_2>0), lambda_p2 = - ( m1_2*alpha1 + m2_2*alpha2 + m3_2*alpha3)./(m1_2*(sqrtWU(U1_p)+eps) + m2_2*(sqrtWU(U2_p)+eps) + m3_2*(sqrtWU(U3_p)+eps));
else lambda_p2 = 0;    end
if ( m1_3 + m2_3 + m3_3>0), lambda_p3 = - ( m1_3*alpha1 + m2_3*alpha2 + m3_3*alpha3)./(m1_3*(sqrtWU(U1_p)+eps) + m2_3*(sqrtWU(U2_p)+eps) + m3_3*(sqrtWU(U3_p)+eps));
else lambda_p3 = 0;  end

U1 = U1_p + dt*(m1_1*lambda_p1 + m1_2*lambda_p2 + m1_3*lambda_p3).*(sqrtWU(U1_p)+eps);
U2 = U2_p + dt*(m2_1*lambda_p1 + m2_2*lambda_p2 + m2_3*lambda_p3).*(sqrtWU(U2_p)+eps);
U3 = U3_p + dt*(m3_1*lambda_p1 + m3_2*lambda_p2 + m3_3*lambda_p3).*(sqrtWU(U3_p)+eps);

 if (mod(n,10)==1)
       imagesc(U3 + 2*U2)
       pause(0.01);
 end
end

 














\end{lstlisting}
	 \caption{Example of a {\bf Matlab} implementation of the scheme described in Section~\ref{sec:def-scheme} in the case 
	 of $N=3$ phases in dimension $2$.}
\label{fig:matlab_code}
\end{table}

\section{Numerical experiments and validation  }

In this section, we report numerical experiments in dimensions $d=2$ and $d=3$, 
with $N=3$ or $N=4$ phases. 
In each case, the computational domain $Q$ is a unit cube $[-0.5,0.5]^d$ discretized  in each direction with $K = 2^8$ nodes in 2D and $K = 2^7$ in 3D. 
We use the classical double-well potential 
$W(s) = \frac{1}{2} s^2 (1-s)^2$. 

\subsection{Validation of the consistency of our approach}

This first test illustrates the consistency of the numerical scheme in the case of 
$N=3$ phases. 
We consider the evolution of two circles by the 
flow~\eqref{eqn:allencahn_mob_additive2} 
associated to the following surface tensions and mobilities:
$$ (\sigma_{12},\sigma_{13},\sigma_{23}) = (1,1,1) 
\quad\text{ and } \quad (m_{12},m_{13},m_{23}) = (1,1,\frac 1 4).$$
{Notice that this set of mobilities is not harmonically additive as} 
$$ \frac{1}{m_1}  =  \frac{1}{2} ( \frac{1}{m_{13}} + \frac{1}{m_{12}} - \frac{1}{m_{23}}) = \frac{1}{2} (1 + 1 - 4) = -1 <0.$$ 
{We use therefore our approach with the canonical decomposition  which reads}
$$(m_{12},m_{13},m_{23}) = \sum_{p=1}^{3}(m^p_{12},m^p_{13},m^p_{23}),$$
with 
$$ (m^1_{12},m^1_{13},m^1_{23}) = (m_{12},0,0), \; (m^2_{12},m^2_{13},m^2_{23}) = (0,m_{13},0),  \; \text{ and }  \; (m^3_{12},m^3_{13},m^3_{23}) = (0,0,m_{23}),$$
where
$$ (m^1_{1},m^1_{2},m^1_{3}) = (2m_{12} ,2m_{12},0),  \; (m^2_{1},m^2_{2},m^2_{3}) = (2m_{13},0,2m_{13}),  \; \text{ and }   \;(m^3_{1},m^3_{1},m^3_{1}) = (0,2m_{23},2m_{23}).$$

Moreover, the initial sets are chosen in the following way:
\begin{itemize}
 \item the phase $u_1$  fills a circle of radius $r_1 = 0.2$ centered 
 at $ x= (-0.25,0,0)$
 \item the phase $u_2$ fills a circle of radius $r_2 = 0.2$ 
 centered  at $ x=  (0.25,0,0)$
\end{itemize}
These initial sets should evolve as circles with radius:
$$ R_1(t) = \sqrt{r_1^2 - 2 \sigma_{13}m_{13}  t} \quad\text{ and }\quad R_2(t) = \sqrt{r_2^2 - 2 \sigma_{23}m_{23}  t}.$$

\begin{figure}[htbp]
\centering
	\includegraphics[width=3.5cm]{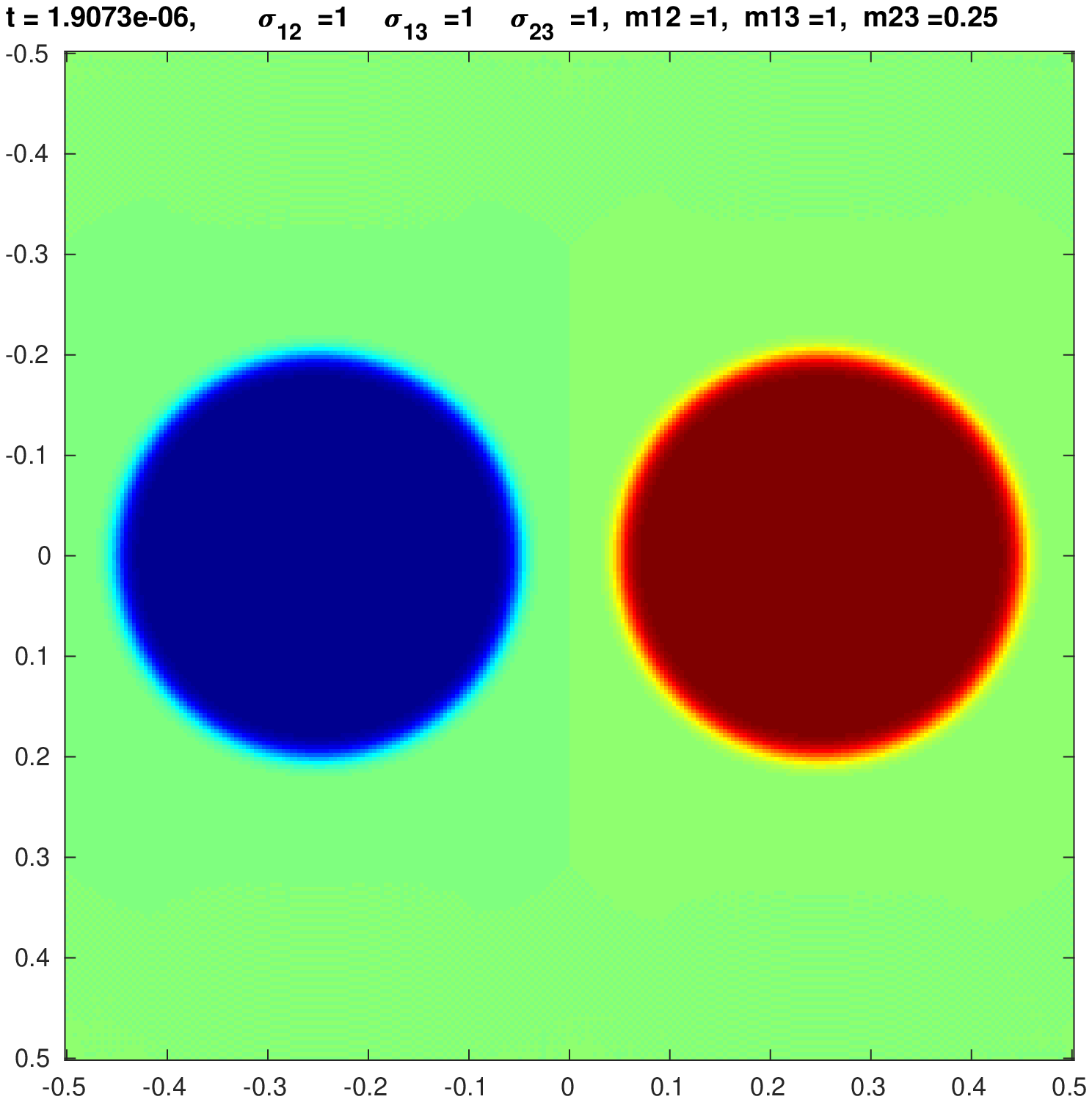}
        \includegraphics[width=3.5cm]{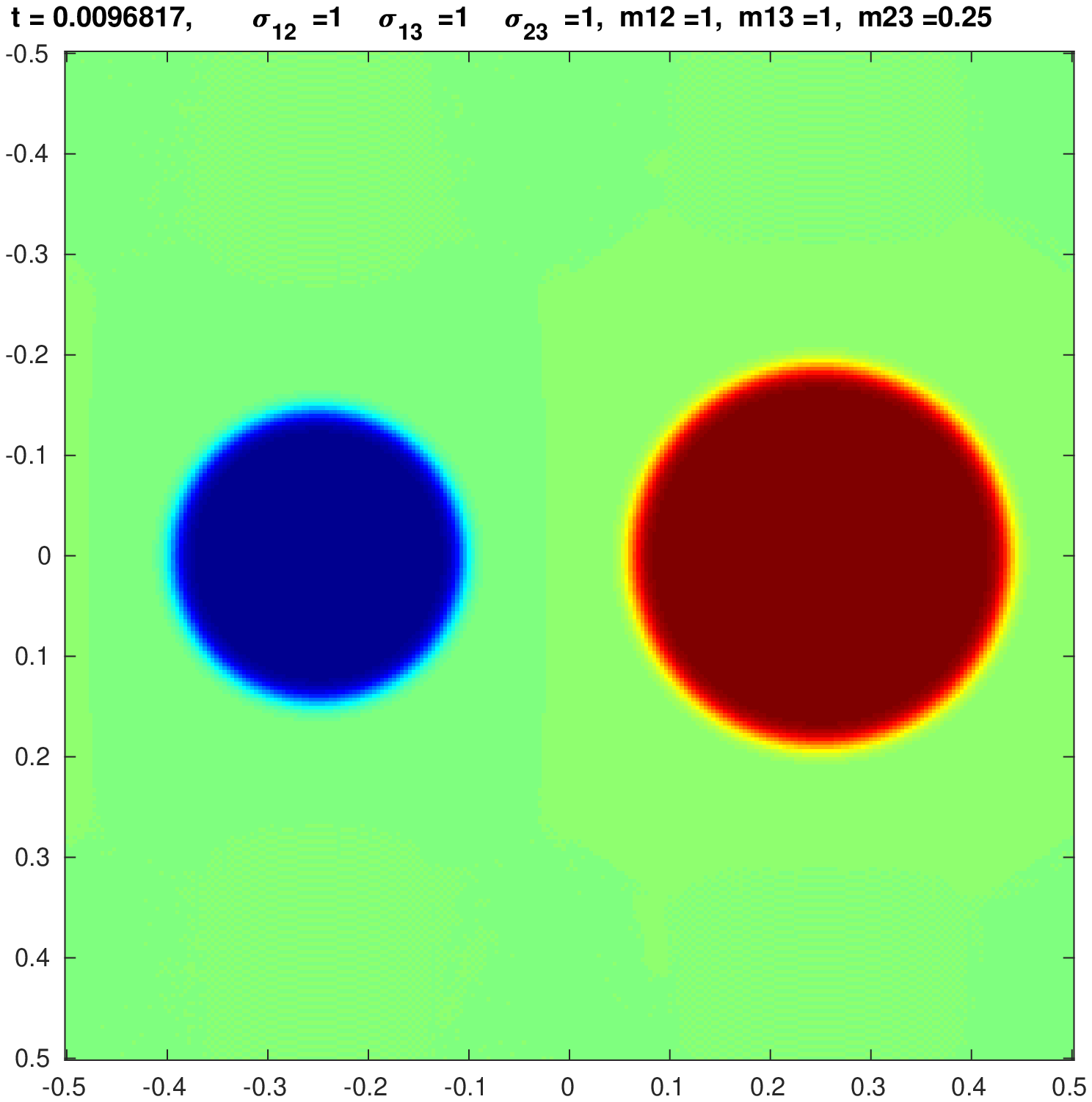}
          \includegraphics[width=3.5cm]{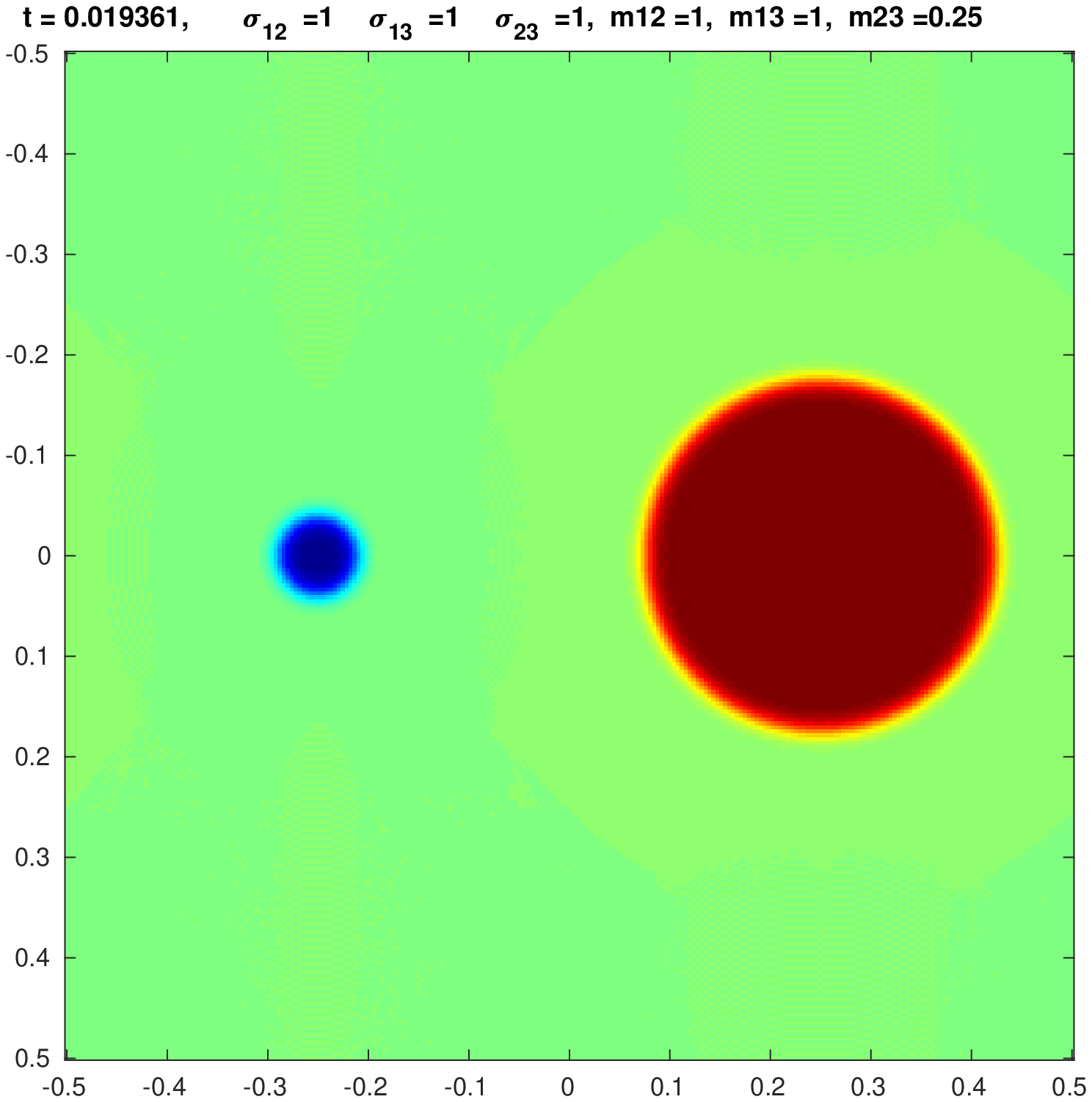}
       \includegraphics[width=3.5cm]{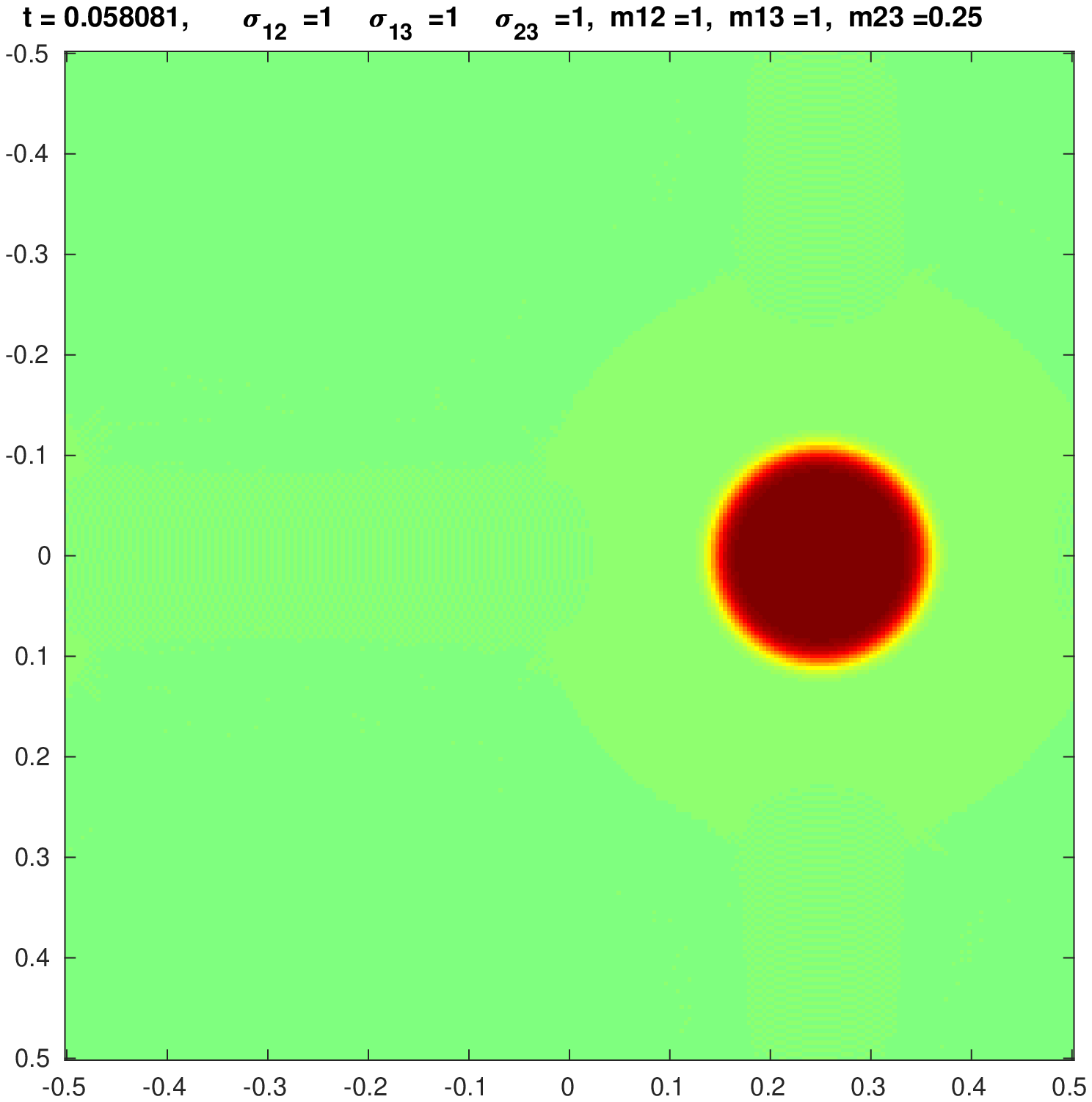}
       
\caption{Mean curvature flow of two circular phases, with $\sigma_{12} = \sigma_{13} =  \sigma_{23} = 1$ and $m_{12}= 1$, $m_{13}= 1$, $m_{23} = \frac 1 4$. 
The images show the values of the function $2u_2 + u_3$ at different times {using a colormap such that} $u_1$, $u_2$, and $u_3$ 
are represented in blue, red and green, respectively.}
\label{fig_circle1}
\end{figure}

\begin{figure}[htbp]
\centering
	\includegraphics[width=6cm]{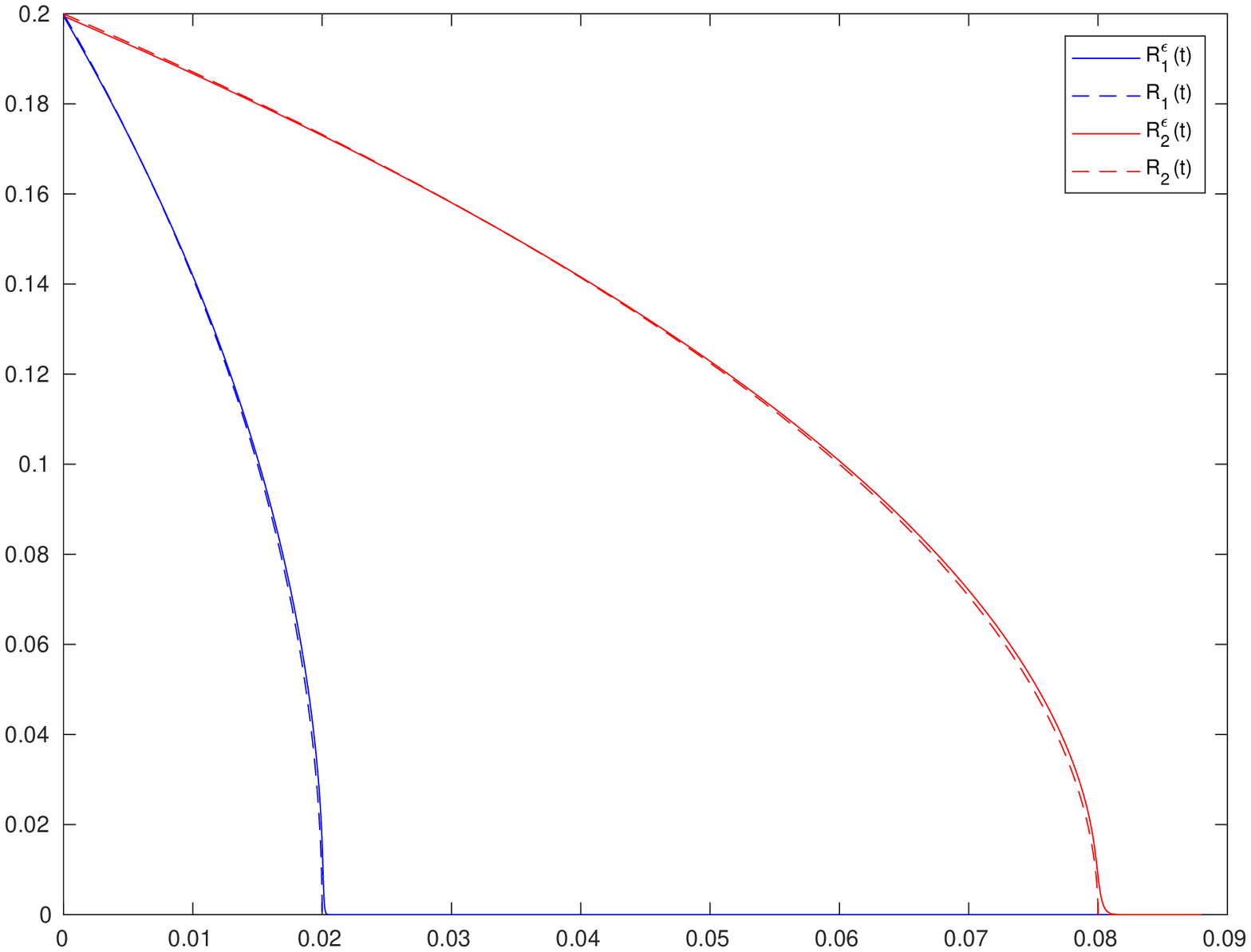}
	\includegraphics[width=6cm]{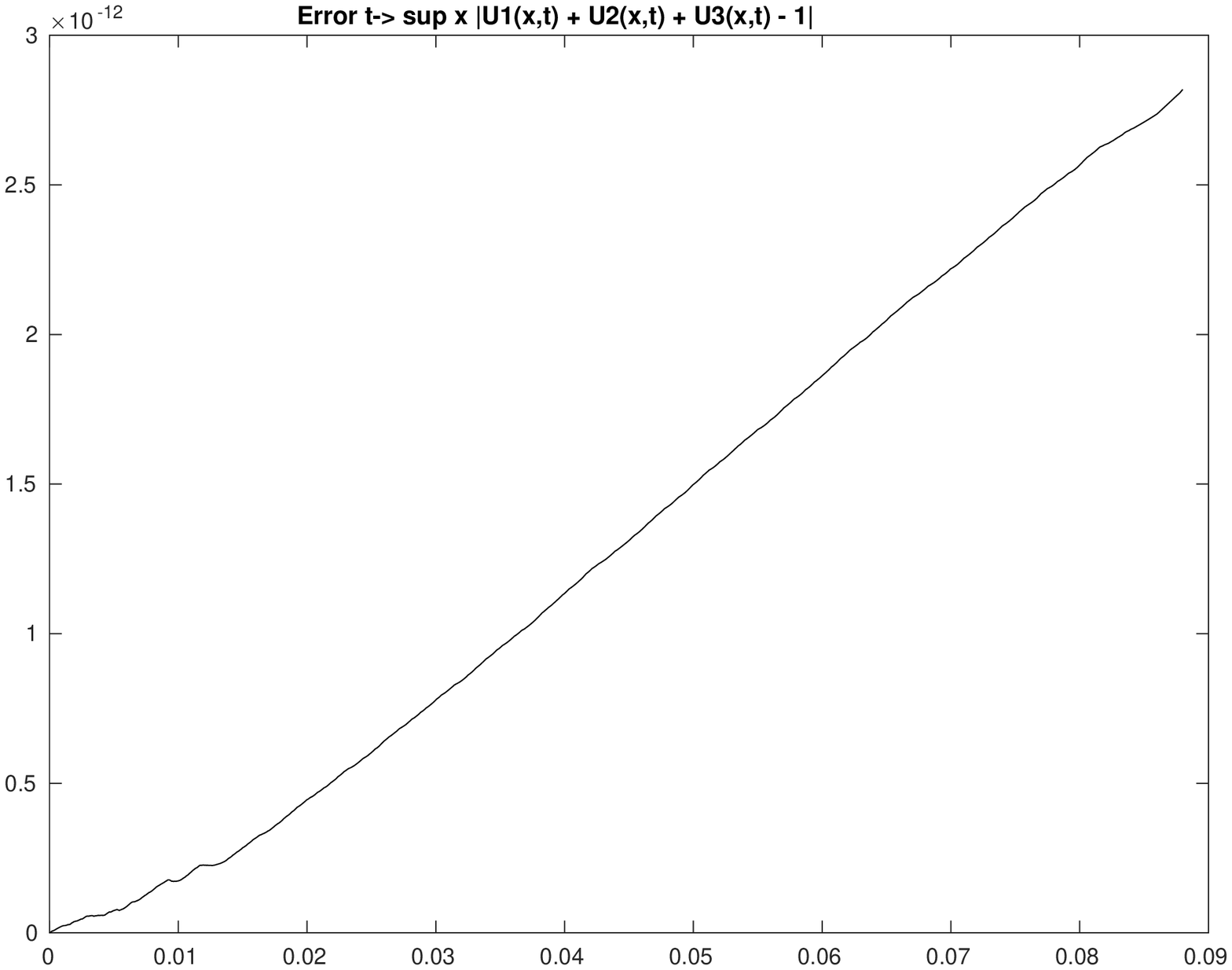}	
	
\caption{Mean curvature flow of two circular phases. Surface tensions are identical, mobilities are $(m_{12},m_{13},m_{23})=(1,1,\frac 1 4)$.
Left : Comparison of the radii $R_k^{\varepsilon}(t)$  and their theoretical values $R_i(t)$
associated with the phases $u_k$, k=1,2.
Right: Plot of $\|1 - \sum_k u_k(.,t)\|_{L^{\infty}}$.}
\label{fig_circle2}
\end{figure}

The following parameters are used for the computations:
$\varepsilon = 1.5/K$, $\delta_t = 0.25/K^2$, and  $\alpha = 0$.
Figure~\ref{fig_circle1} shows the numerical multiphase solution 
${\bf u}^{\varepsilon} = (u_1^{\varepsilon},u_2^{\varepsilon},u_3^{\varepsilon})$ 
at different times. 
The first graph in Figure~\ref{fig_circle2} shows a very good agreement
between the approximative radii $R_1^{\varepsilon}$ and $R_2^{\varepsilon}$ 
and their expected theoretical values. 
The second graph in Figure~\ref{fig_circle2} shows that 
the numerical error on the constraint $\sum_k u_k = 1$ is of the order of $10^{-12}$
in this context. \\

\subsection{Influence of the choice of a particular decomposition of mobilities}

The decomposition~\eqref{m-additive} is not unique and it is therefore legitimate 
to question its effect on the numerical approximation of the flow.
We consider here the simplest case using $N=3$ phases,   
homogeneous surface tensions $\sigma_{i,j} = 1$ and homogeneous mobility coefficients $m_{ij} = 1$. 
We then compare the numerical approximations associated with the following
decompositions of the mobilities:
\begin{itemize}
 \item {the canonical choice with $P=3$:} 
 $$ (m^1_{12},m^1_{13},m^1_{23}) = (1,0,0), \; (m^2_{12},m^2_{13},m^2_{23}) = (0,1,0)  \; \text{ and }  \; (m^3_{12},m^3_{13},m^3_{23}) = (0,0,1)$$
where
 $$ (m^1_{1},m^1_{2},m^1_{3}) = (2,2,0),  \; (m^2_{1},m^2_{2},m^2_{3}) = (2,0,2)  \; \text{ and }   \;(m^3_{1},m^3_{1},m^3_{1}) = (0,2,2).$$ 
 \item {a sparse  decomposition with $P=1$:} 
 $$(1,1,1) = (m_{12},m_{13},m_{23}) = \sum_{p=1}^{1}(m^p_{12},m^p_{13},m^p_{23}) = (m^1_{12},m^1_{13},m^1_{23}),$$
 where $(m^1_{1},m^1_{2},m^1_{3}) = (2,2,2)$. Notice that $(m_{12},m_{13},m_{23}) = (1,1,1)$ is indeed harmonically additive 
 which explains why we can use $P=1$.
\end{itemize}

\begin{figure}[htbp]
\centering
\includegraphics[width=3.5cm]{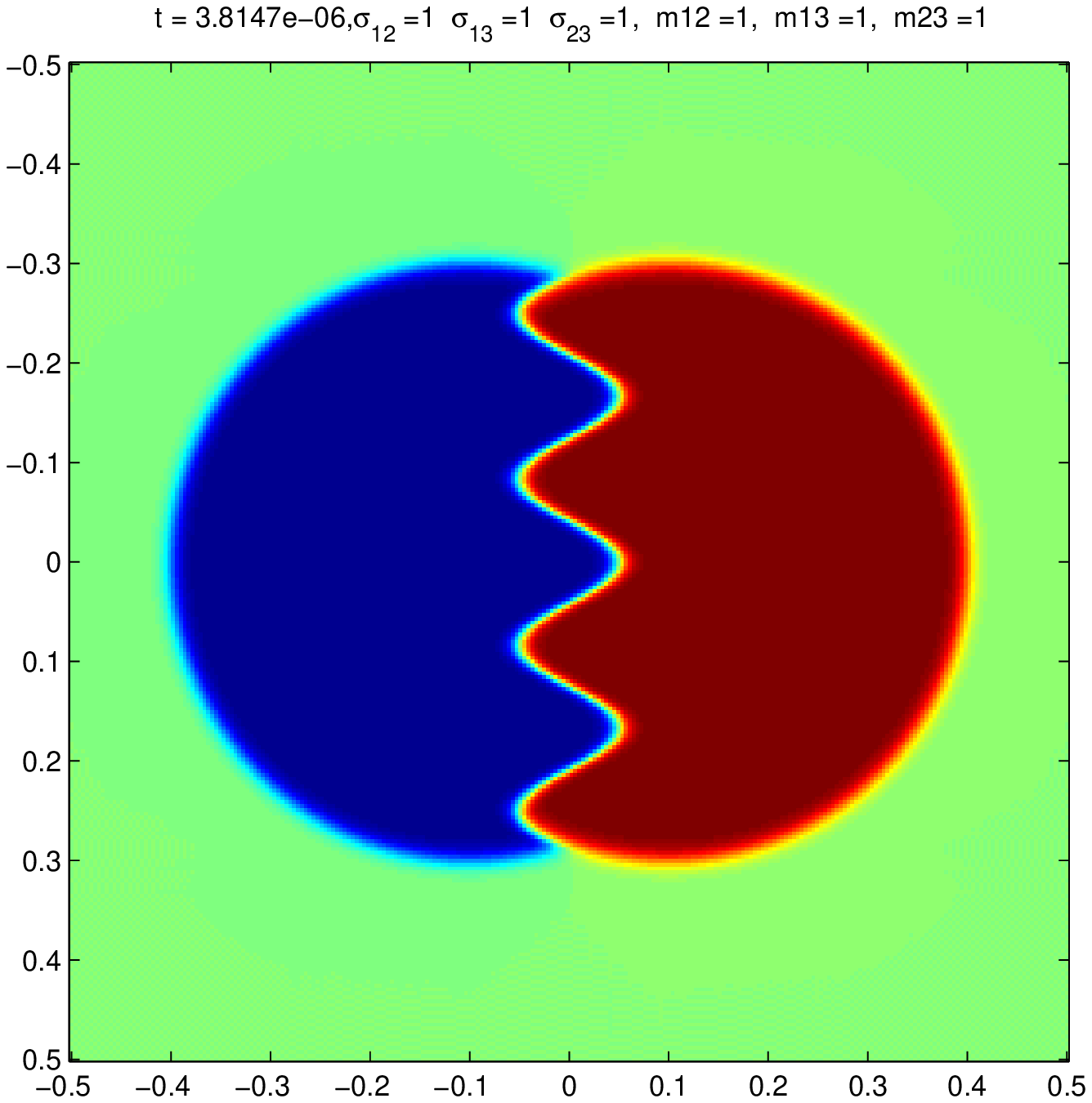}
\includegraphics[width=3.5cm]{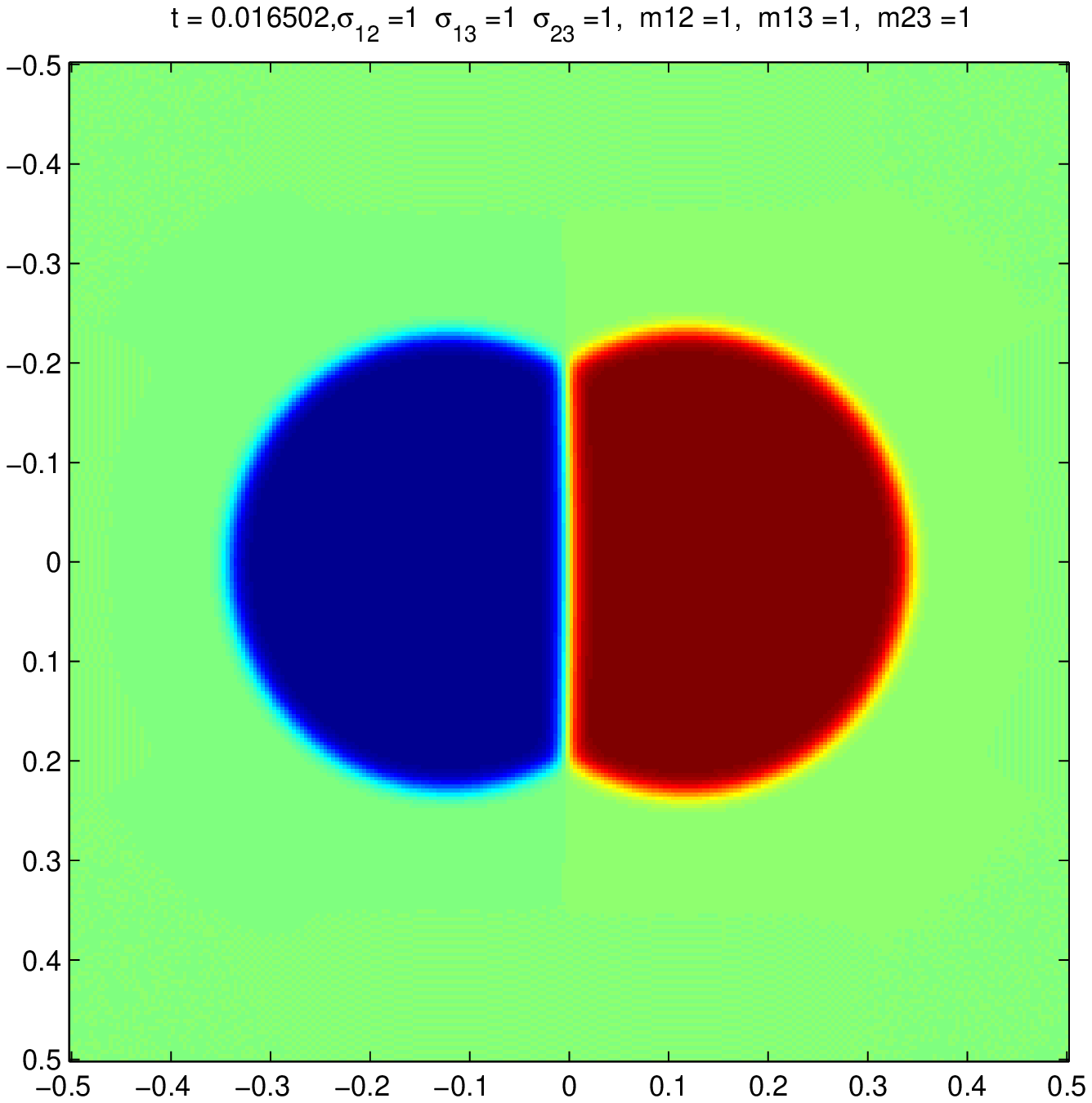}
\includegraphics[width=3.5cm]{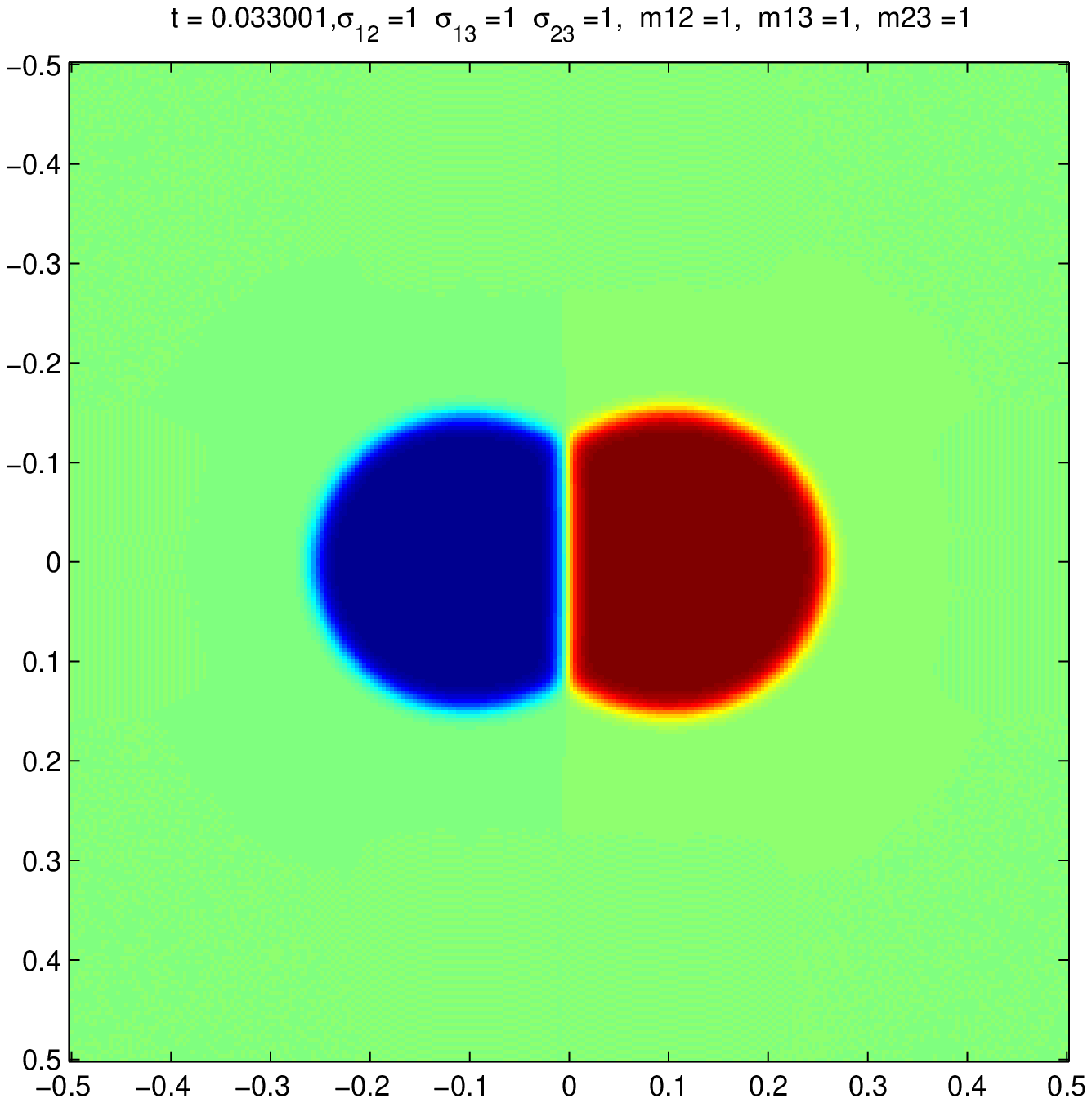}
\includegraphics[width=3.5cm]{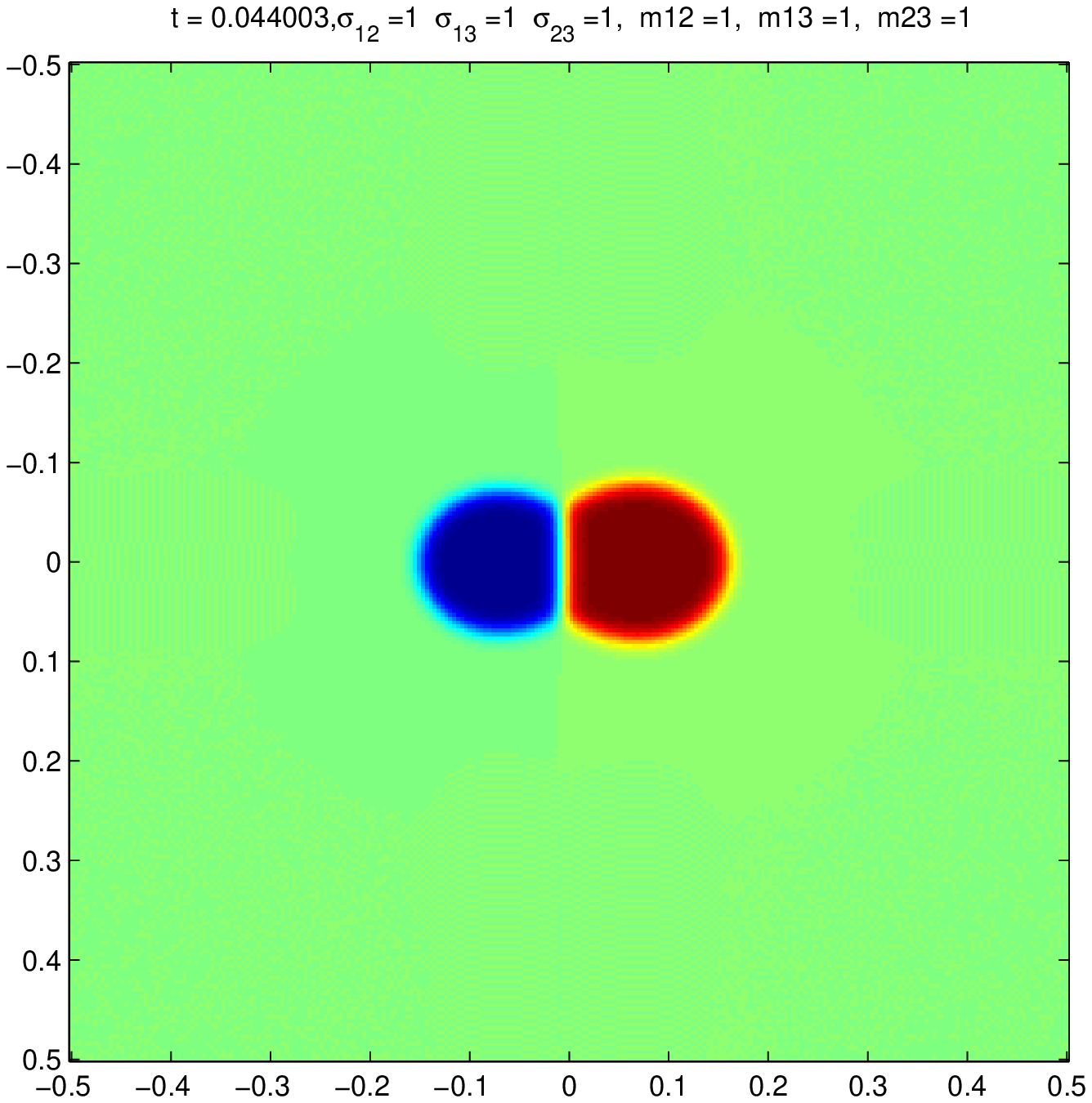} \\
\includegraphics[width=3.5cm]{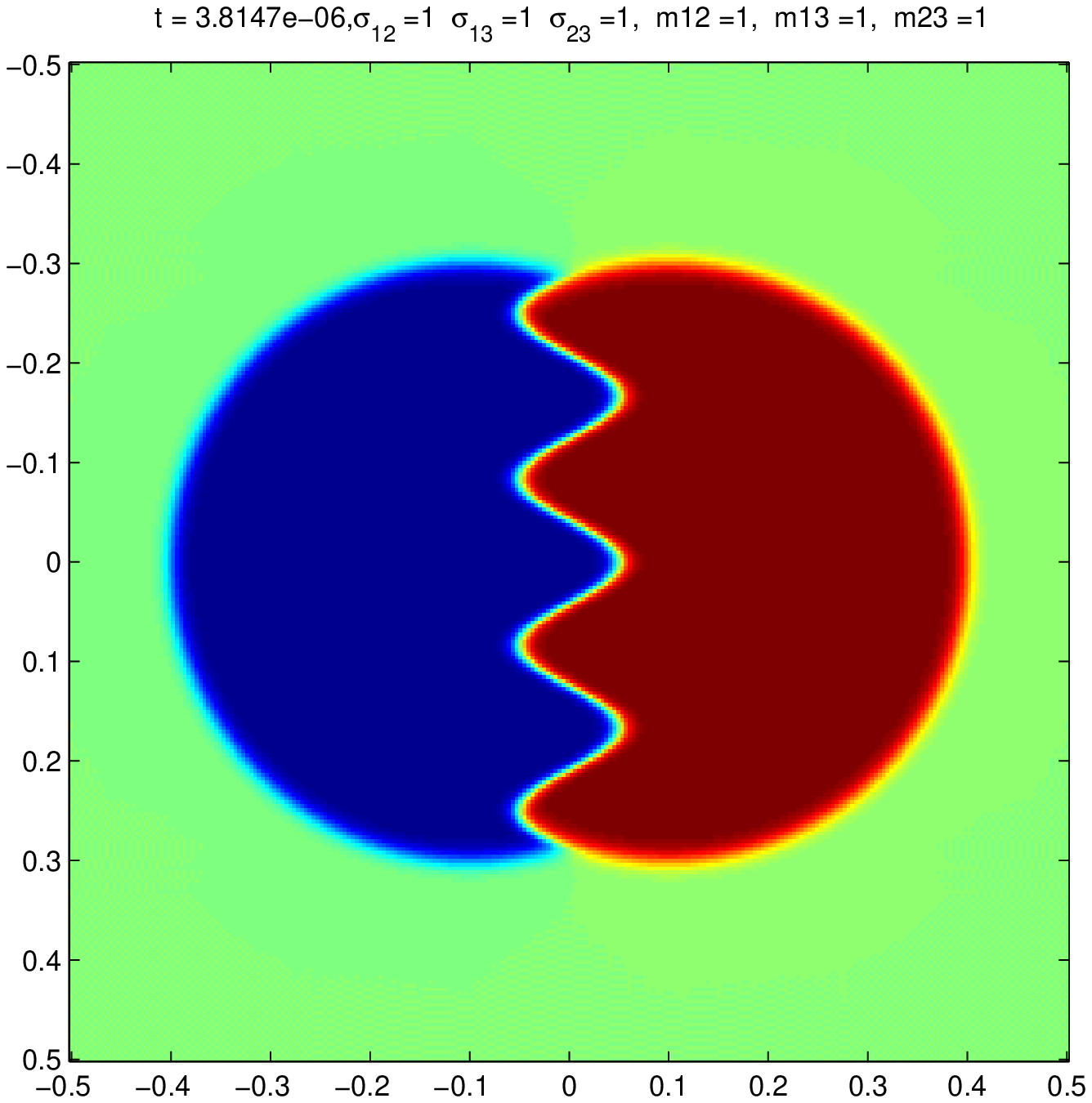}
\includegraphics[width=3.5cm]{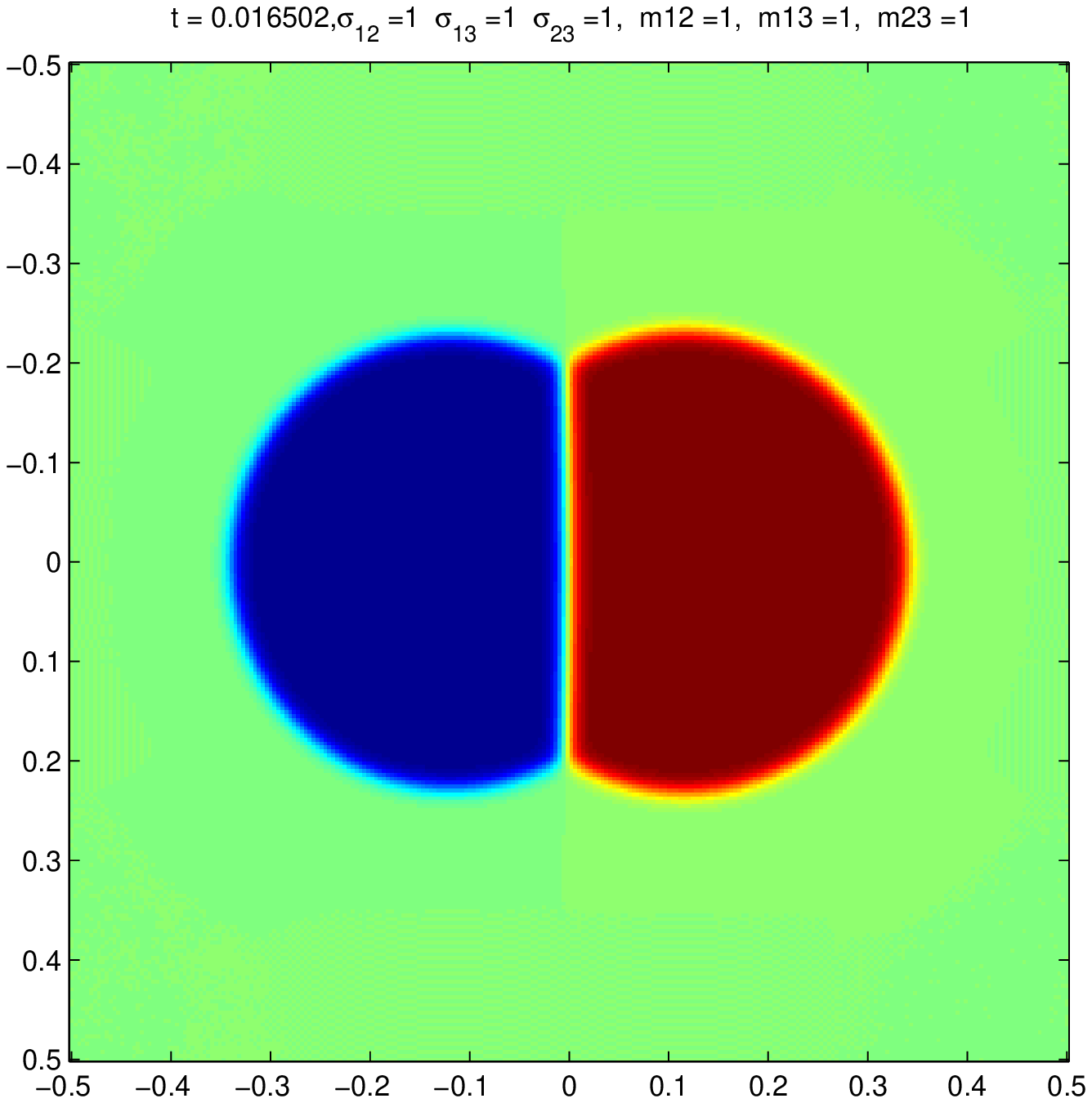}
\includegraphics[width=3.5cm]{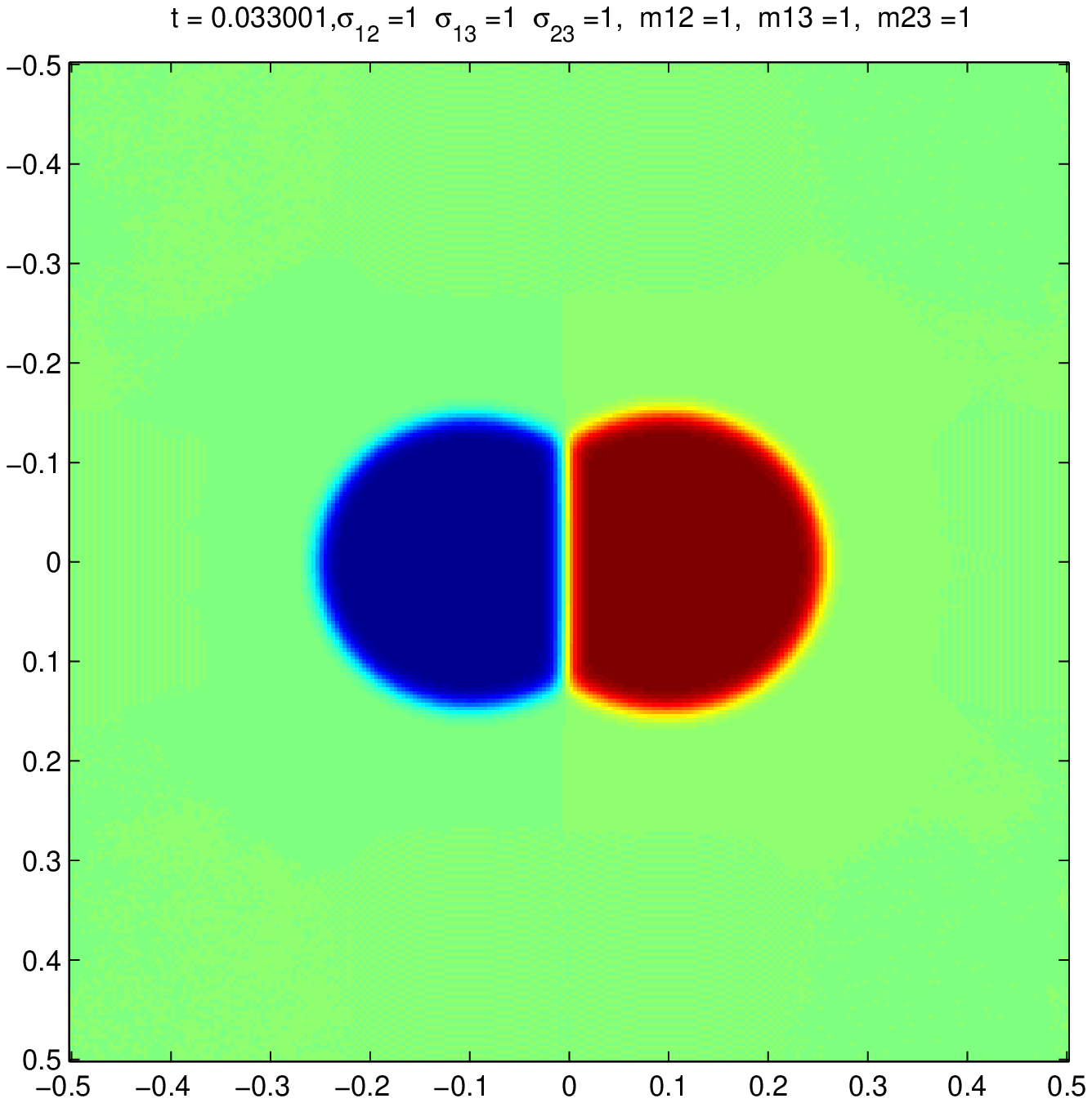}
\includegraphics[width=3.5cm]{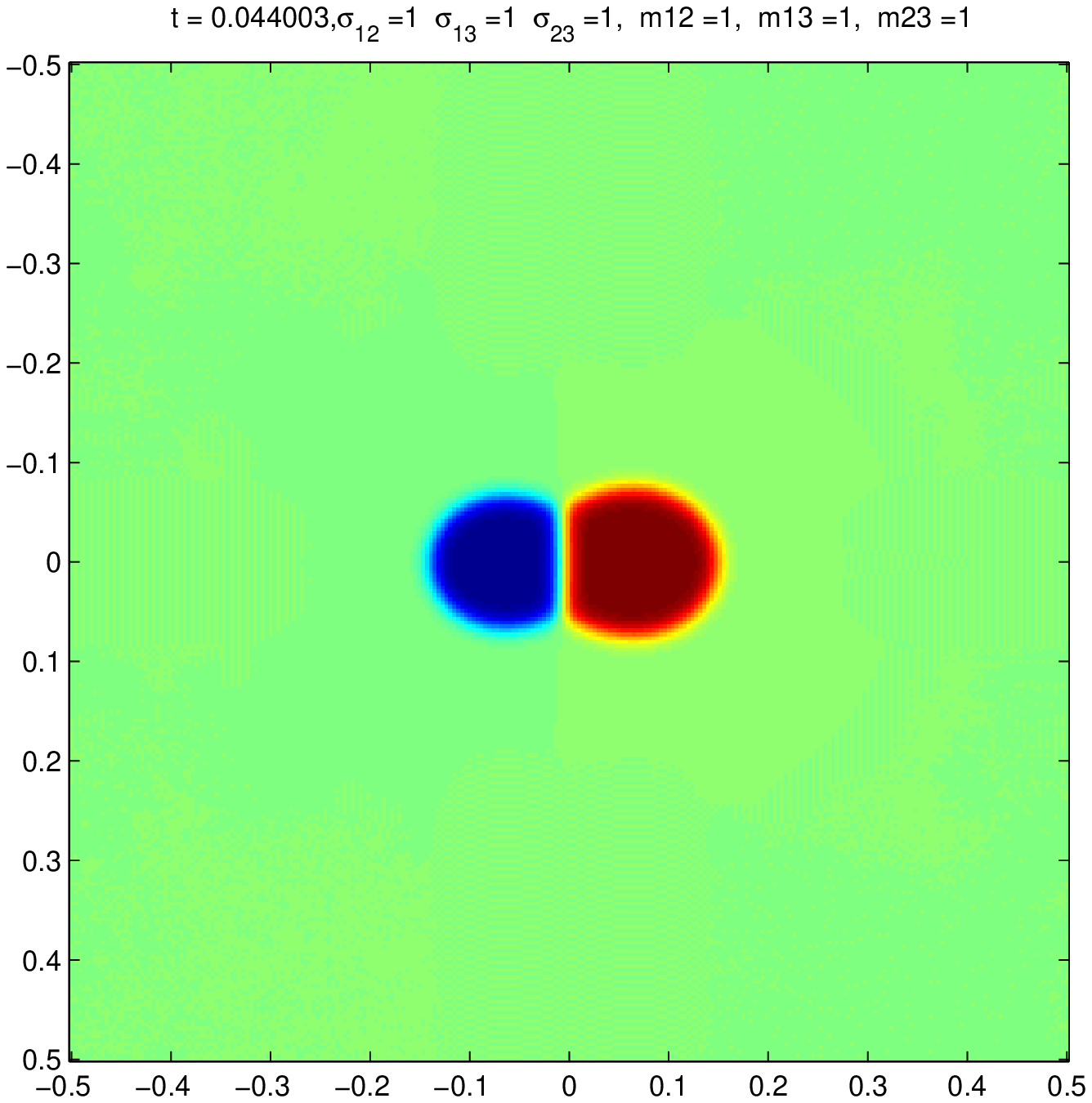} \\

\caption{Influence of the choice of a particular decomposition of the mobilities: 
the canonical decomposition is used on the first row, a sparse decomposition on the second.
The images show the values of $2u_2 + u_3$ at different times, {with the suitable colormap so that} $u_1$, $u_2$, and $u_3$ are represented in blue, red and green, respectively. }
\label{fig_hardecom1}
\end{figure}

The following numerical parameters $\varepsilon = 1.5/K$, $\delta_t = 0.25/K^2$, and  $\alpha = 0$
are used.
Figure \ref{fig_hardecom1} shows the numerical multiphase solution 
${\bf u}^{\varepsilon} = (u_1^{\varepsilon},u_2^{\varepsilon},u_3^{\varepsilon})$ 
at different times.
The rows correspond to the canonical and sparse decomposition of 
the $m_{ij}$'s, respectively .  
We observe that the two flows are quite similar, which suggests that the choice of a particular decomposition has little influence.  

\subsection{Validation of our approach for highly contrasted mobilities}

\begin{figure}[htbp]
\centering
	\includegraphics[width=3.5cm]{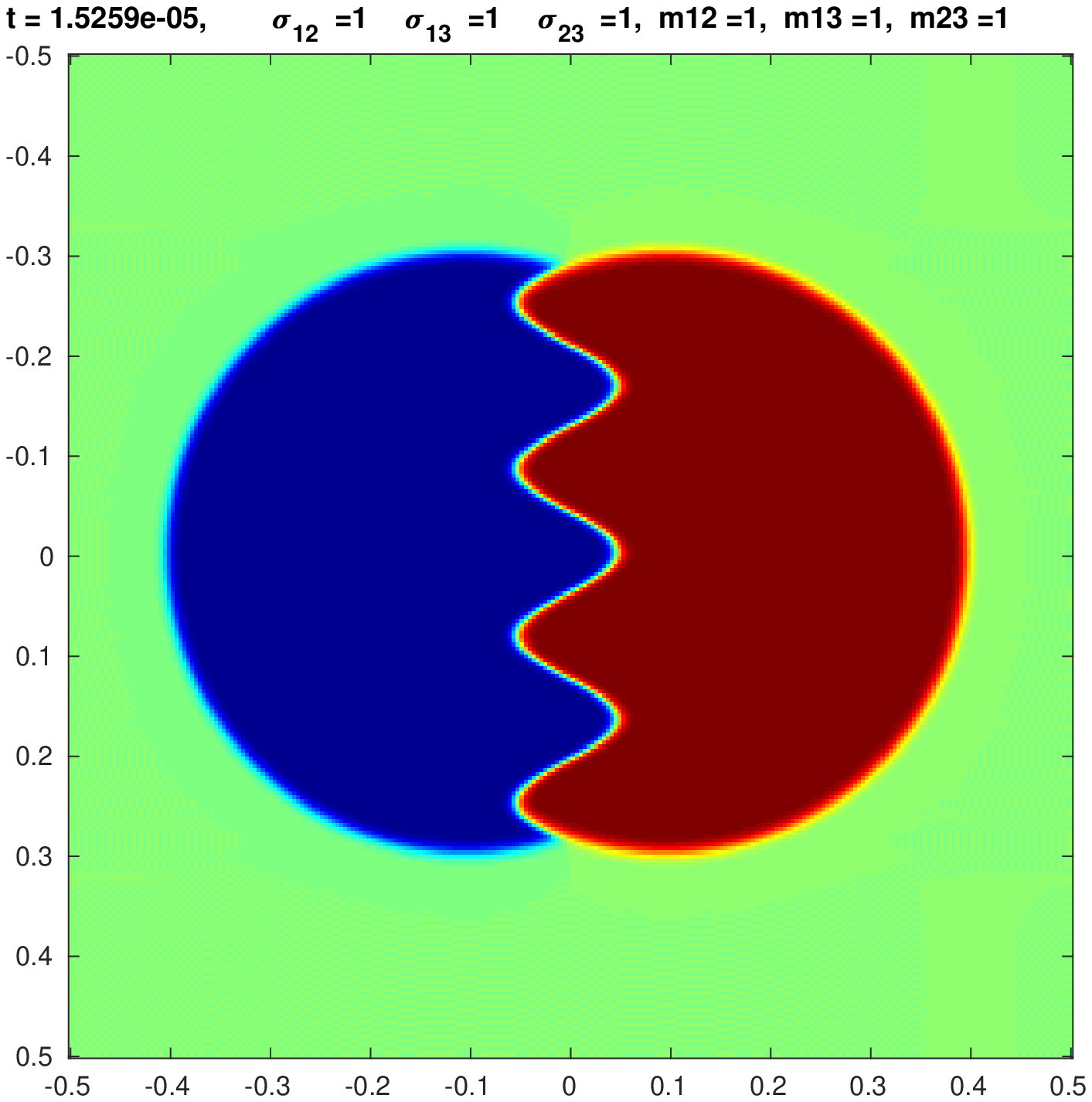}
	\includegraphics[width=3.5cm]{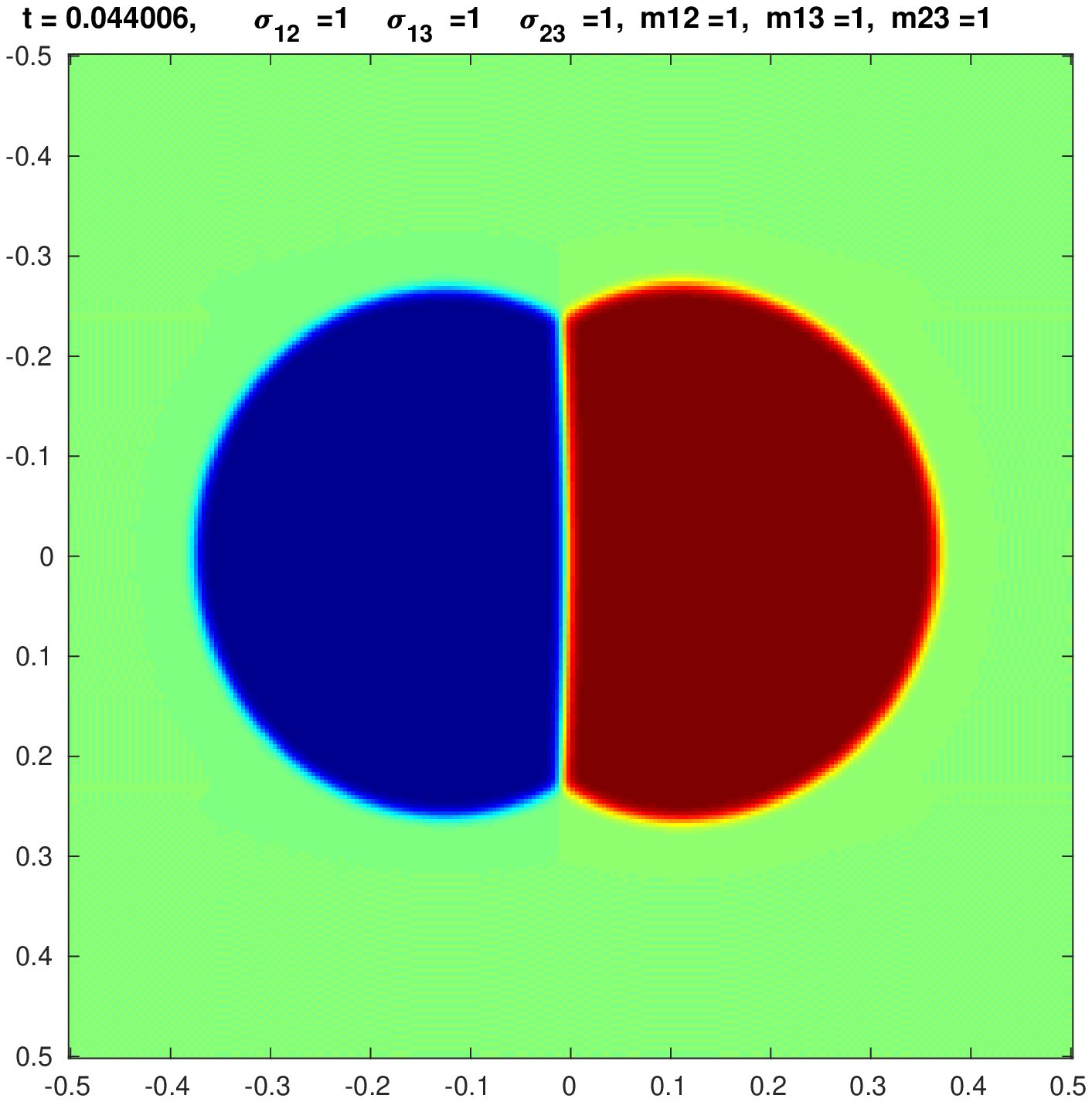}
	\includegraphics[width=3.5cm]{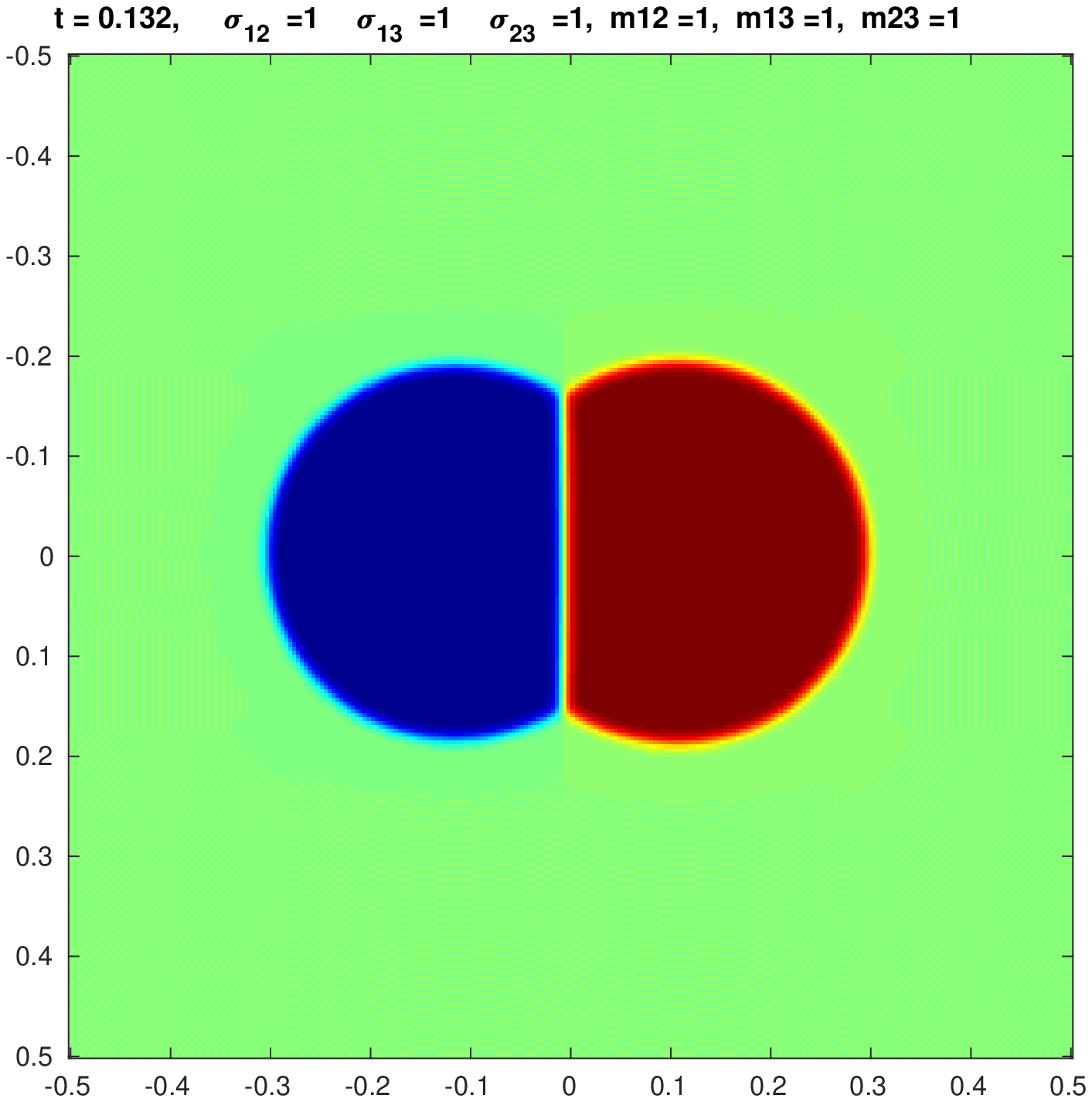}
	\includegraphics[width=3.5cm]{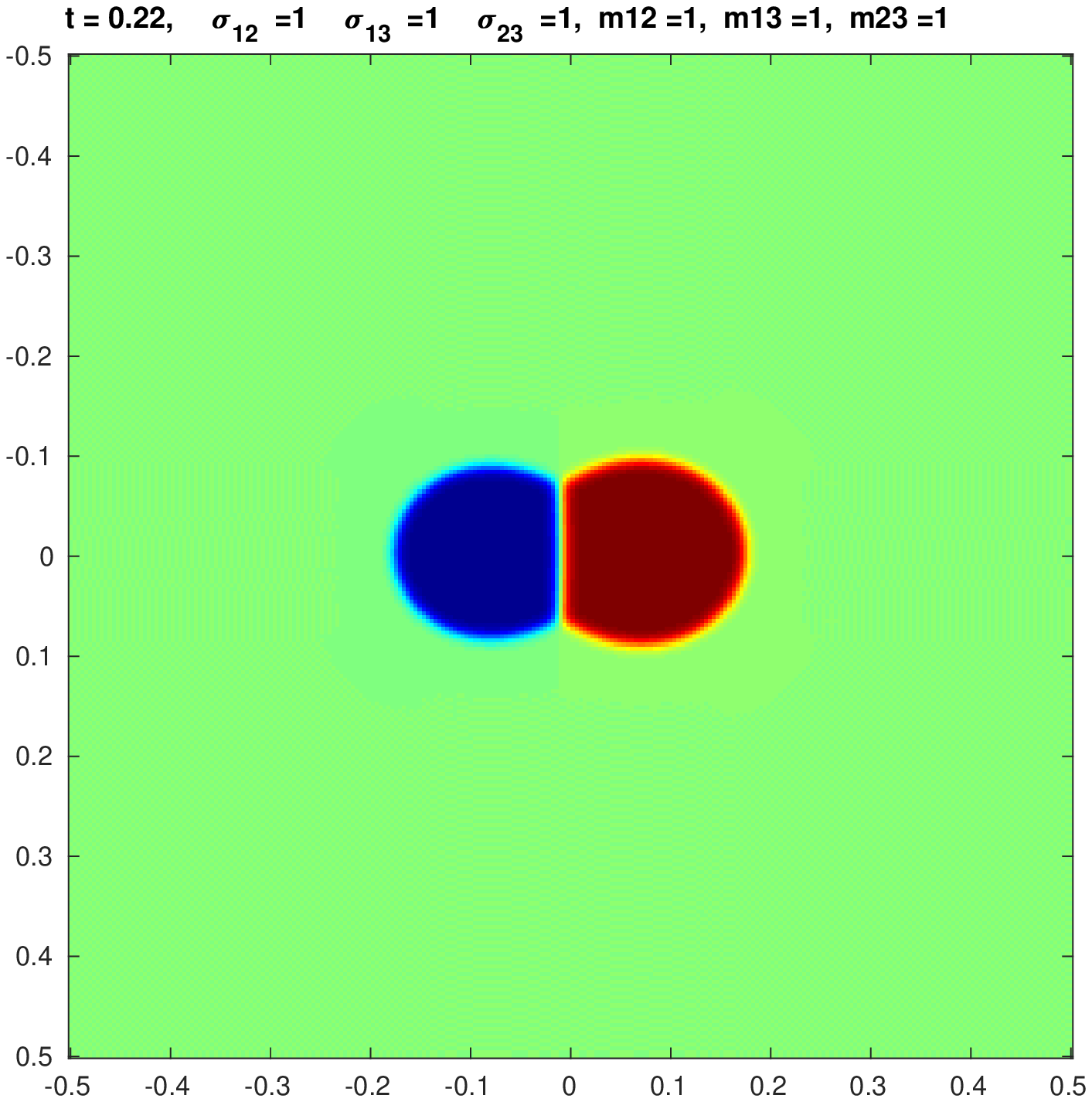} \\
	\includegraphics[width=3.5cm]{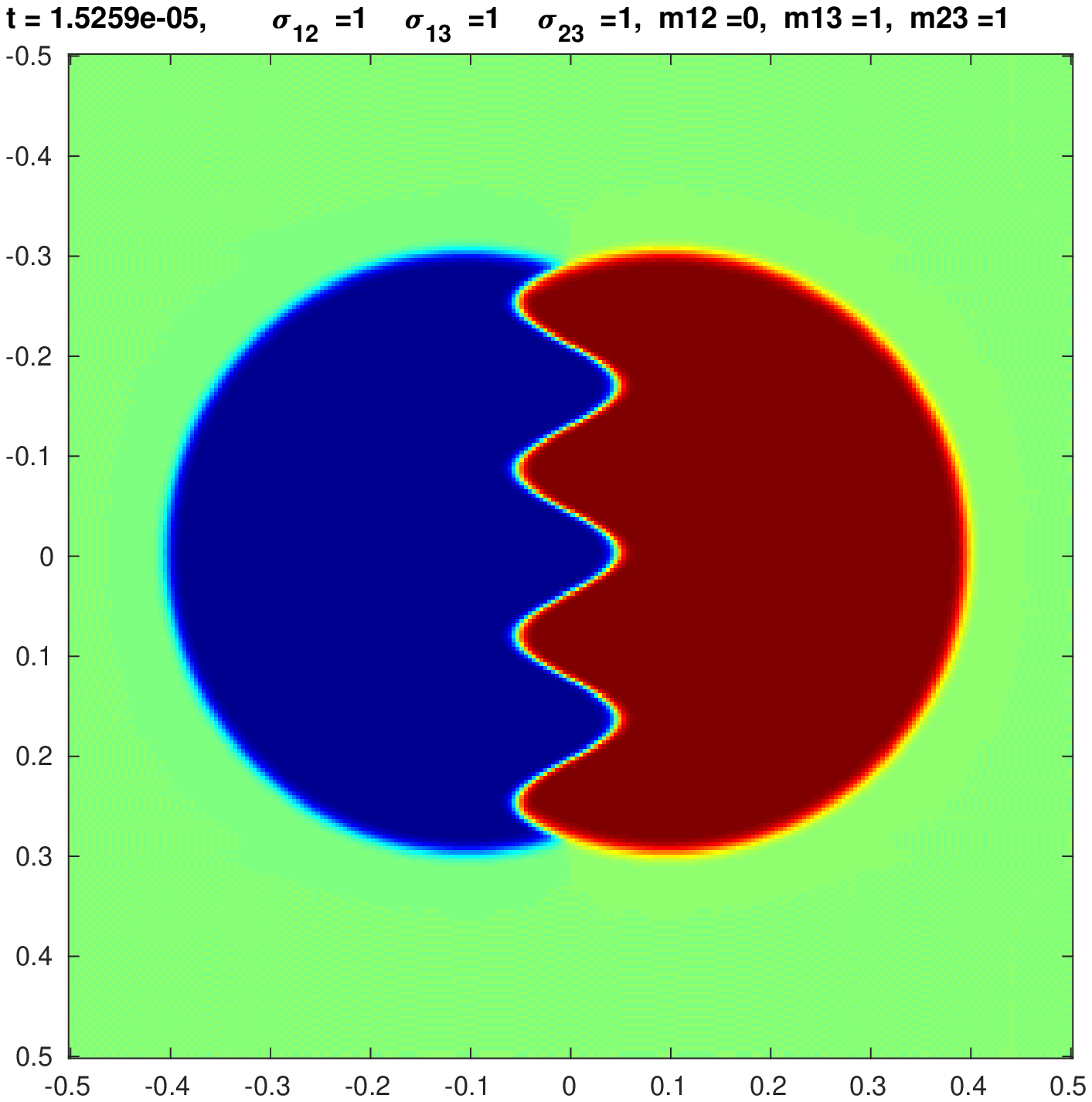}
	\includegraphics[width=3.5cm]{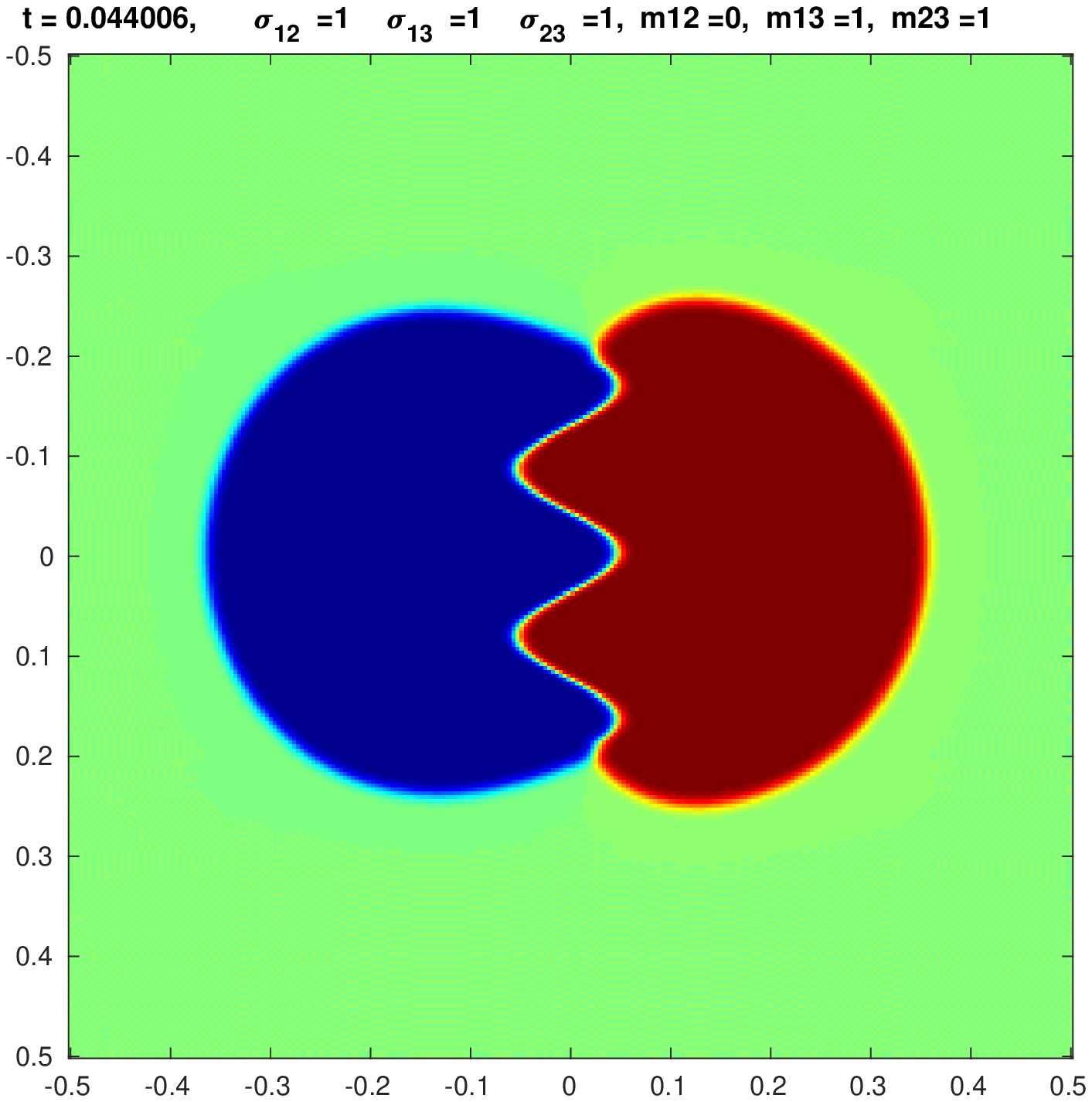}
	\includegraphics[width=3.5cm]{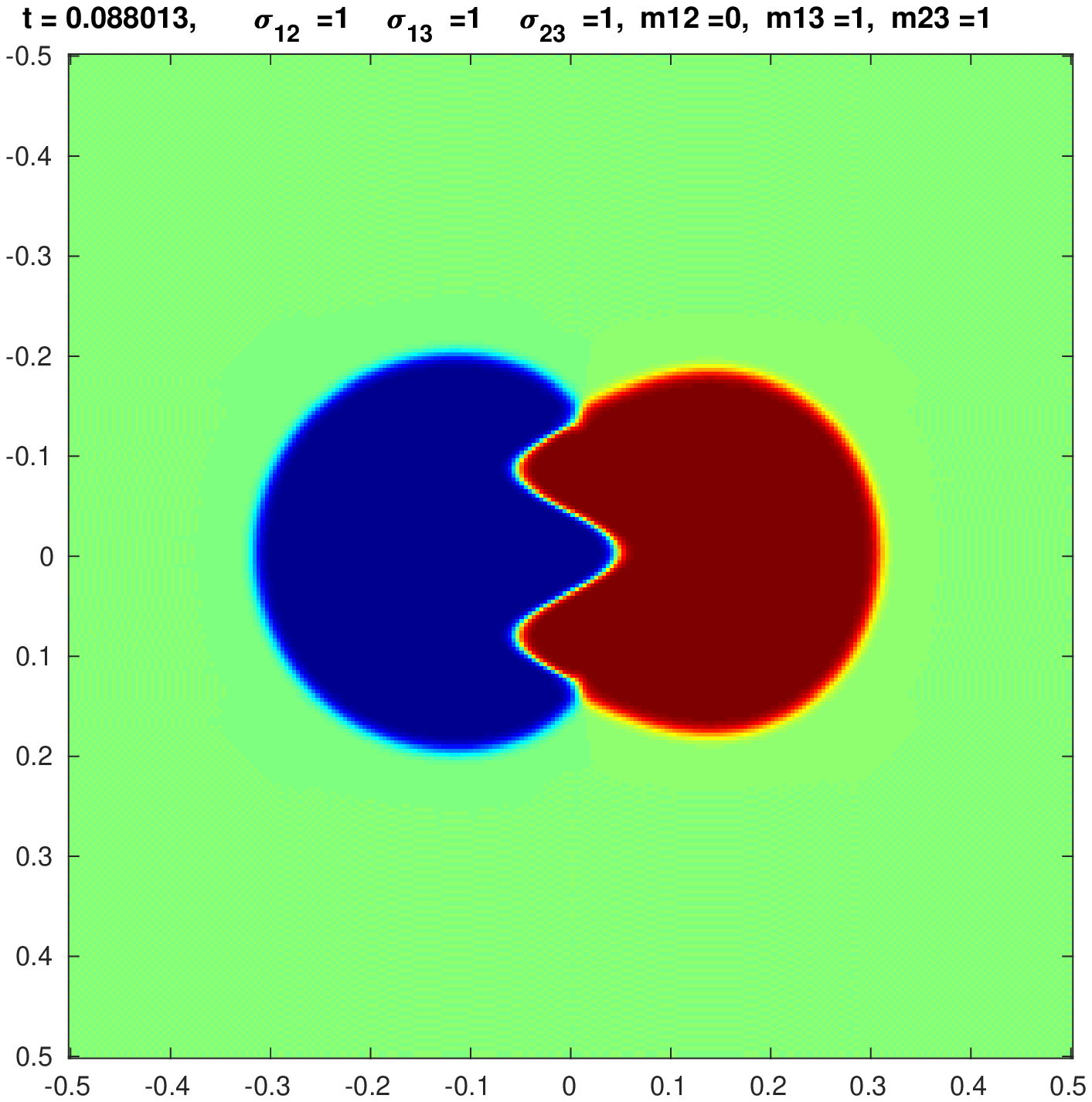}
	\includegraphics[width=3.5cm]{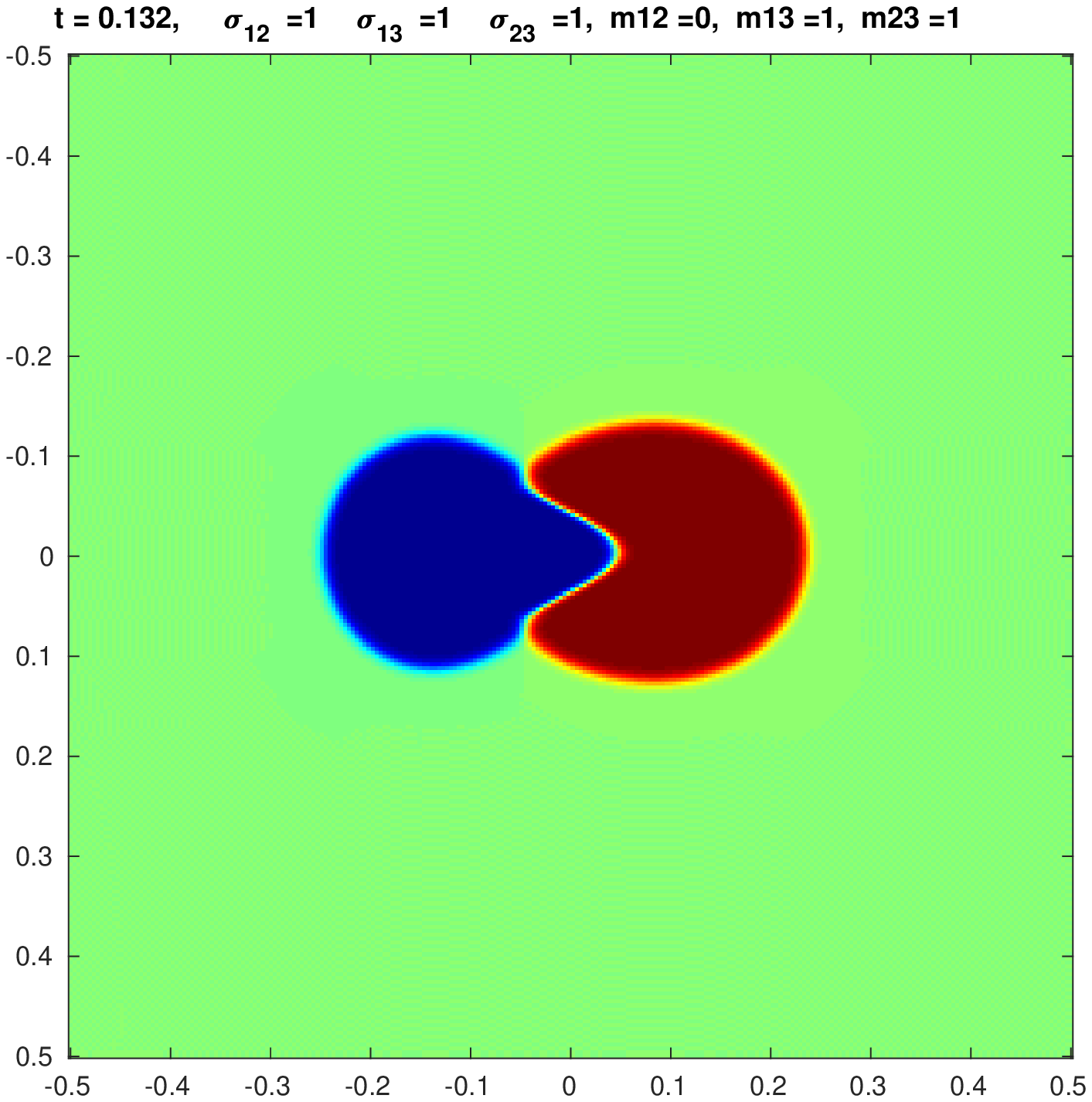} \\
	\includegraphics[width=3.5cm]{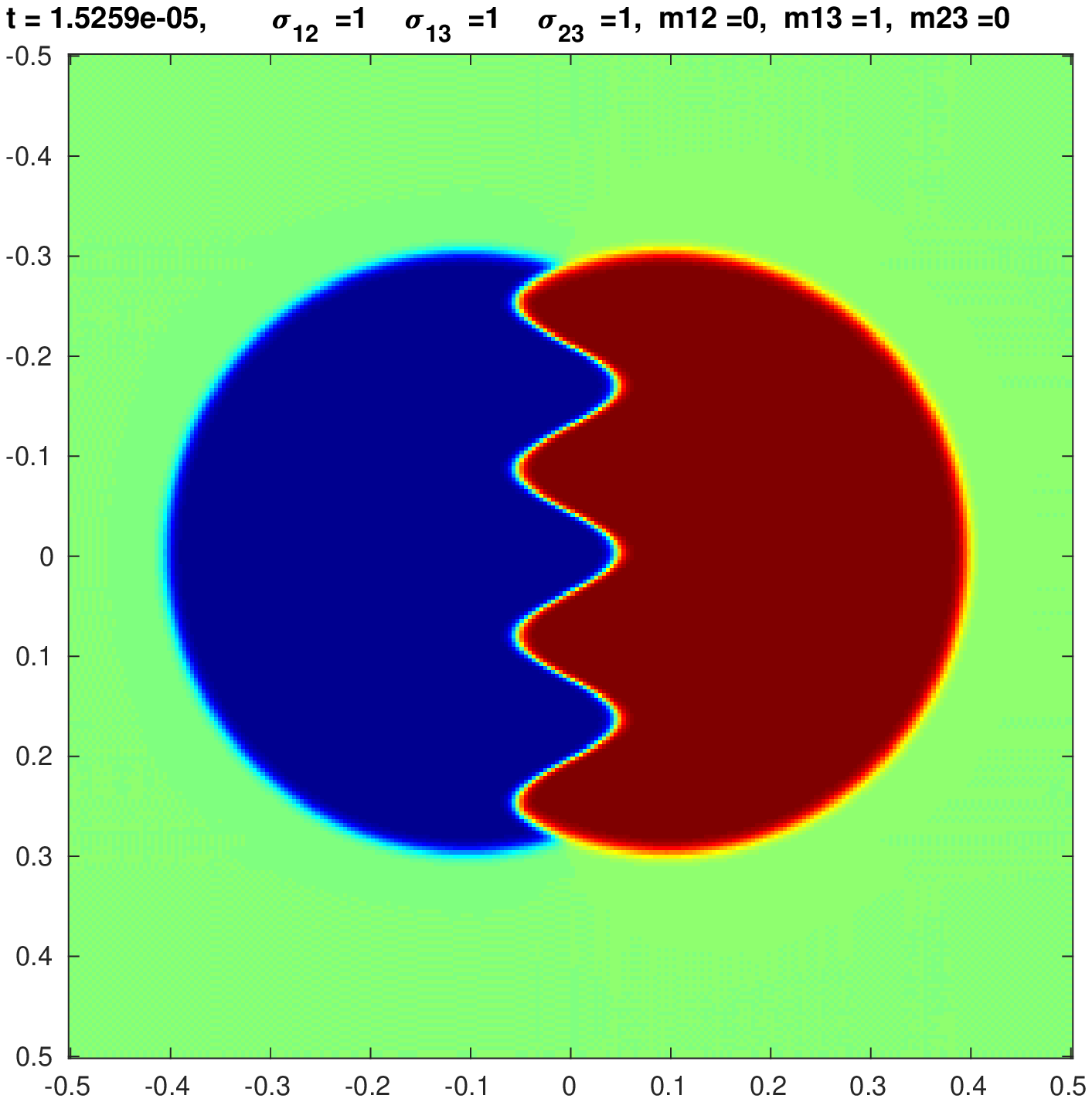}
	\includegraphics[width=3.5cm]{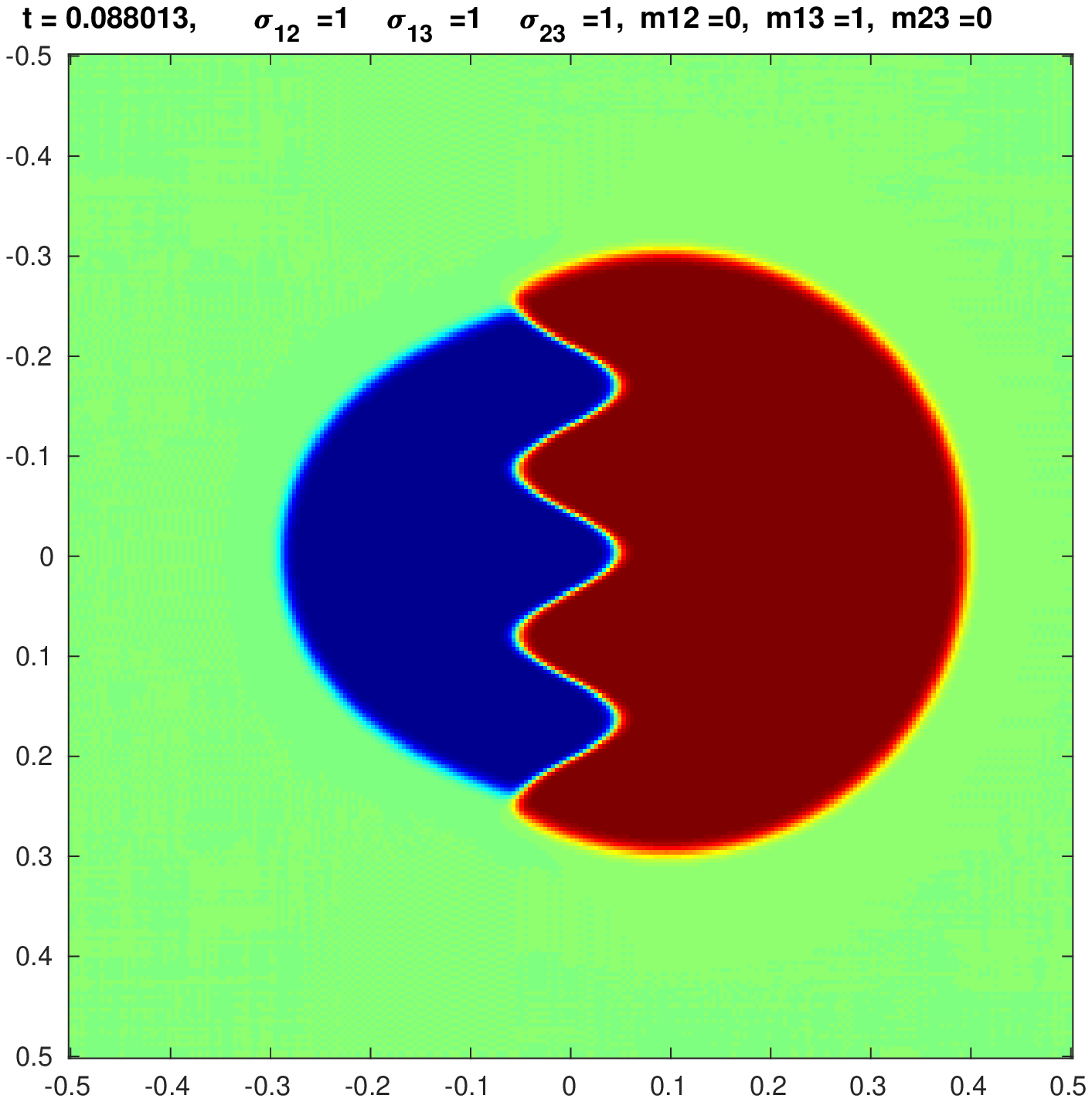}
	\includegraphics[width=3.5cm]{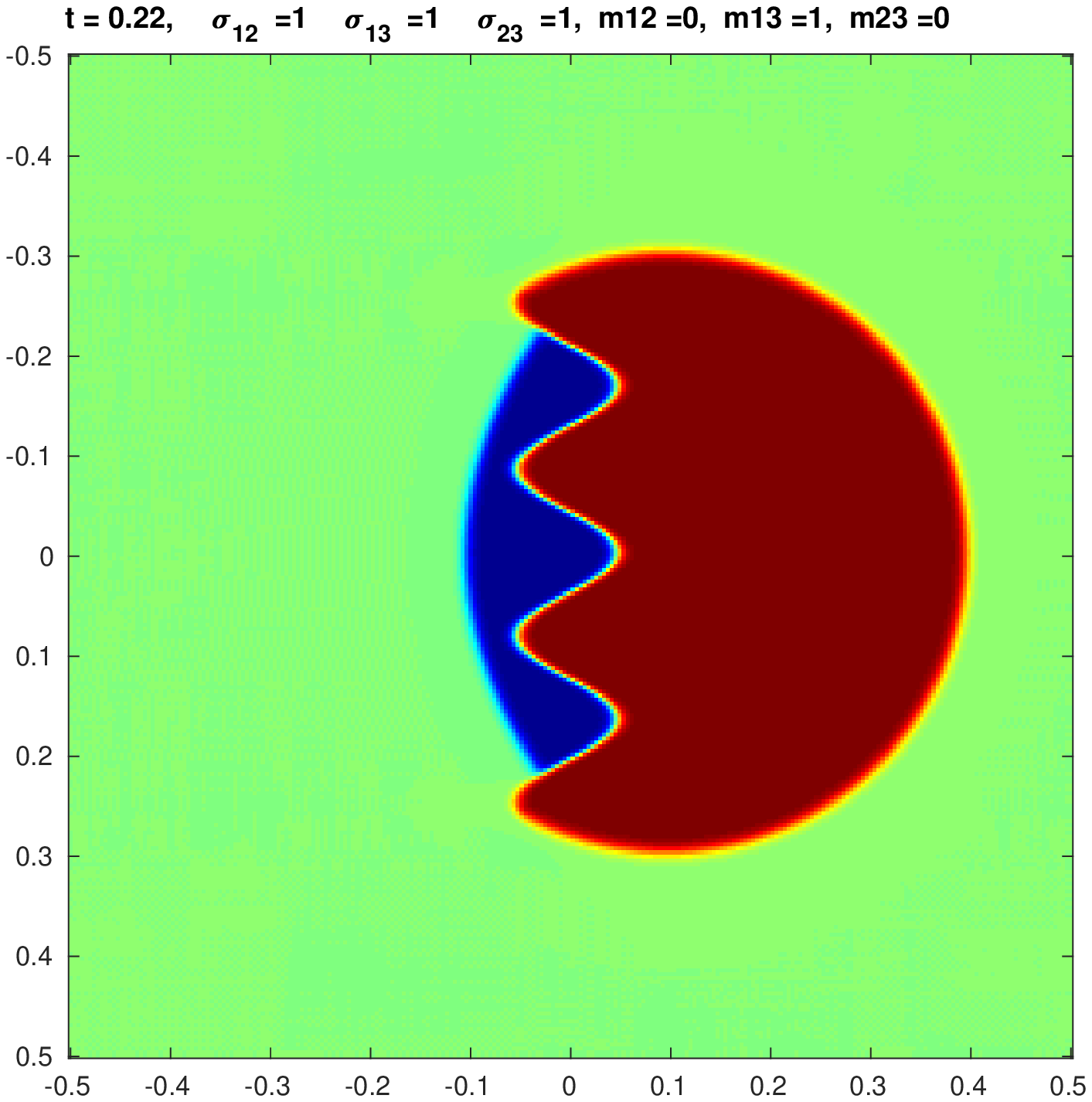}
	\includegraphics[width=3.5cm]{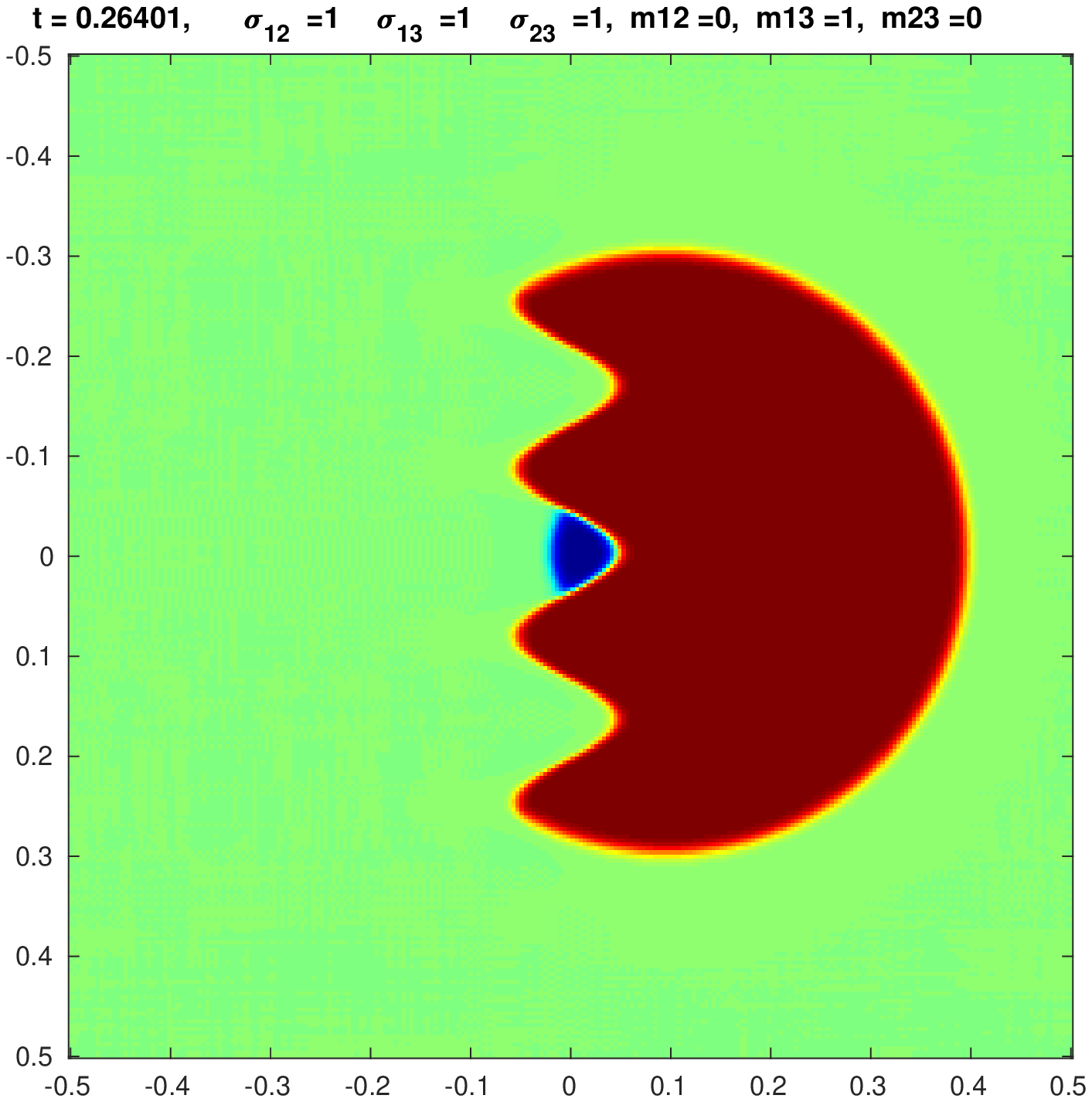} \\
\caption{Mean curvature flows with highly contrasted mobilities,
identical surface tensions.
The rows correspond to $(m_{12},m_{13},m_{23}) = (1,1,1)$, $(m_{12},m_{13},m_{23}) = (0,1,1)$, and   
$(m_{12},m_{13},m_{23}) = (0,1,0)$, respectively.
Images show the values of  the function $2u_2 + u_3$ at different times {using a colormap such that} $u_1$, $u_2$, and $u_3$ appear in blue, red and green, respectively. }
\label{fig_compar1}
\end{figure}

Our next tests show that our approach can handle highly contrasted mobilities.
One expects that when $m_{ij}$ is small (or vanishes) the corresponding interface
$\Gamma_{ij}$ hardly moves.
The tests also show that mobilities are parameters that may strongly affect the flow.
The computations have been performed with  
$\varepsilon = 1/K$, $\delta_t = 1/K^2$, and  $\alpha = 2$. 
Figure \ref{fig_compar1} represents a first series of numerical experiments 
in which $\sigma_{12} = \sigma_{13} = \sigma_{23} = 1$. 
The rows depict the flow associated with
the mobilities $(m_{12},m_{13},m_{23})=(1,1,1)$, $(0,1,1)$, and  $(0,1,0)$ respectively,
with the same initial condition.
On each image, the phases $u_1$ and $u_2$  are plotted in blue and red respectively.  
As expected, the blue-red interface $\Gamma_{12}$
does not move when $m_{12}=0$ (second line),
or when $m_{23}$ (third line). 
\\

\begin{figure}[htbp]
\centering
	\includegraphics[width=3.5cm]{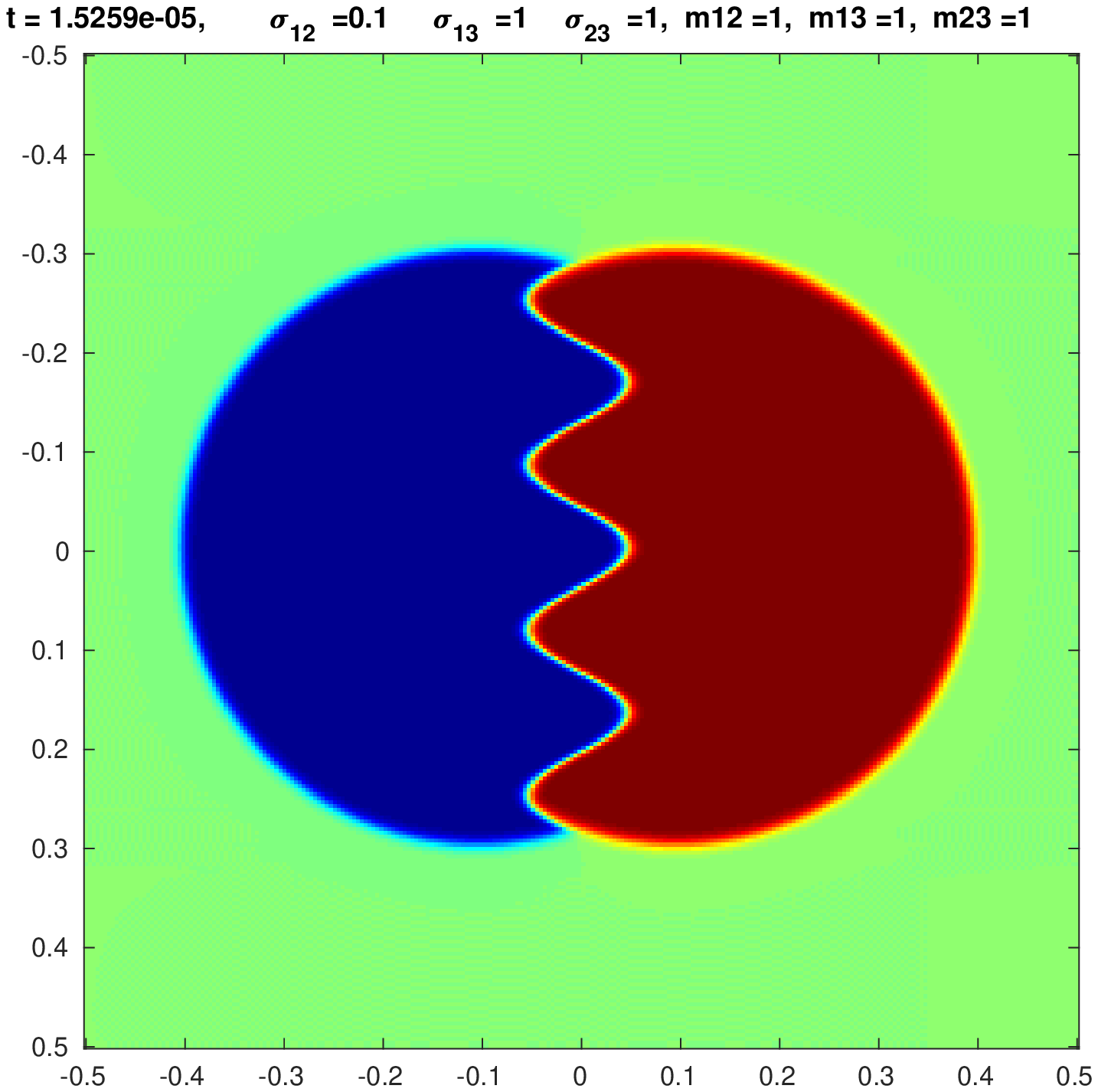}
	\includegraphics[width=3.5cm]{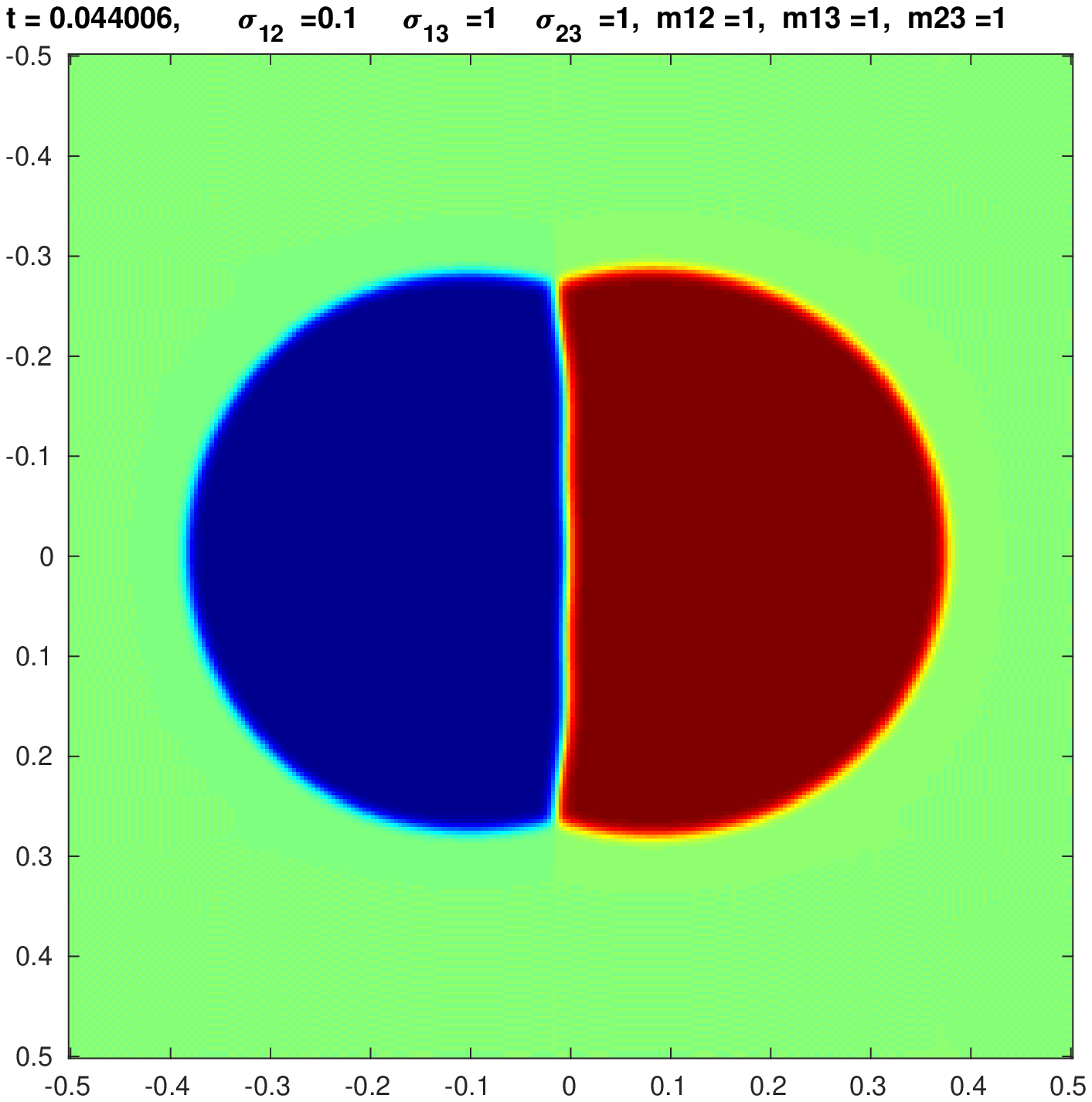}
	\includegraphics[width=3.5cm]{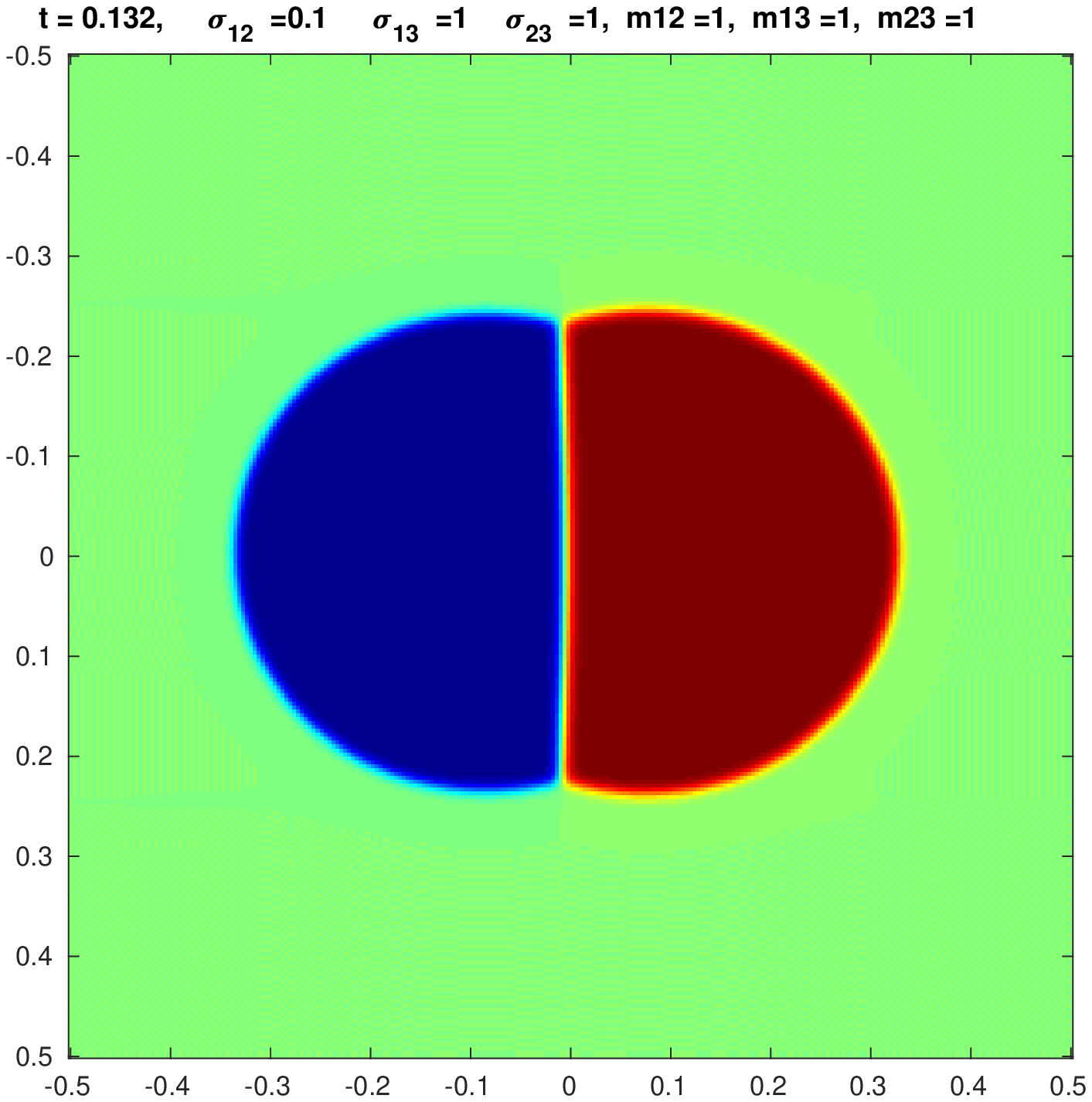}
	\includegraphics[width=3.5cm]{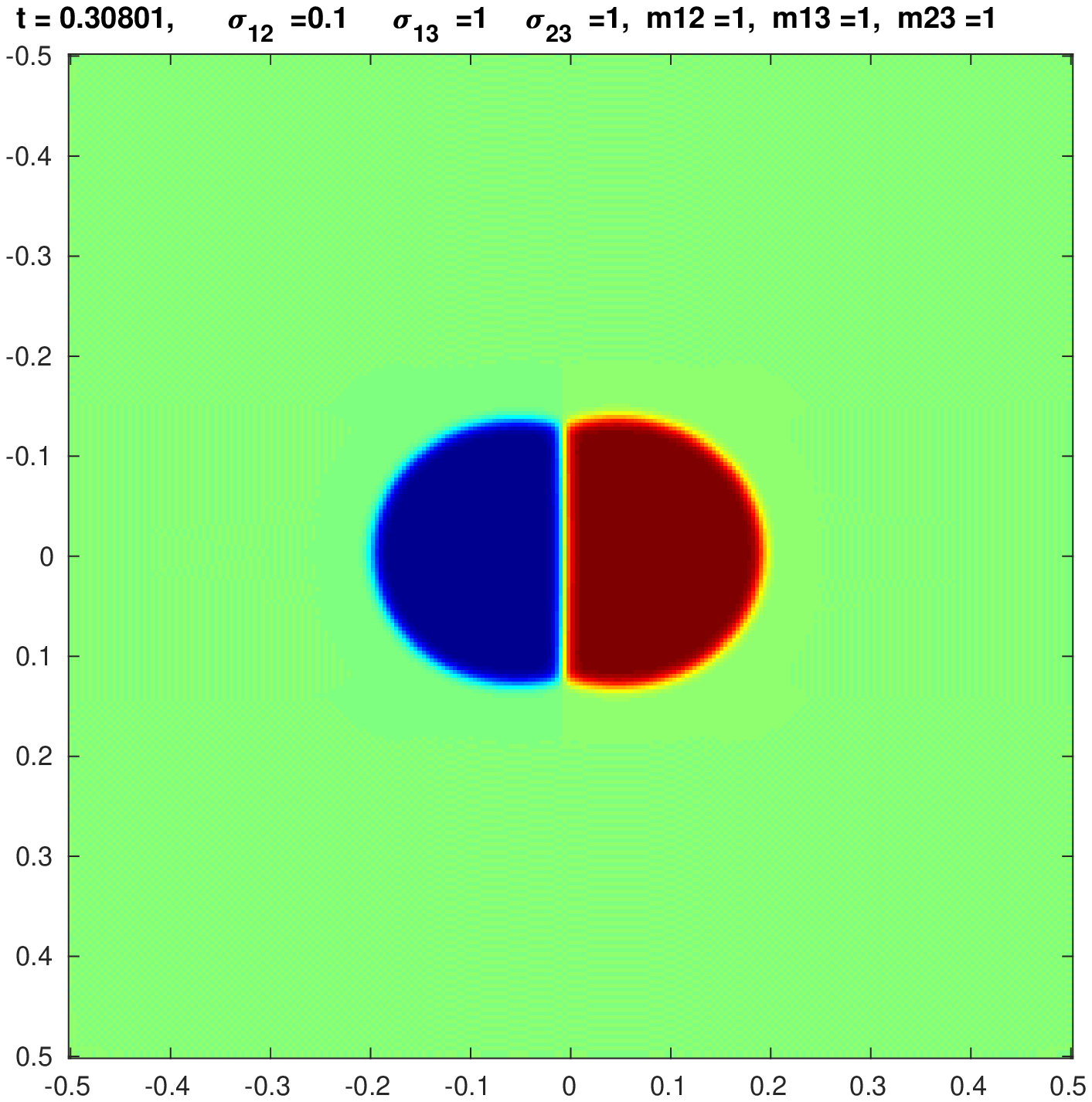} \\
	\includegraphics[width=3.5cm]{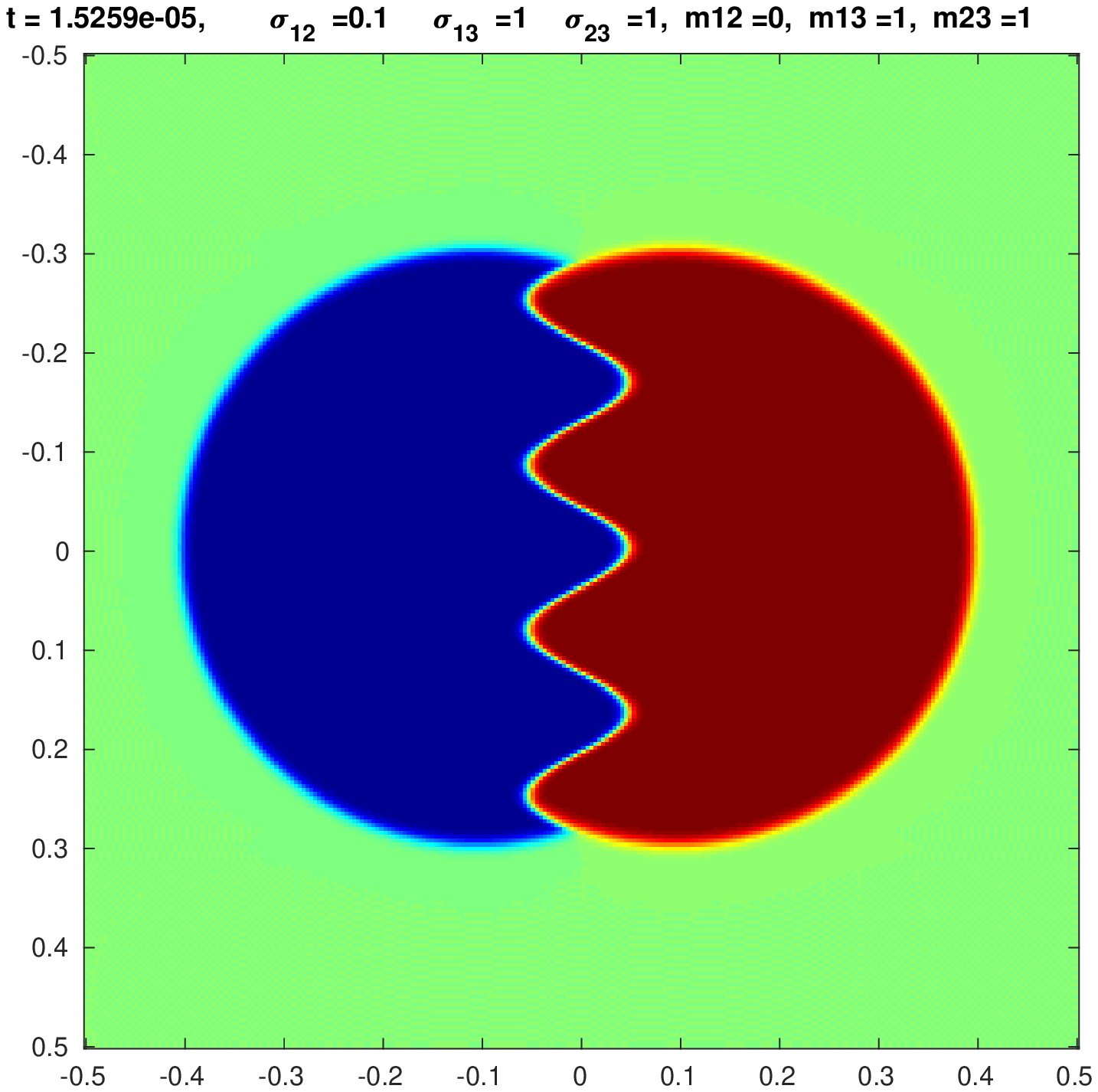}
	\includegraphics[width=3.5cm]{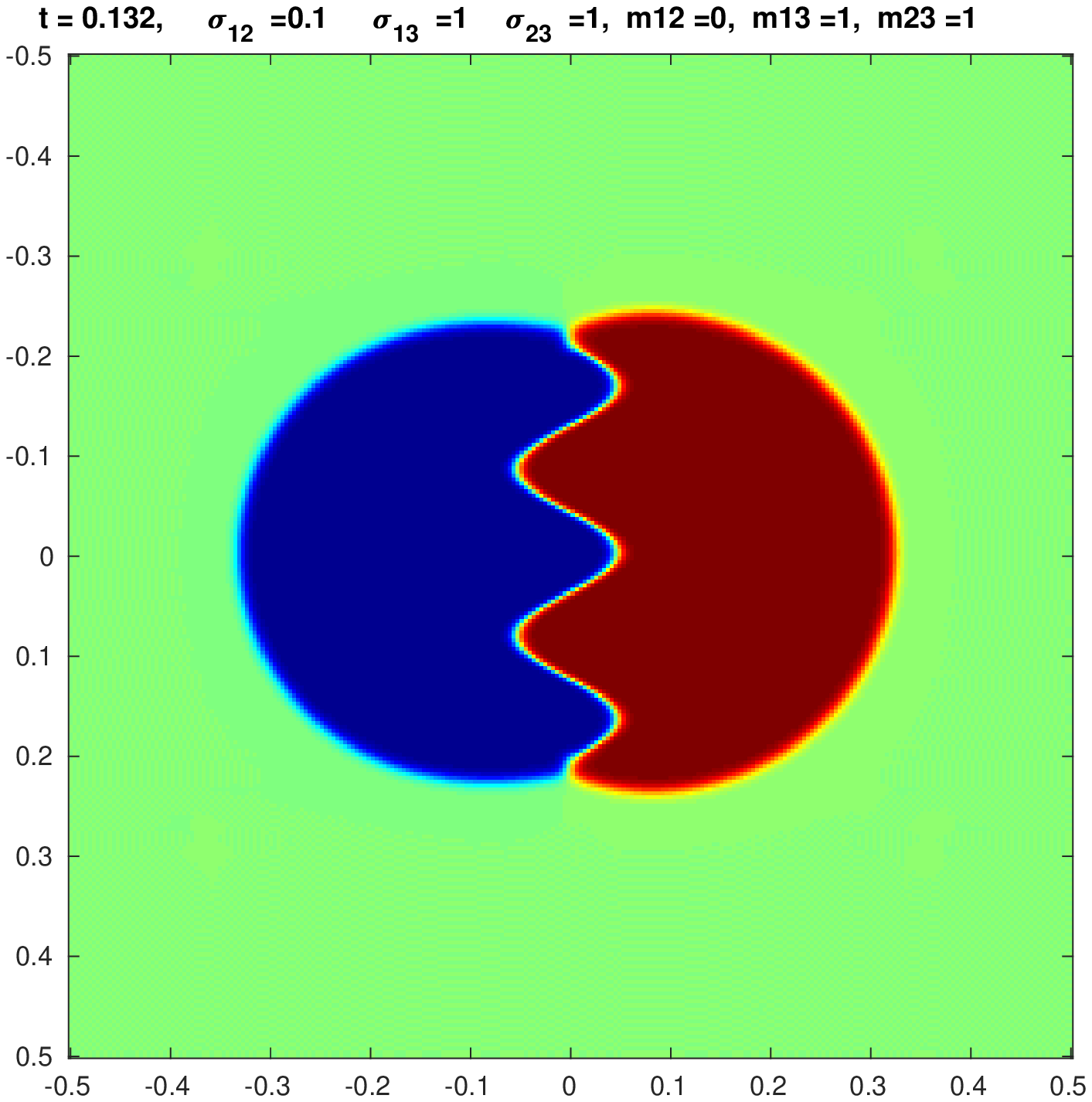}
	\includegraphics[width=3.5cm]{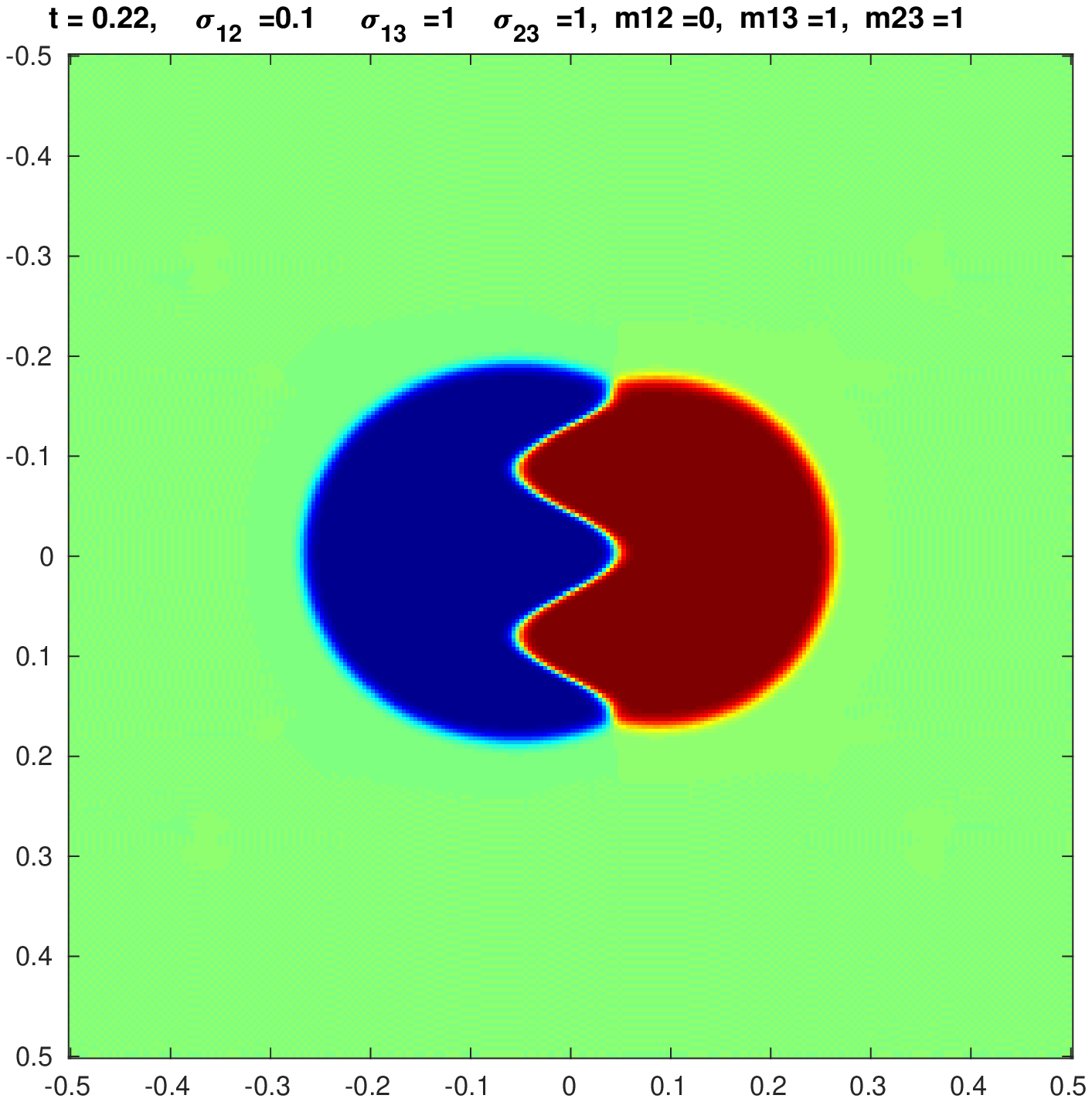}
	\includegraphics[width=3.5cm]{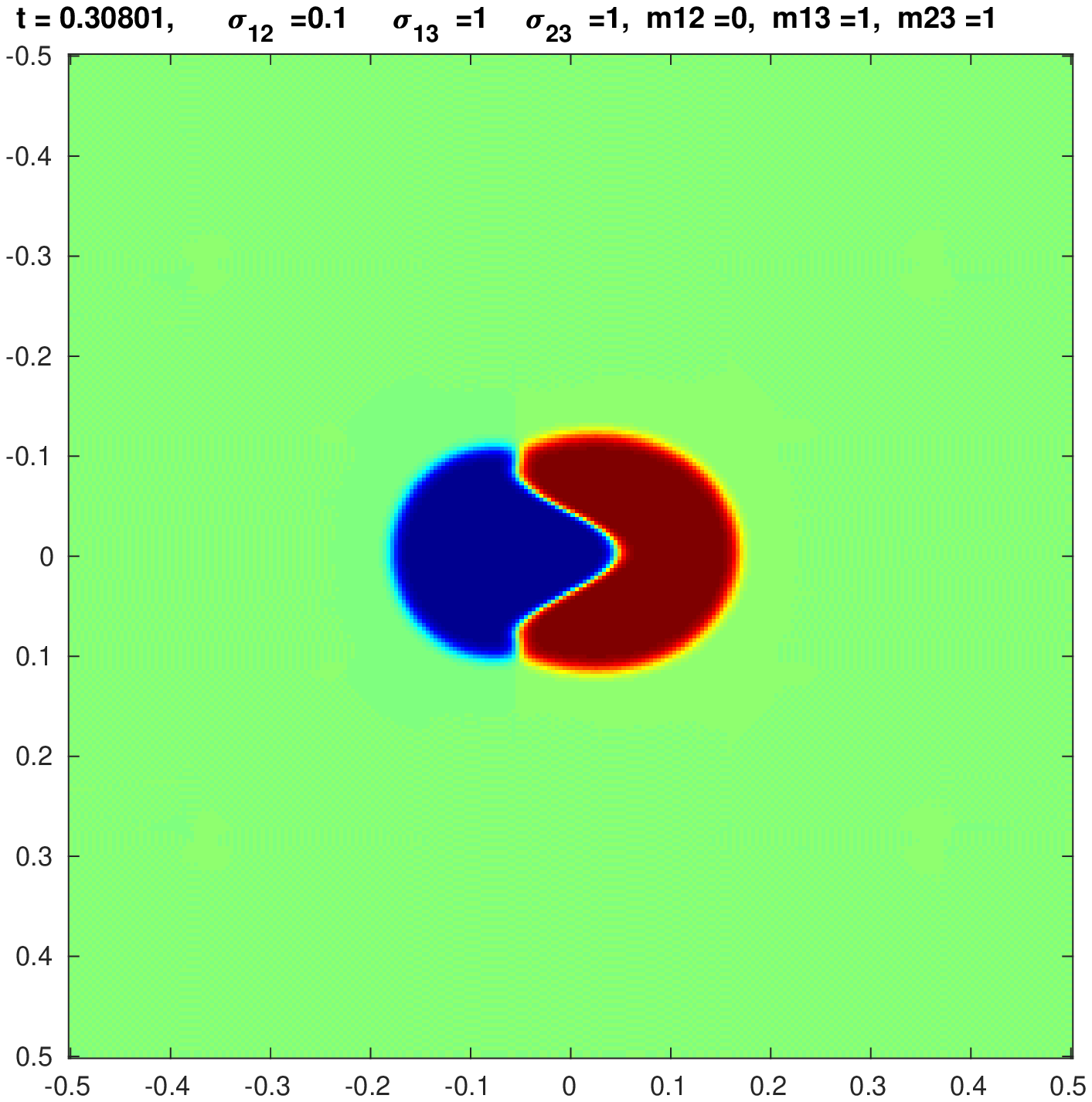} \\
	\includegraphics[width=3.5cm]{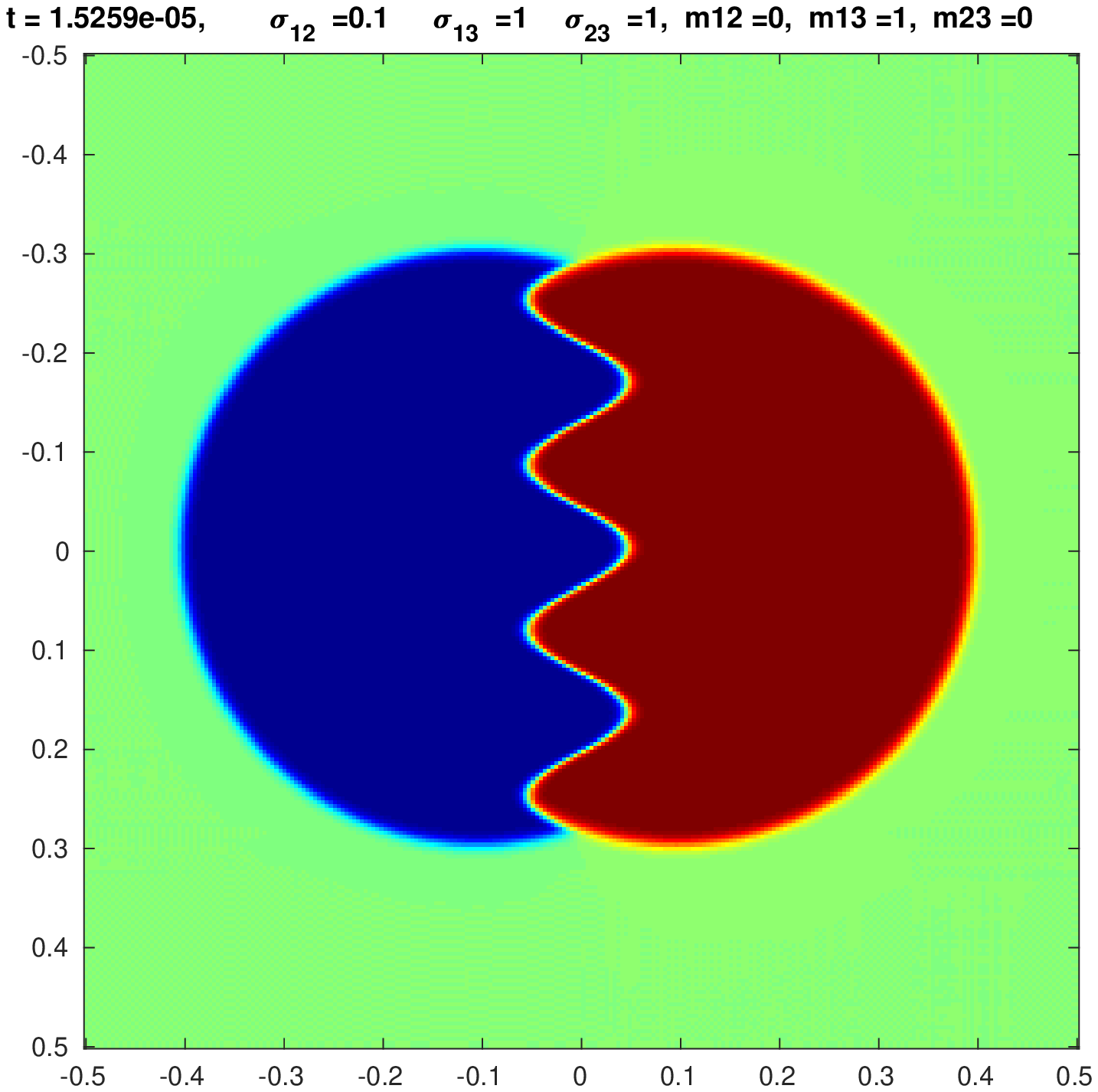}
	\includegraphics[width=3.5cm]{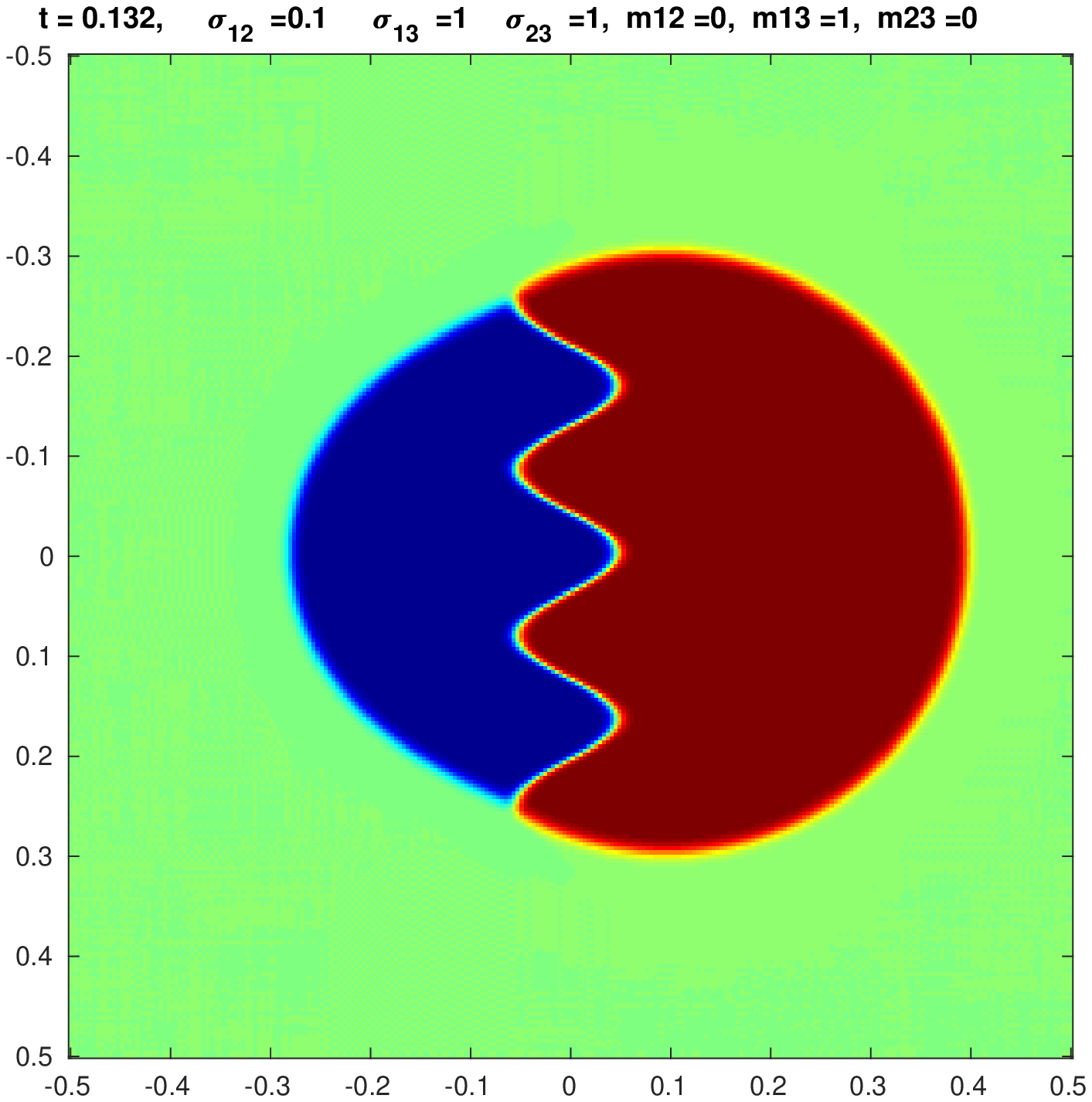}
	\includegraphics[width=3.5cm]{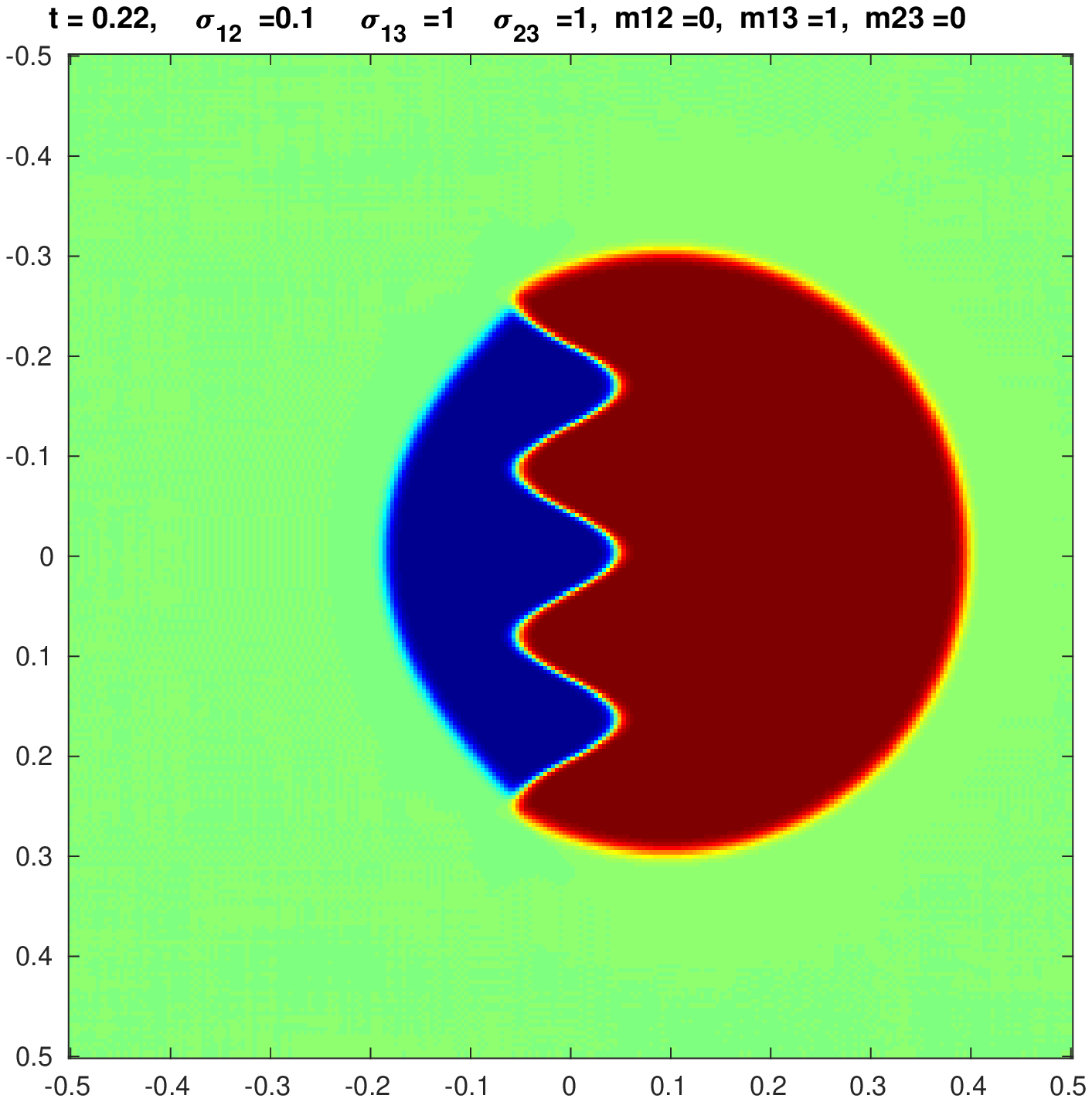}
	\includegraphics[width=3.5cm]{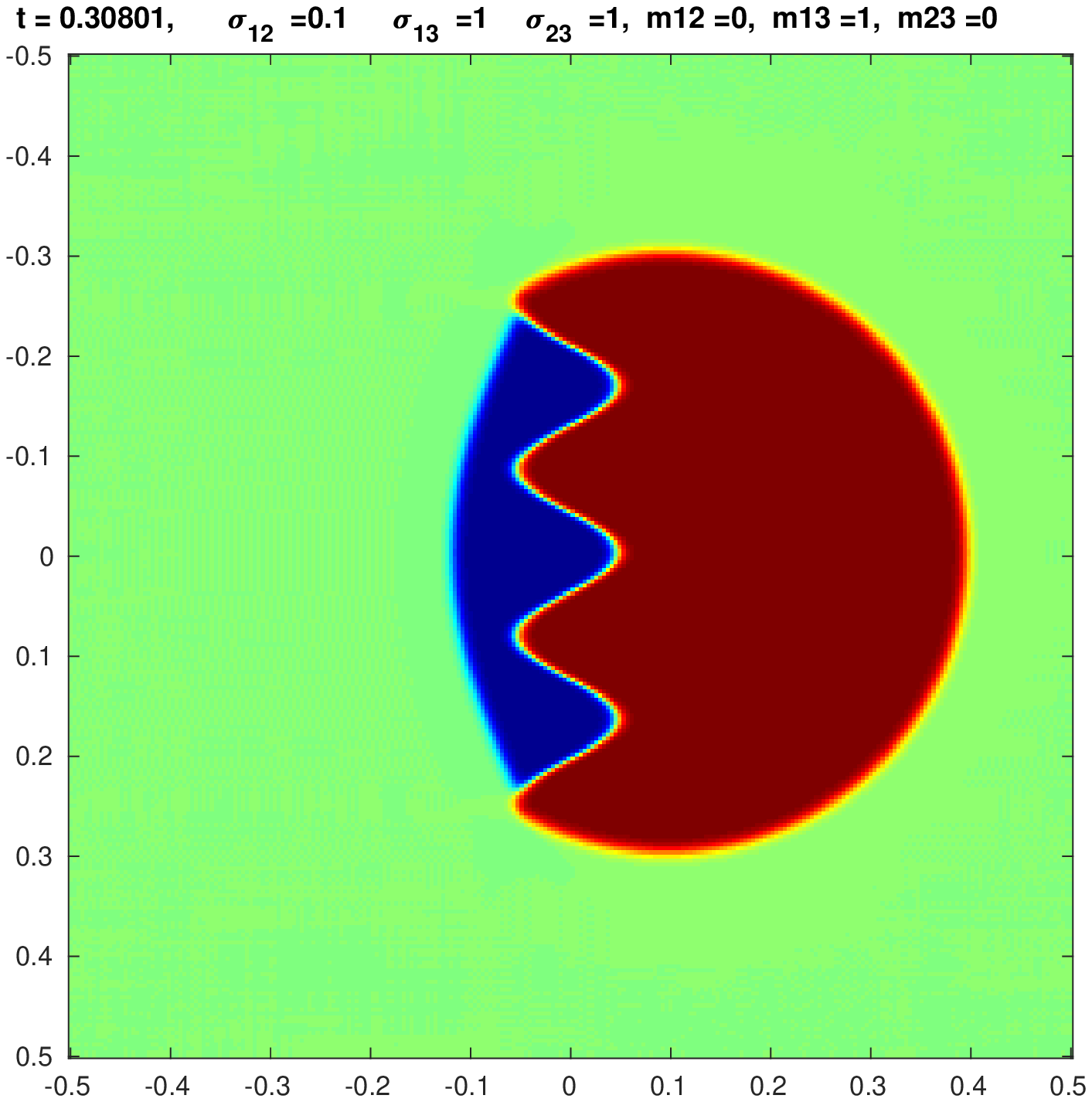} \\
	\caption{Mean curvature flows with highly contrasted mobilities,
non-identical surface tensions.
The rows correspond to $(m_{12},m_{13},m_{23}) = (1,1,1)$, $(m_{12},m_{13},m_{23}) = (0,1,1)$, and   
$(m_{12},m_{13},m_{23}) = (0,1,0)$, respectively.
Images show the values of the function $2u_2 + u_3$ at different times {using a colormap such that} $u_1$, $u_2$, and $u_3$ appear in blue, red and green,  respectively. }
\label{fig_compar2}
\end{figure}

Figure \ref{fig_compar2} represent similar experiments with the non-identical
surface tensions $\sigma_{12} = 0.1$ and $\sigma_{13} = \sigma_{23} = 1$. 
The same conclusions hold. 
 
\subsection{Numerical experiments with $N=4$ phases}
\begin{figure}[htbp]
\centering
	\includegraphics[width=3.5cm]{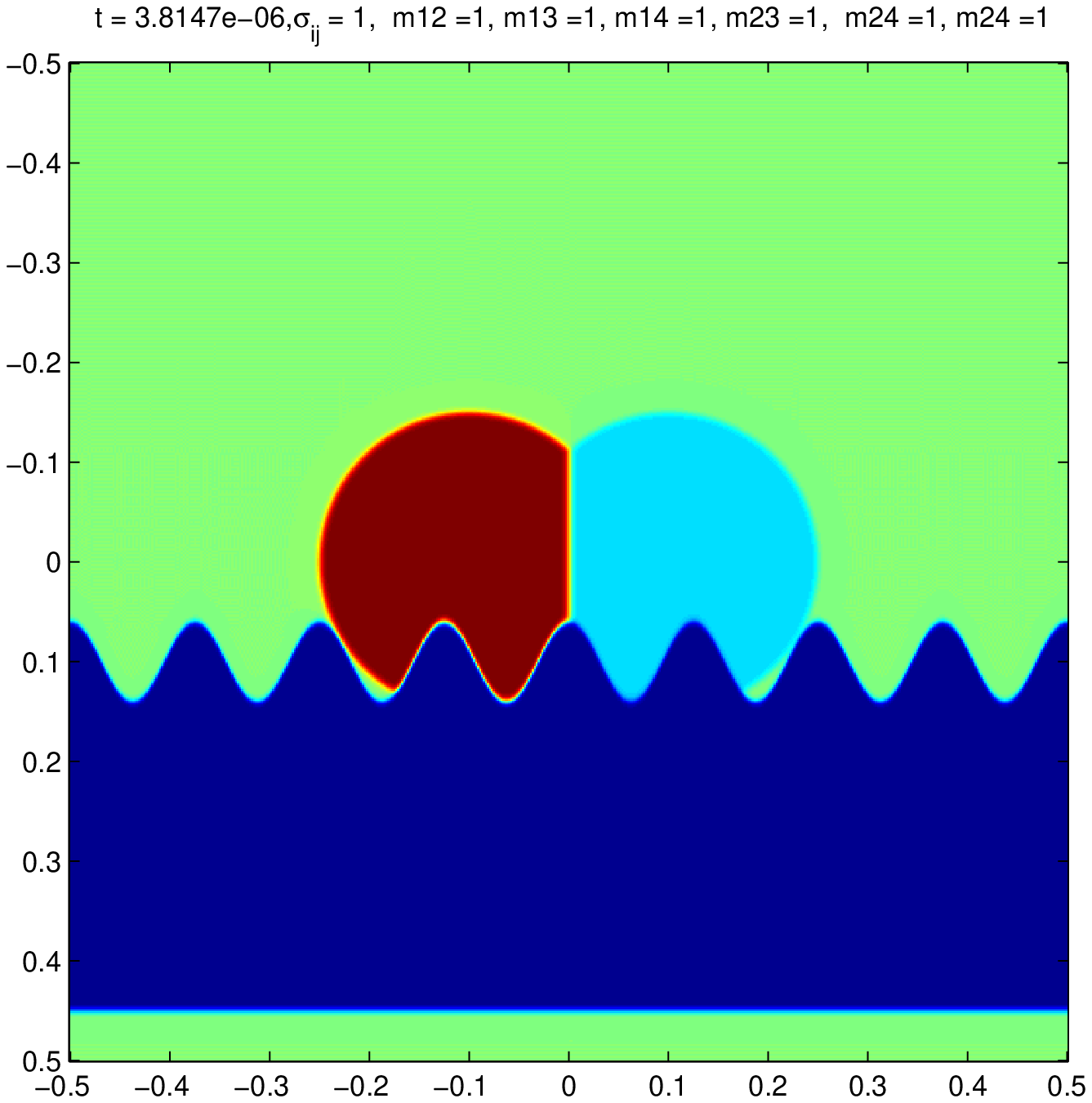}
	\includegraphics[width=3.5cm]{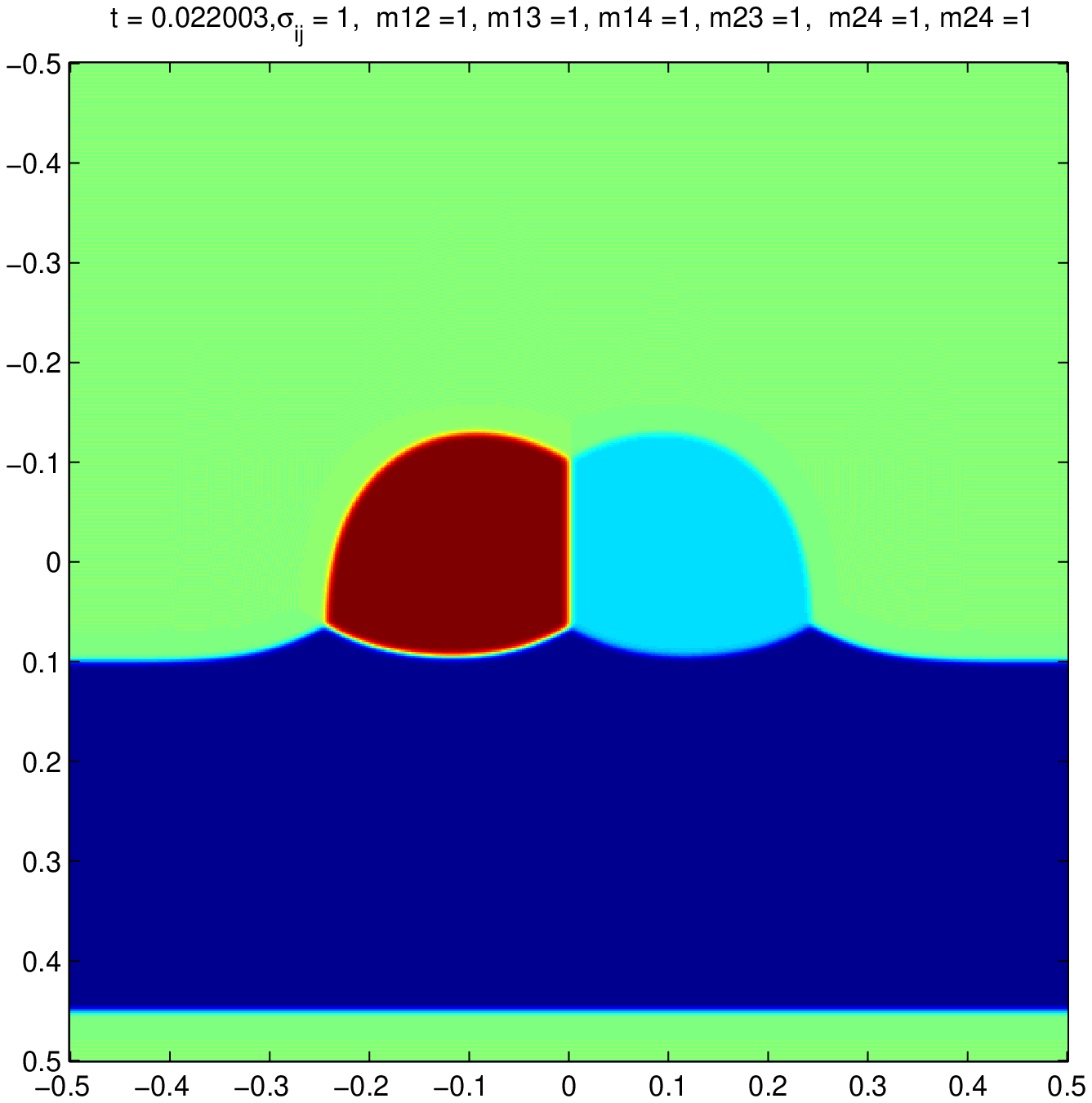}
	\includegraphics[width=3.5cm]{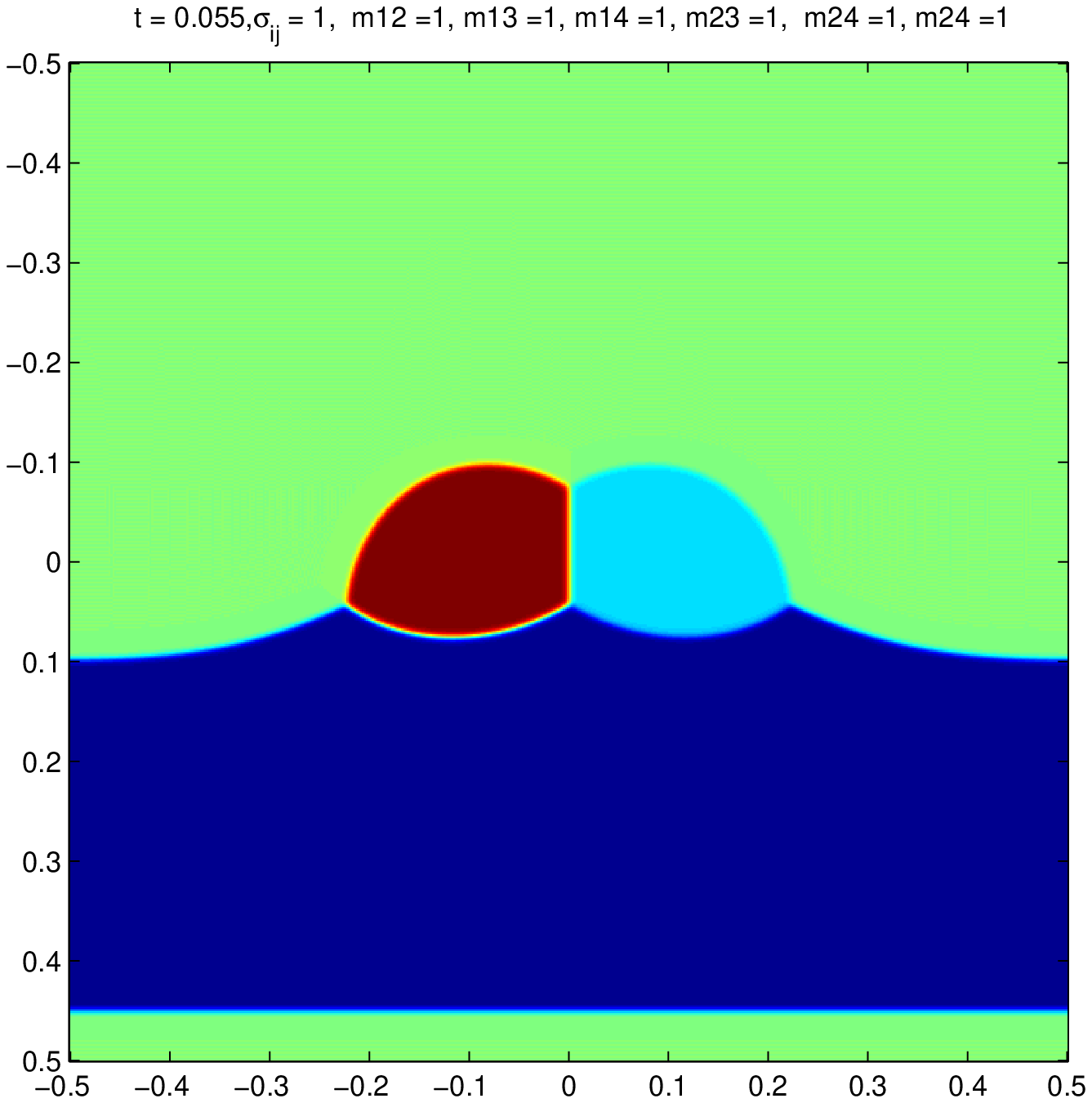}
	\includegraphics[width=3.5cm]{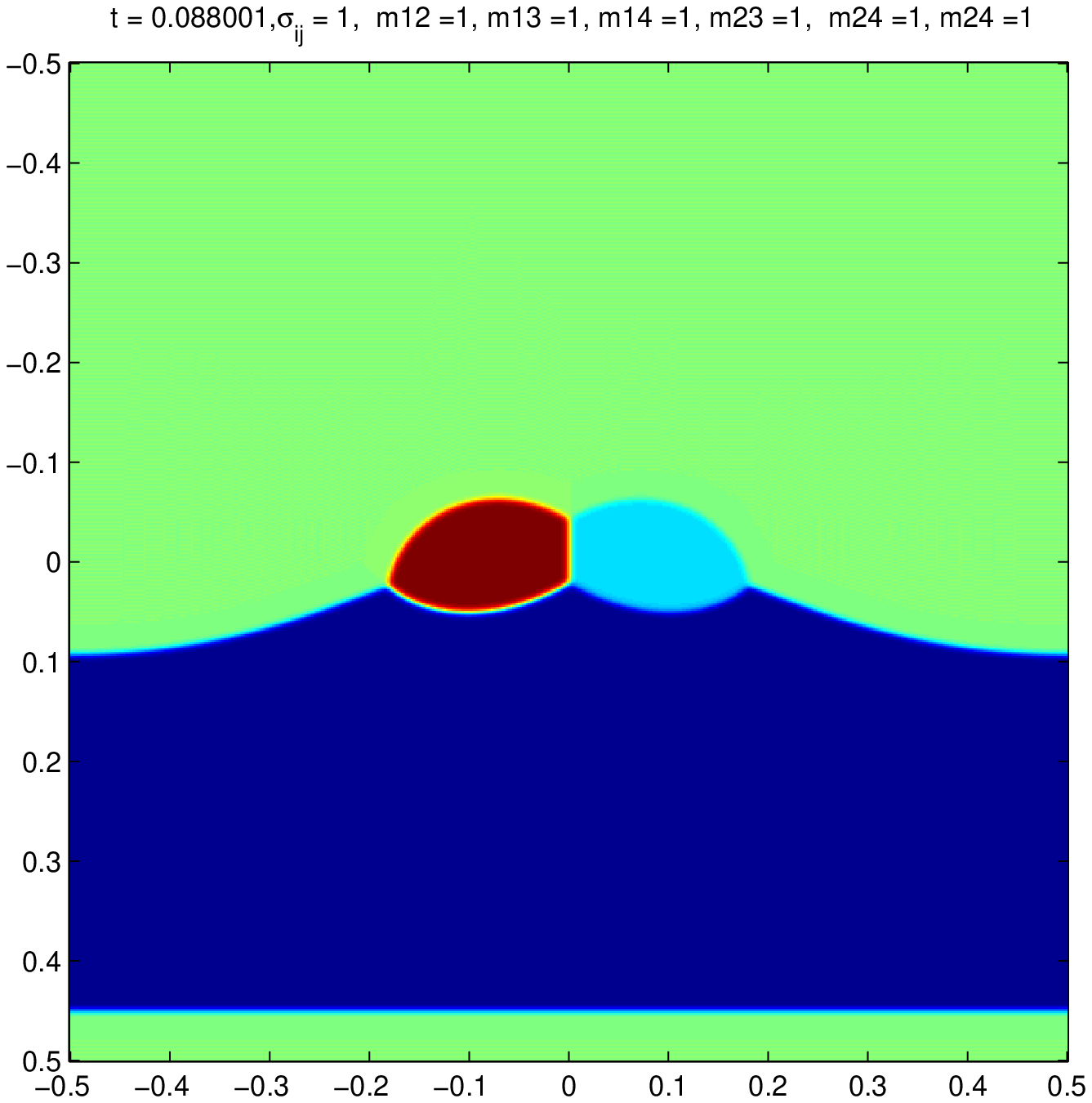} \\
	\includegraphics[width=3.5cm]{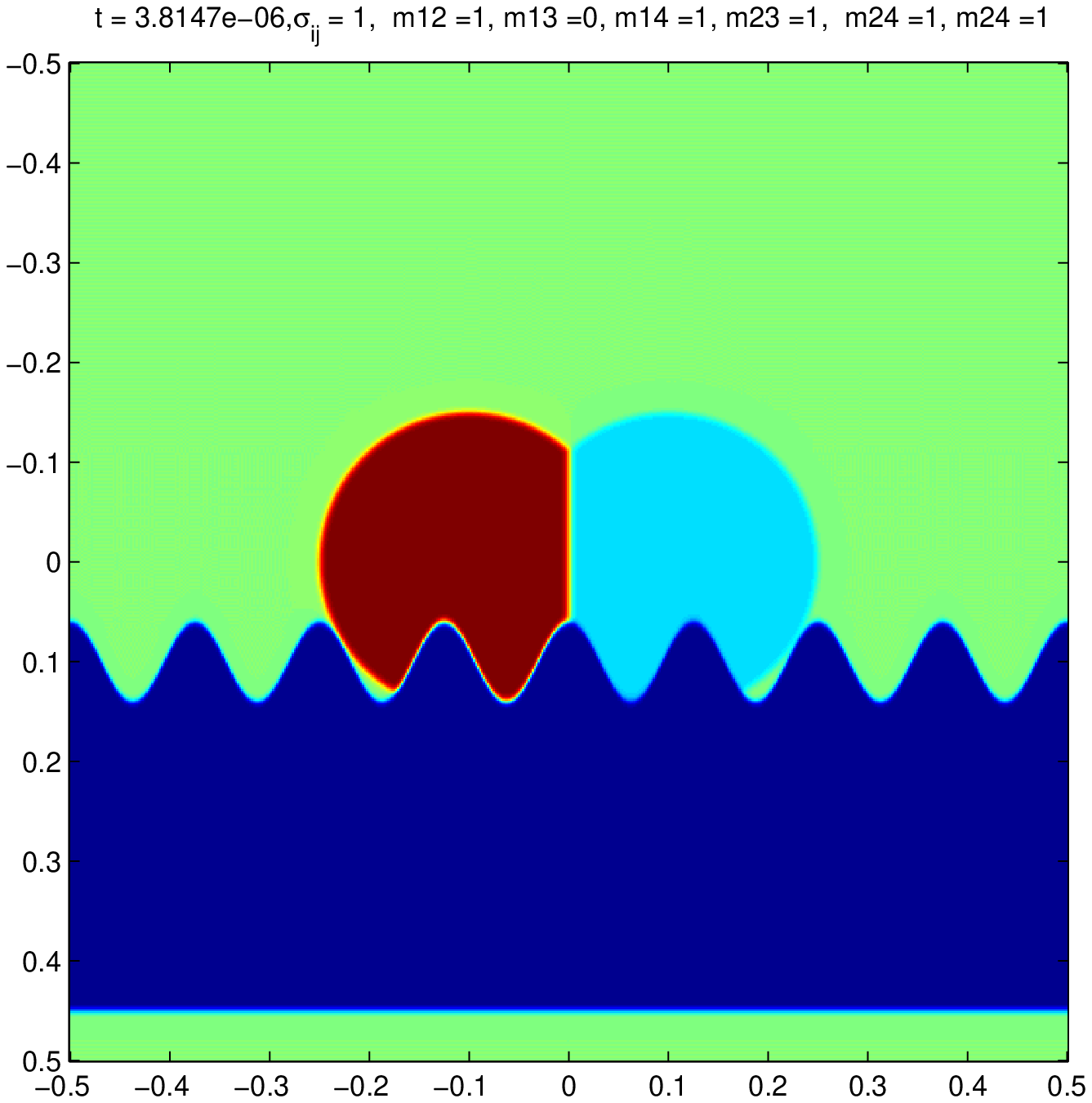}
	\includegraphics[width=3.5cm]{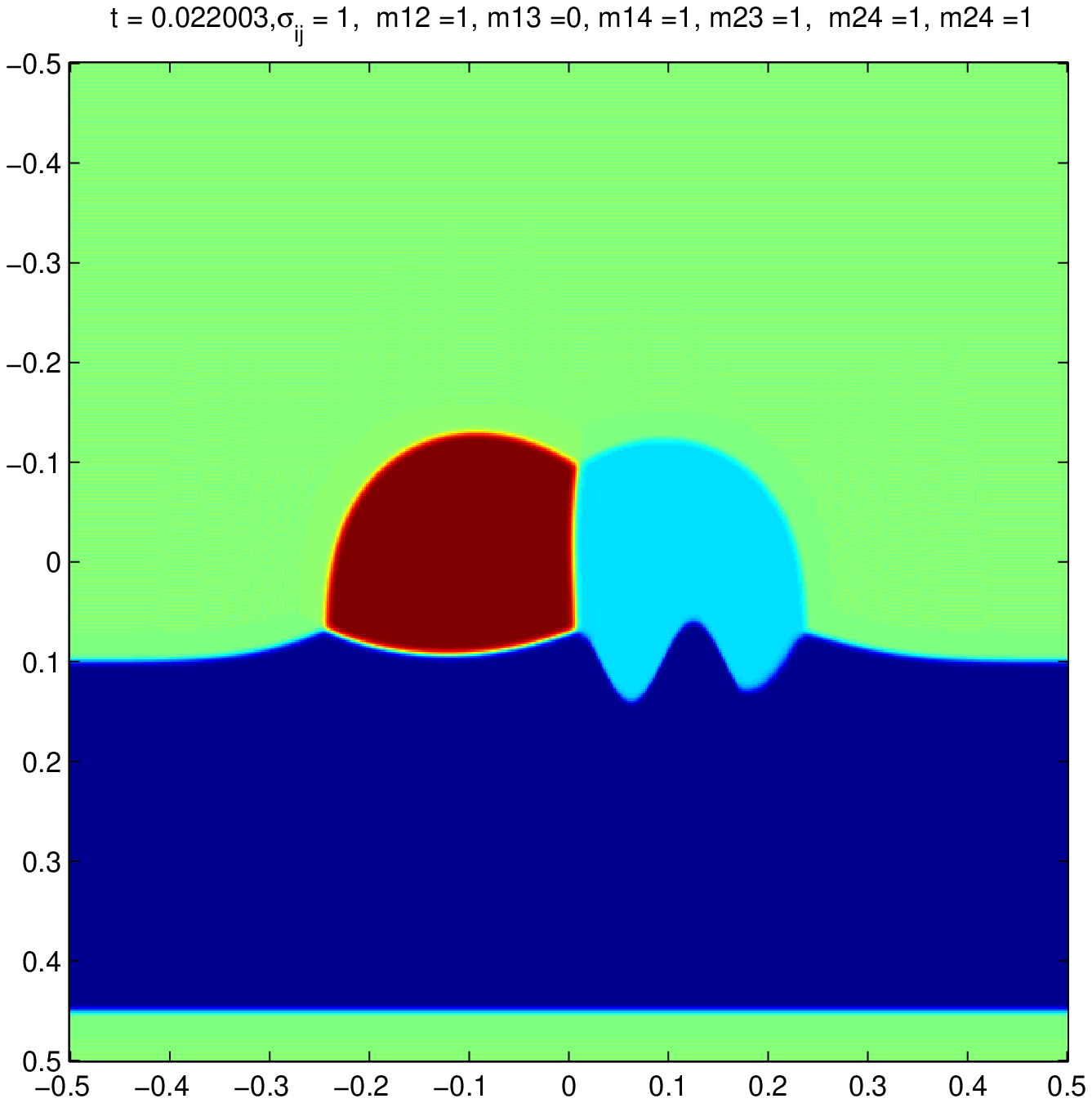}
	\includegraphics[width=3.5cm]{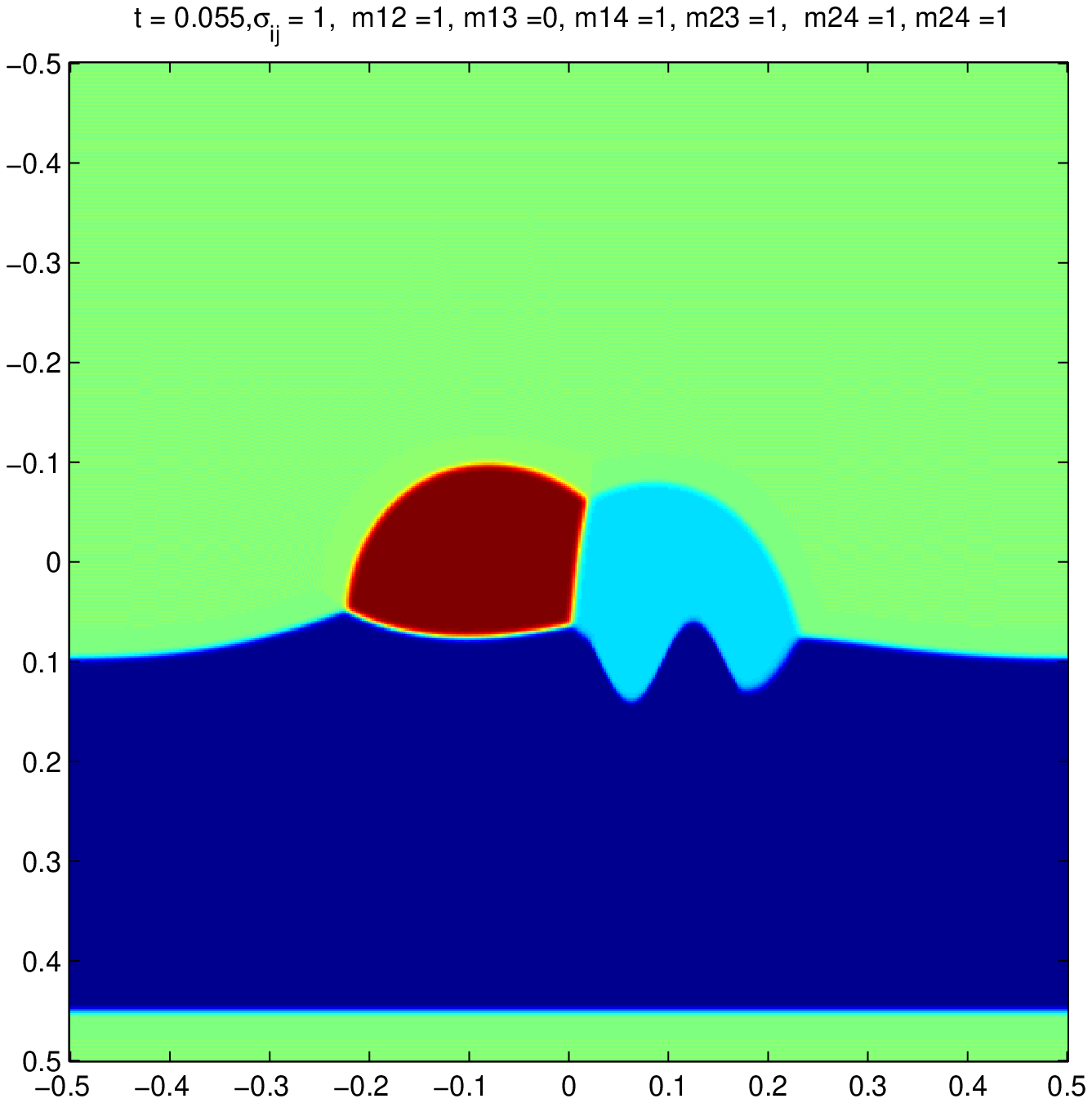}
	\includegraphics[width=3.5cm]{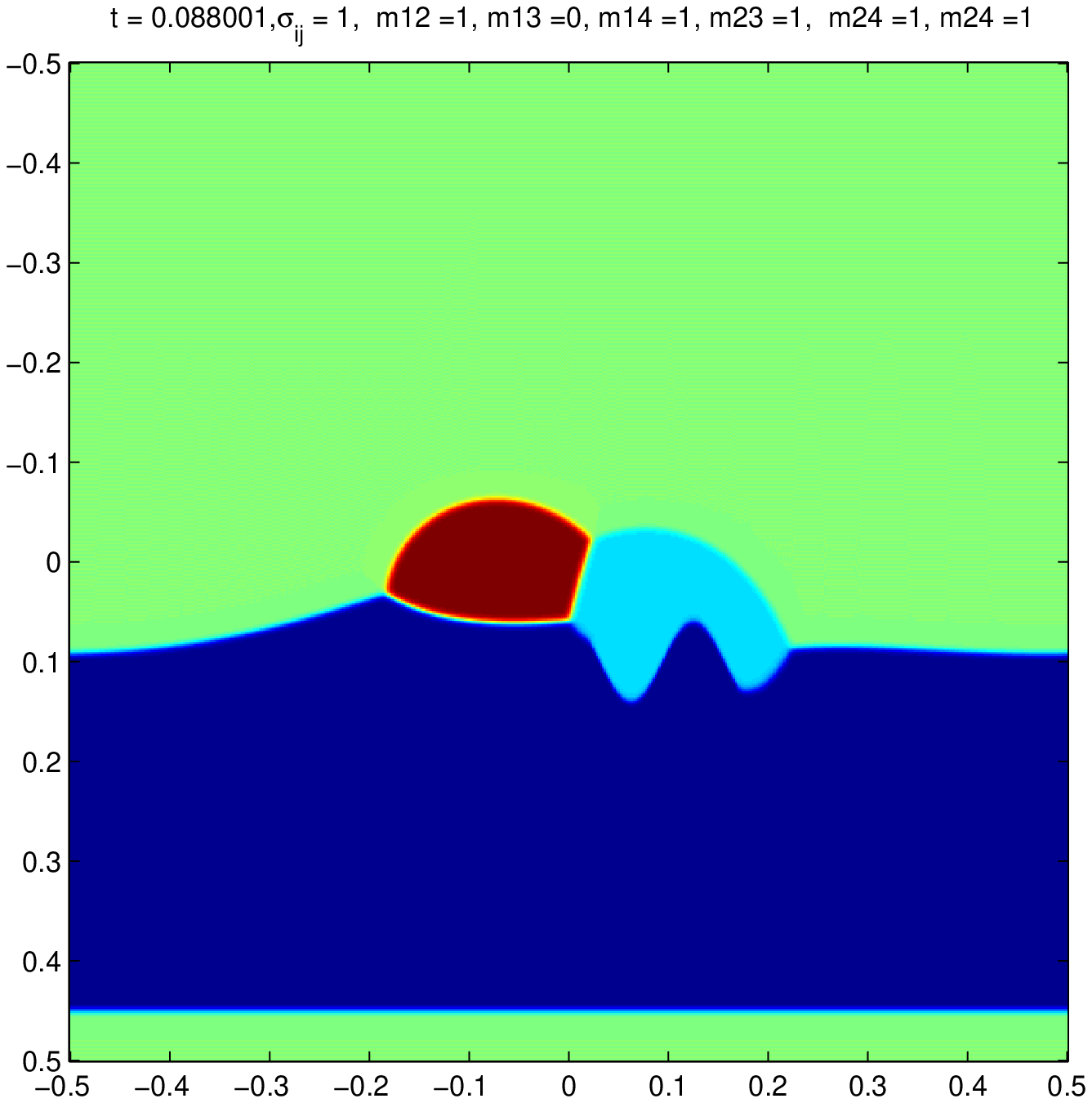} \\
	\includegraphics[width=3.5cm]{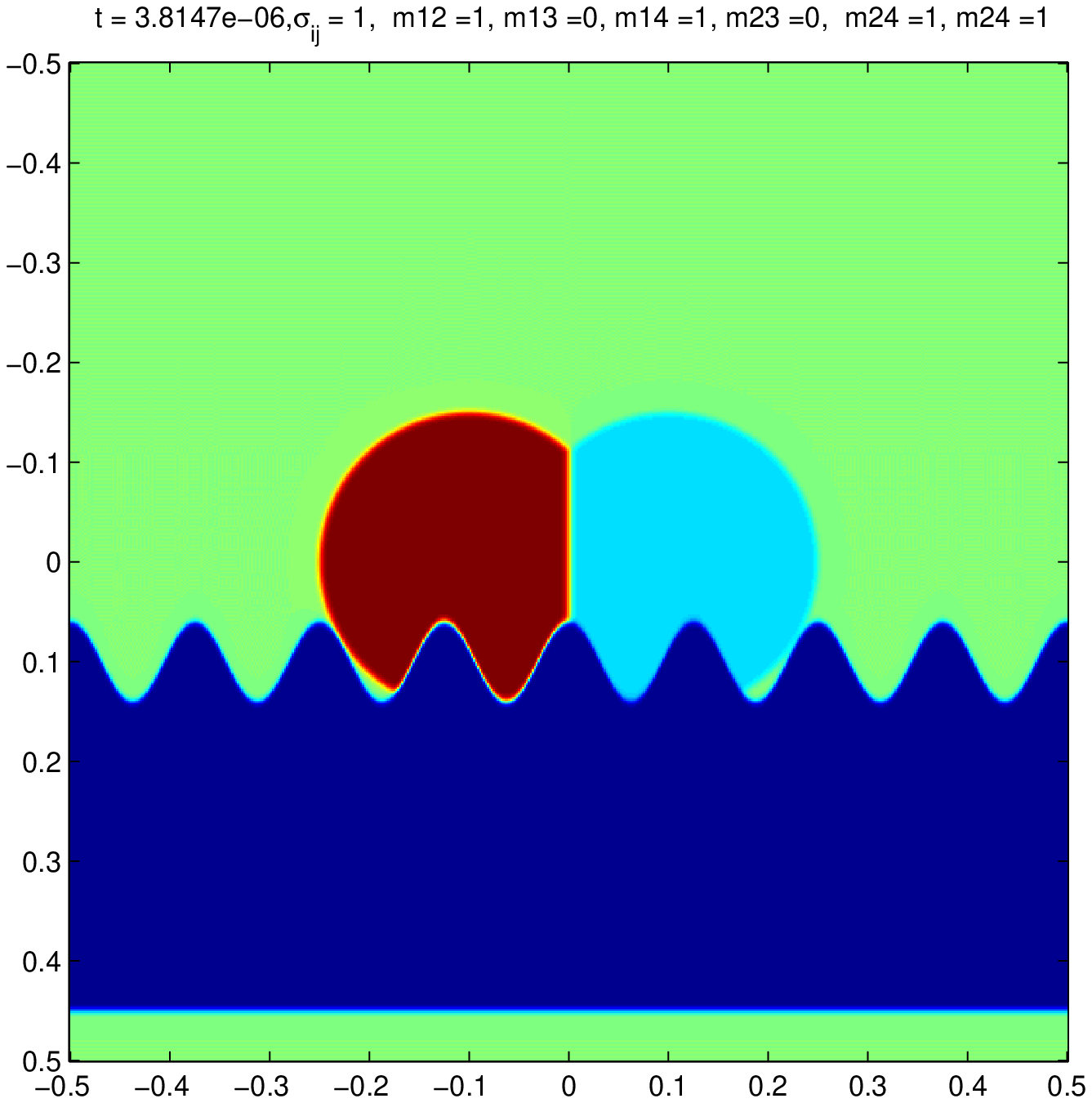}
	\includegraphics[width=3.5cm]{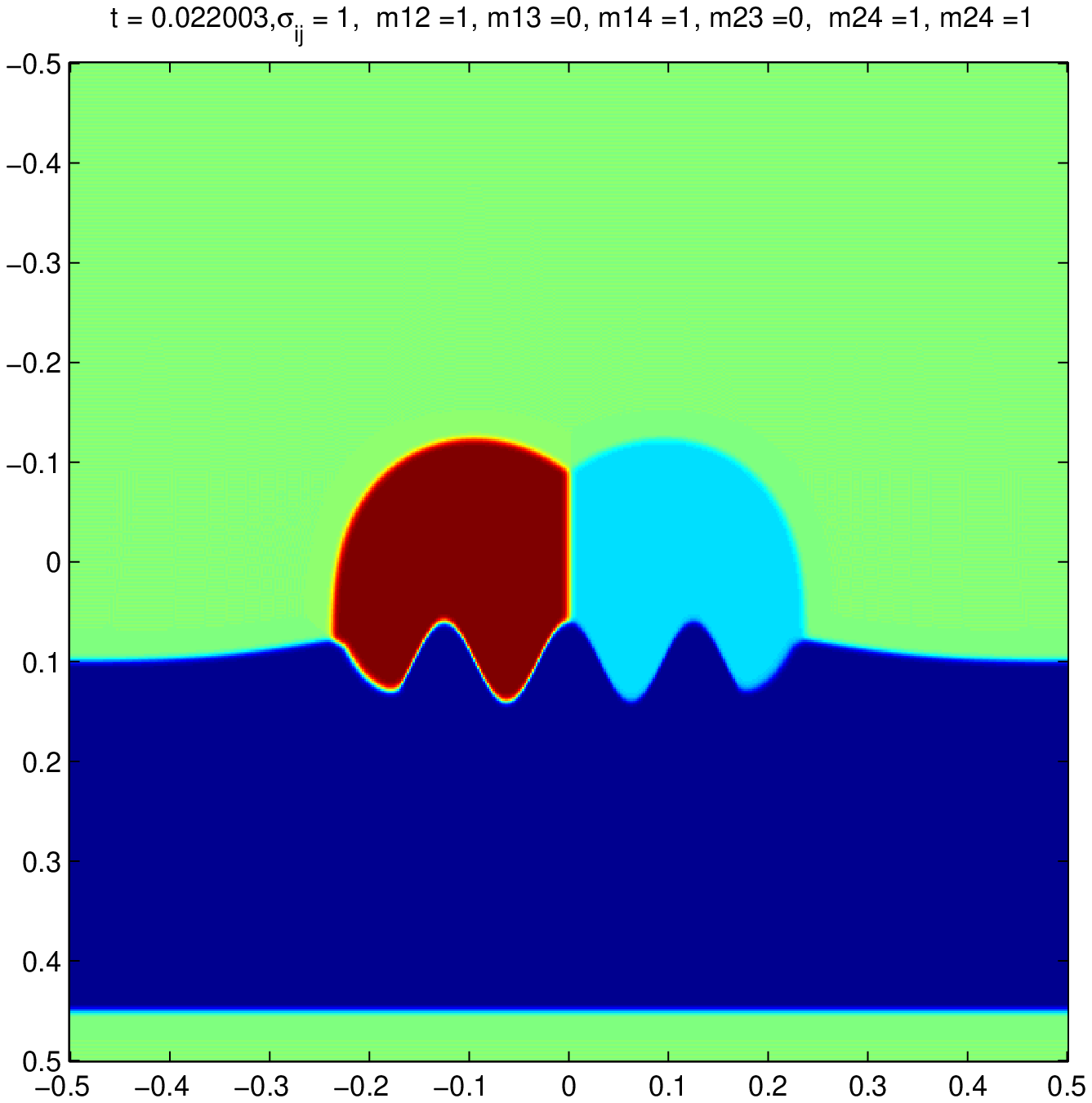}
	\includegraphics[width=3.5cm]{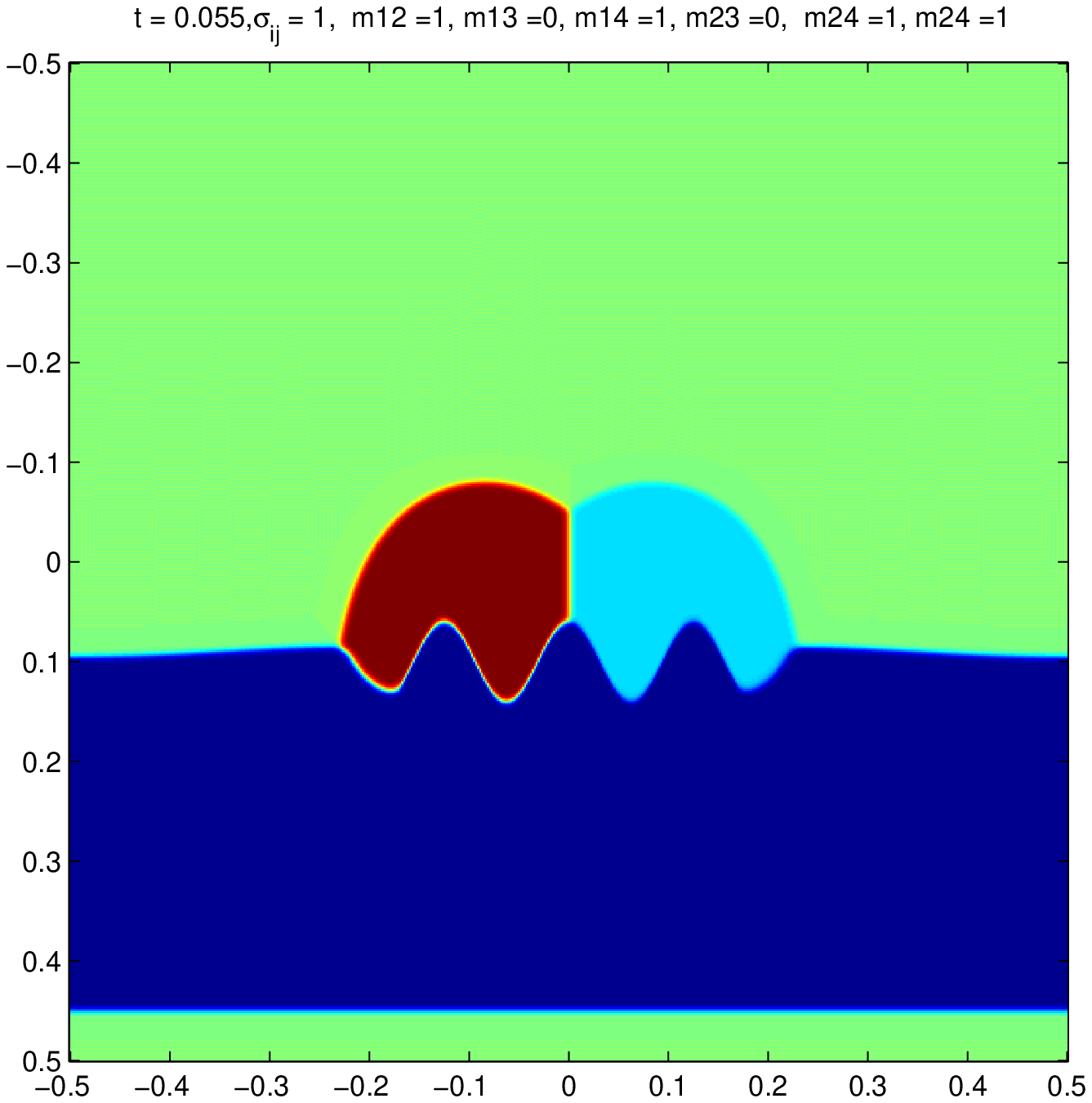}
	\includegraphics[width=3.5cm]{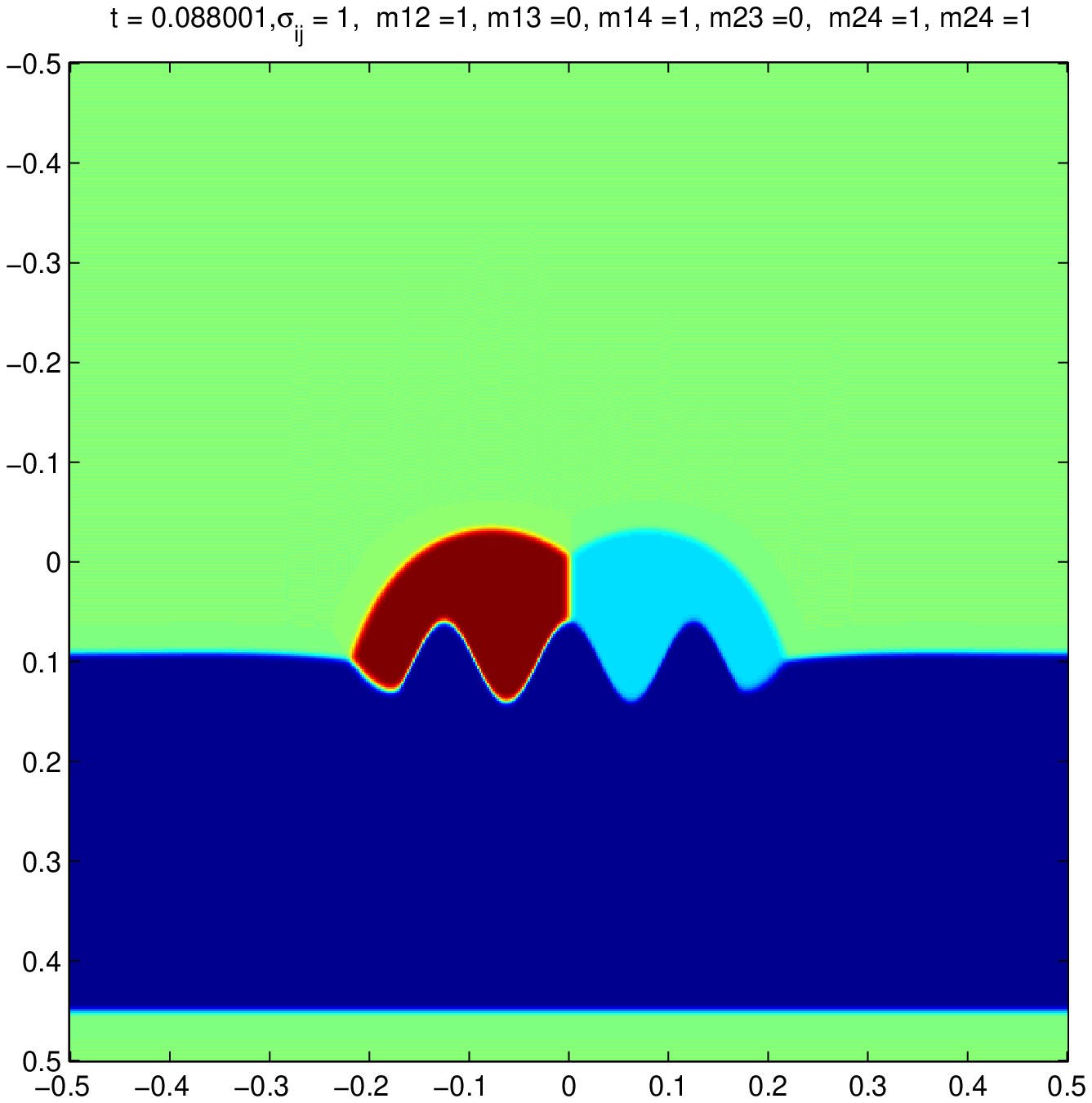} \\	
\caption{Numerical experiments with $N=4$ phases and homogeneous surface tensions $\sigma_{ij} = 1$.
The rows correspond to $(m_{12},m_{13},m_{14},m_{23},m_{24},m_{34})$ equal to 
$(1,1,1,1,1,1)$, $(0,1,1,1,1,1)$, and  $(0,0,1,1,1,1)$, respectively. 
The  {images} represent the values of $u_1 + 3 u_2 + 1.5 u_4$ at different times, {with a suitable colormap so that} $u_1$, $u_2$, $u_3$, and $u_4$ are shown in light blue,  red, blue and green, respectively.}
\label{fig_N4}
\end{figure}

We show now that our method can handle flows involving more than 3 phases. We consider a configuration with 4 phases and a canonical decomposition of the mobilities 
$ {\bf m} =  (m_{12},m_{13},m_{14},m_{23},m_{24},m_{34})$, which takes the form 
 $$(m_{12},m_{13},m_{14},m_{23},m_{24},m_{34}) = \sum_{p=1}^{6} (m^p_{12},m^p_{13},m^p_{14},m^p_{23},m^p_{24},m^p_{34}),$$
 where
 $$ (m^p_{12},m^p_{13},m^p_{14},m^p_{23},m^p_{24},m^p_{34}) = \begin{cases}
                                                        (m_{12},0,0,0,0,0,0) &\text{ if } p = 1\\
                                                         (0,m_{13},0,0,0,0,0) &\text{ if } p = 2\\
                                                           \hspace*{15mm}\vdots &  \\
                                                         (0,0,0,0,0,0,m_{34}) &\text{ if } p = 6\\  
                                                              \end{cases}
 $$ 
Figure \ref{fig_N4} shows a series of numerical experiments using 
$$ (\sigma_{12},\sigma_{13},\sigma_{14},\sigma_{23},\sigma_{24},\sigma_{34}) = (1,1,1,1,1,1).$$
The rows correspond to 
${\bf m} = (1,1,1,1,1,1)$, ${\bf m} = (0,1,1,1,1,1)$, and  ${\bf m} = (0,0,1,1,1,1)$, respectively.  
In each image, the phases $u_1, u_2, u_3, u_4$ are depicted in
light blue, red, blue, and green, respectively.
These results show good agreement with the expected theoretical flows.

%


\subsection{Numerical experiments in dimension $3$} 

\begin{figure}[htbp]
\centering
	\includegraphics[width=5.2cm]{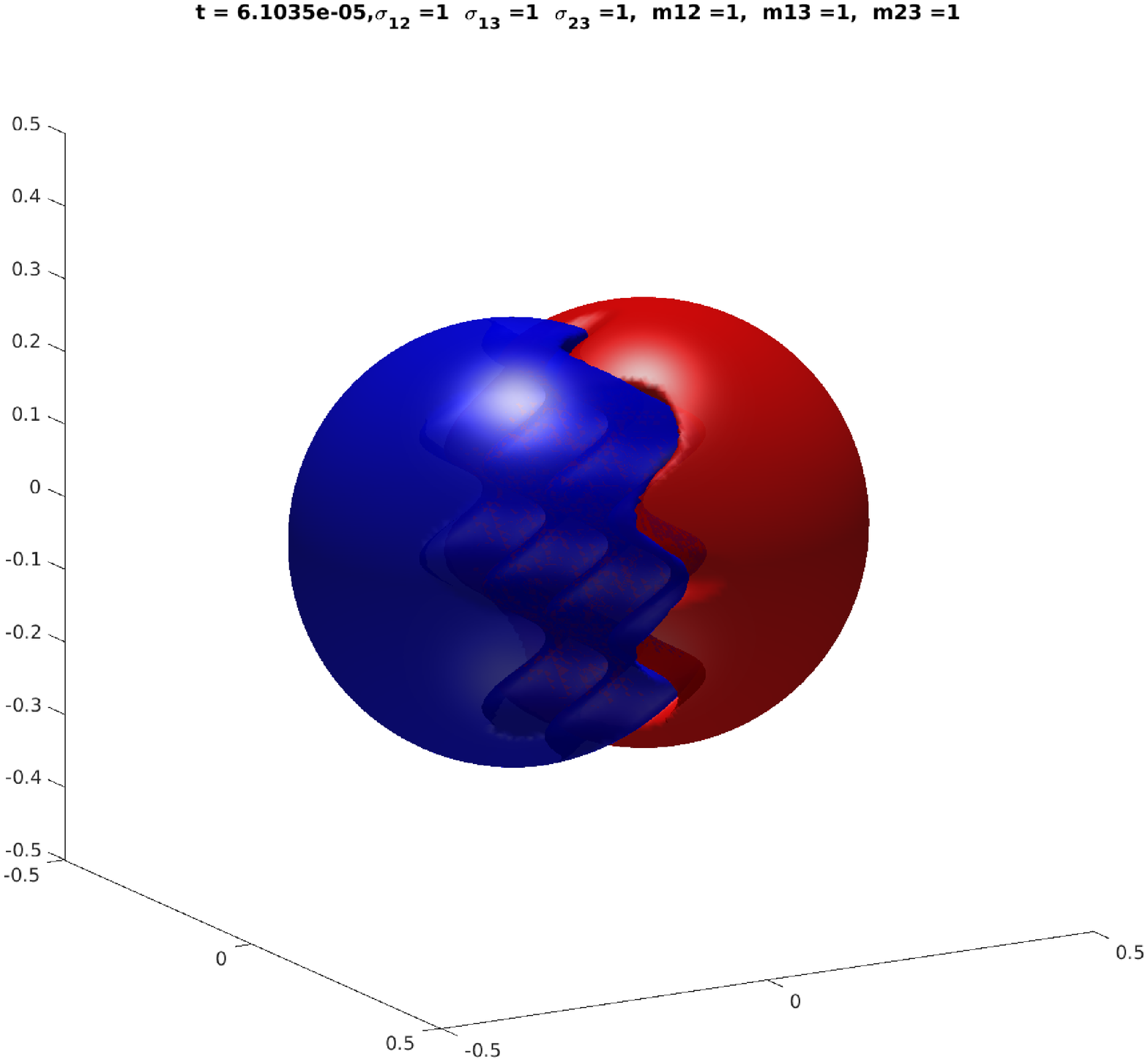}
	\includegraphics[width=5.2cm]{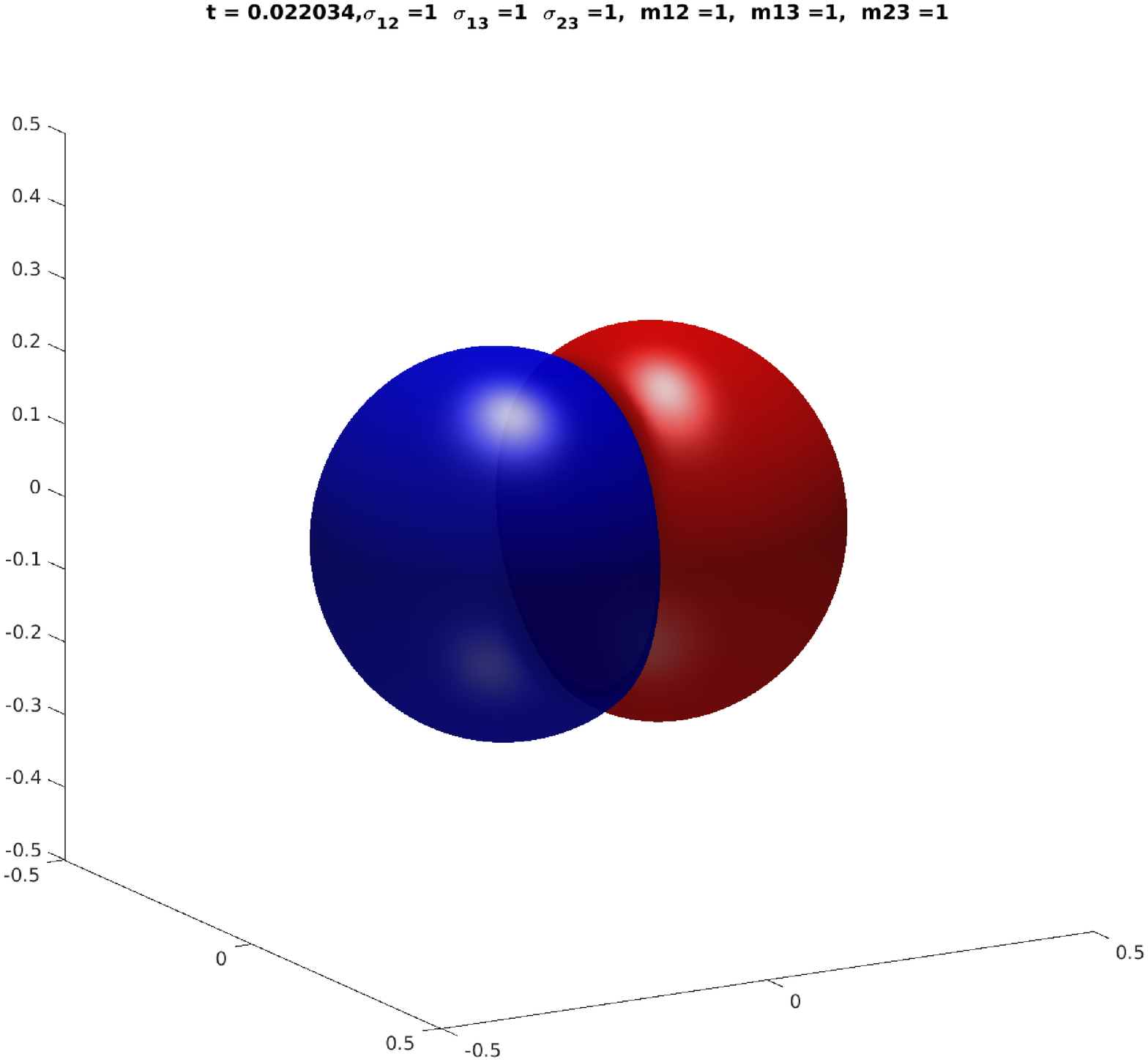}
	\includegraphics[width=5.2cm]{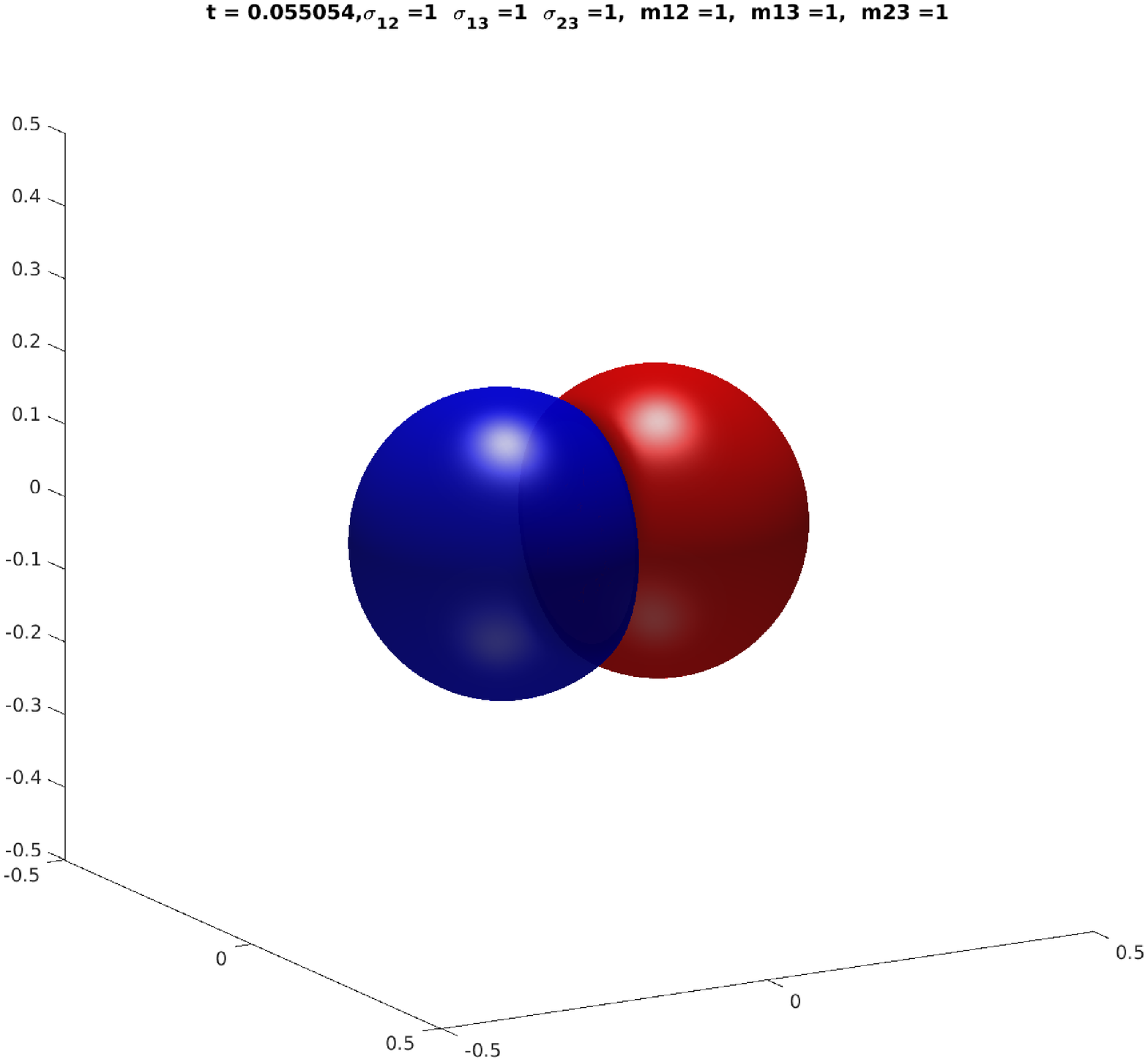} \\
	\includegraphics[width=5.2cm]{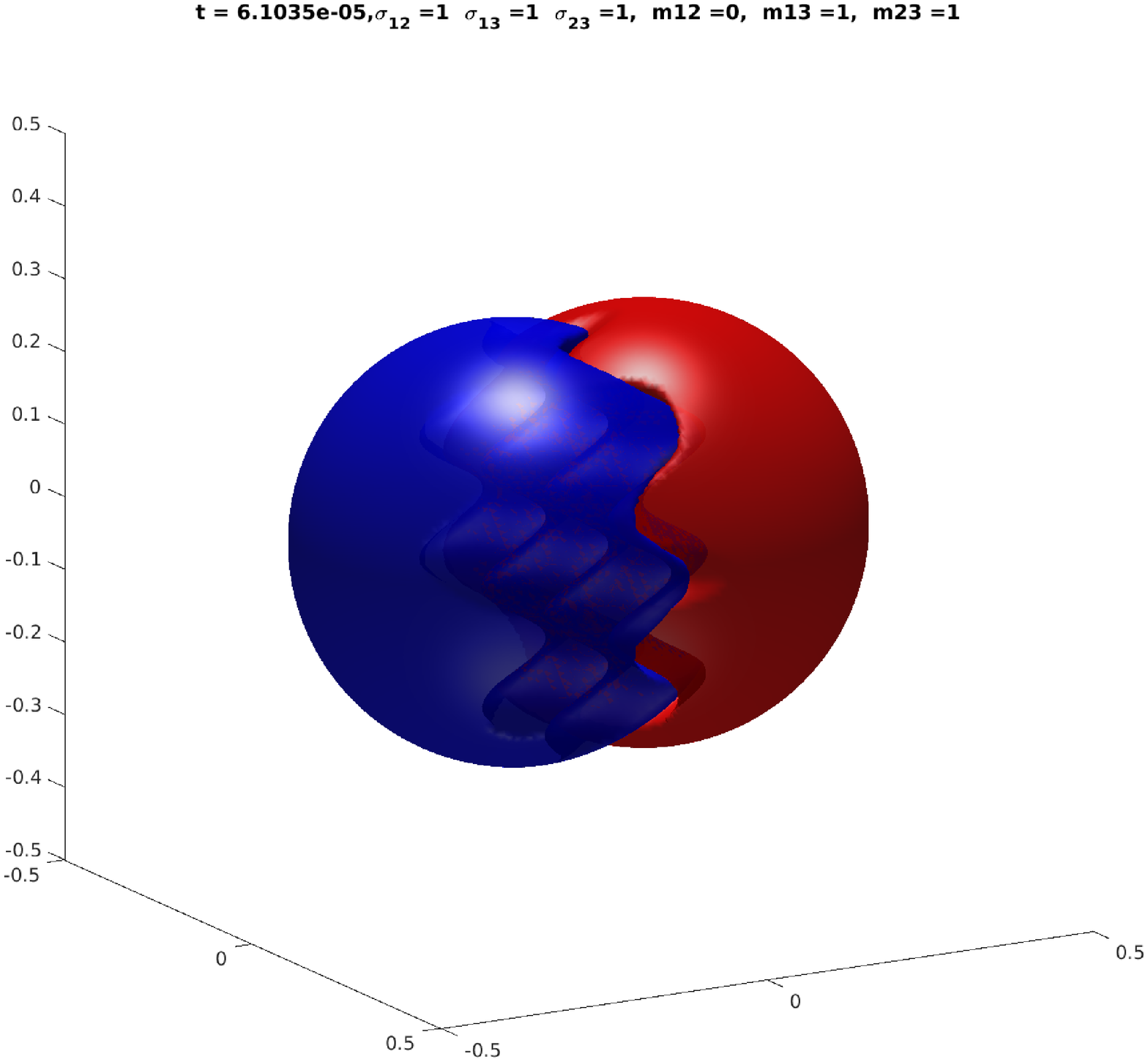}
	\includegraphics[width=5.2cm]{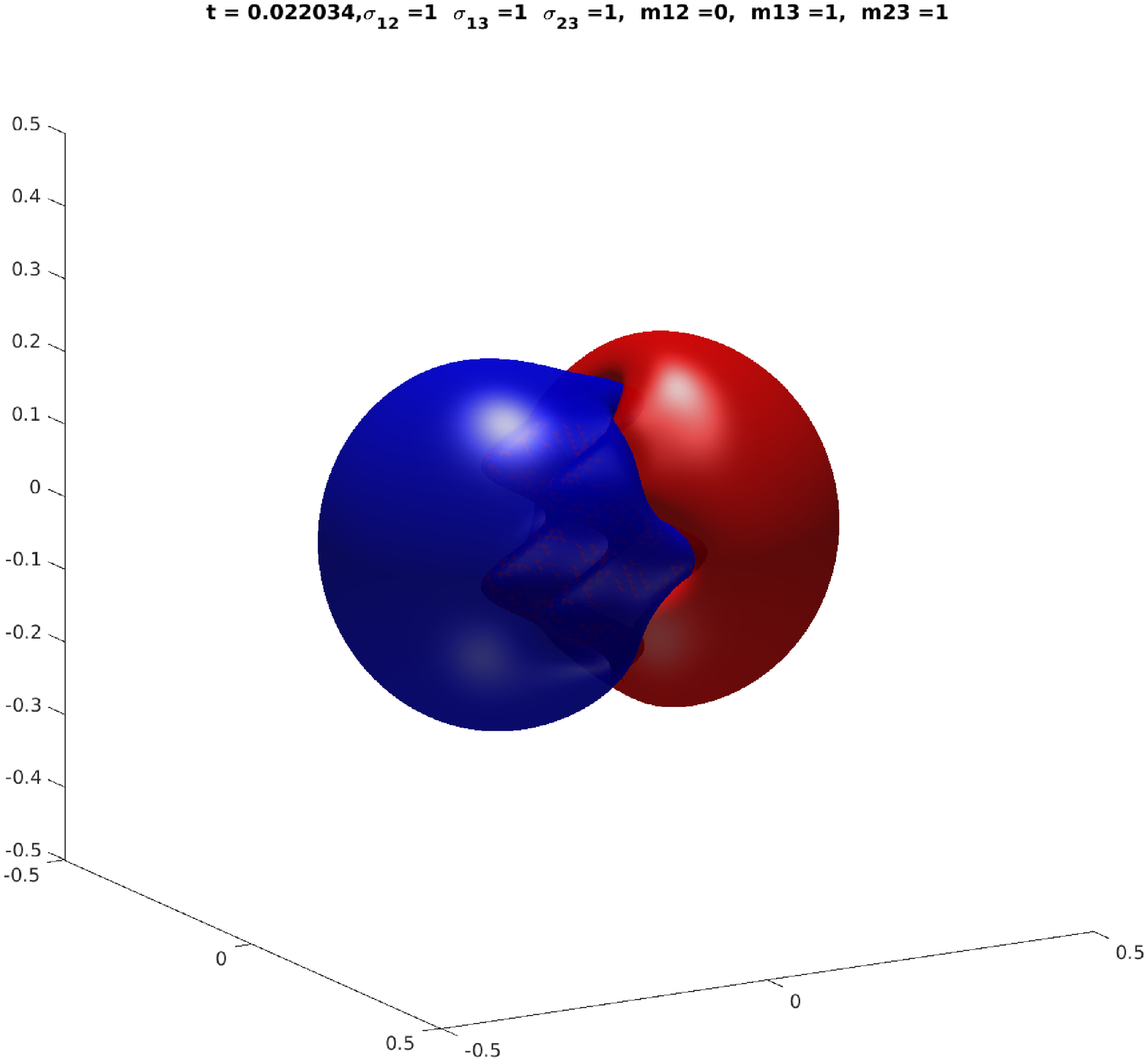}
	\includegraphics[width=5.2cm]{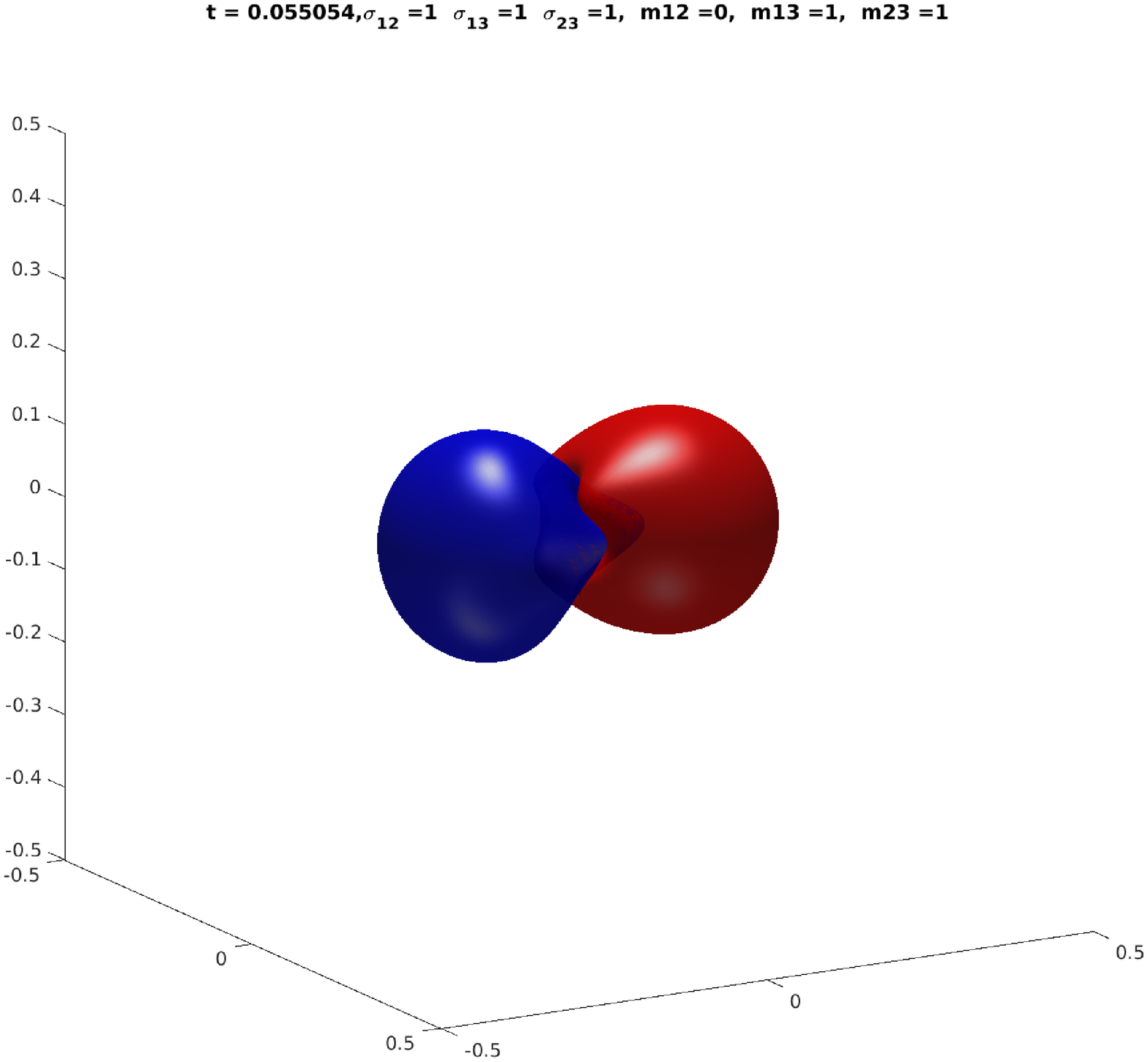} \\
	\includegraphics[width=5.2cm]{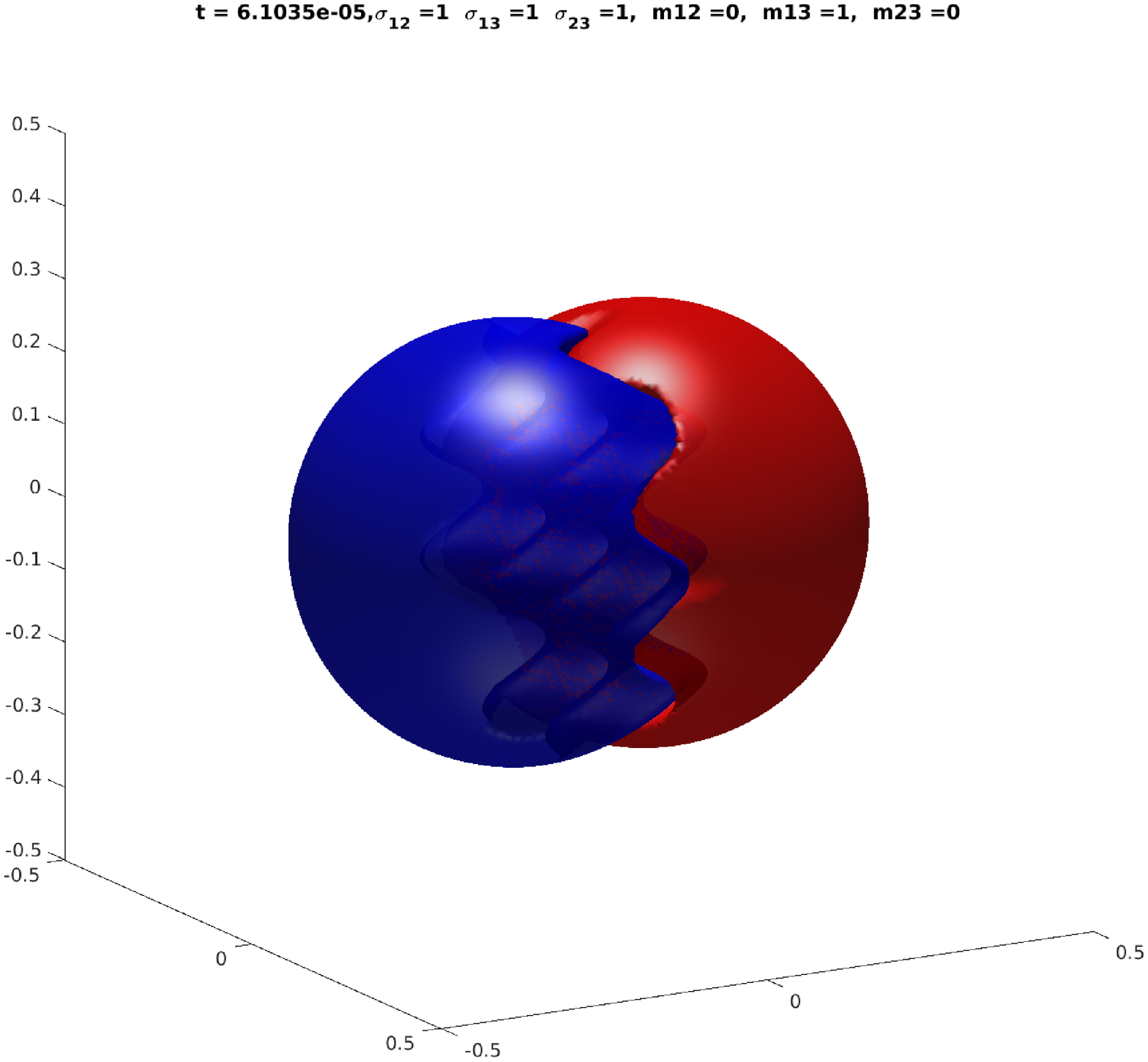}
	\includegraphics[width=5.2cm]{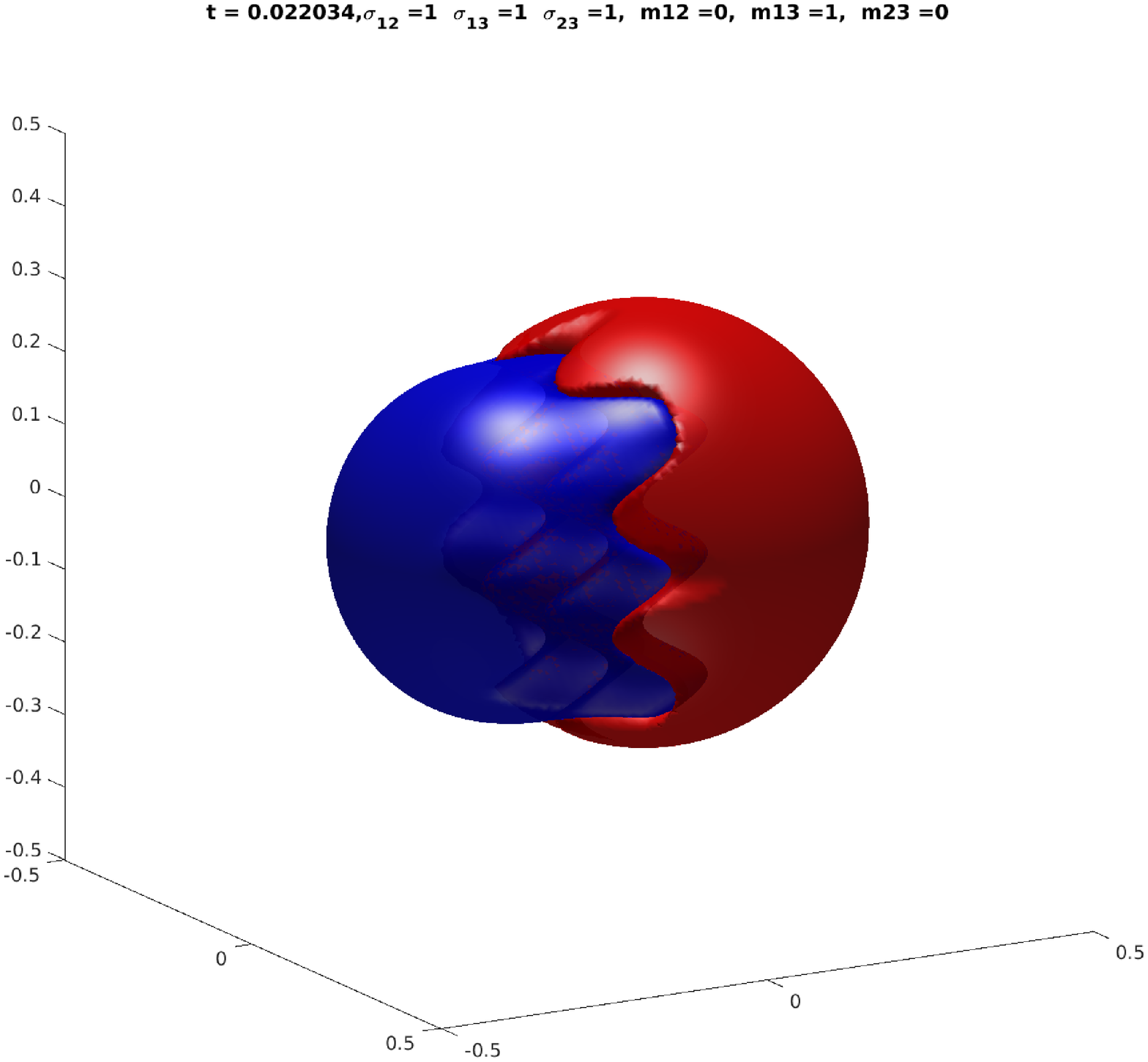}
	\includegraphics[width=5.2cm]{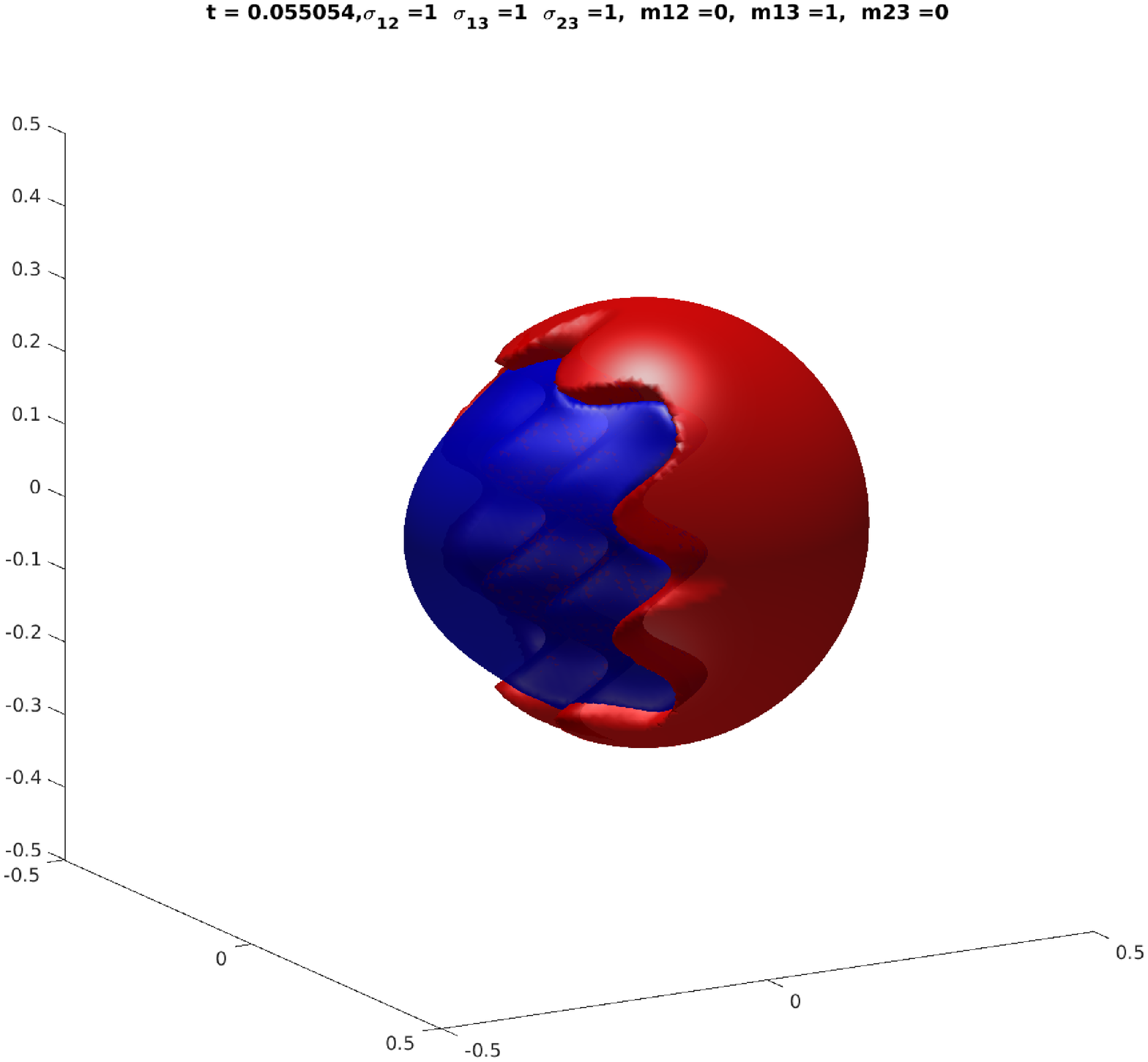} \\
\caption{3D multiphase mean curvature flows with homogeneous surface tensions $\sigma_{12} = \sigma_{13} = \sigma_{23} = 1$.
The rows correspond to $(m_{12},m_{13},m_{23}) = (1,1,1)$, $(m_{12},m_{13},m_{23}) = (0,1,1)$, and   $(m_{12},m_{13},m_{23}) = (0,1,0)$, respectively.
{The images show the level set $\{u_1 = 1/2\}$ in blue,
and the level set $\{u_2=1/2\}$ in red, at different times}.
}
\label{fig_3D1}
\end{figure}

Figure \ref{fig_3D1} shows the 3D version of the 2D computations reported in
Figure~\ref{fig_compar1}. 
The surface tensions are identical, $\sigma_{ij} = 1$.
The rows represent the evolutions from the same initial condition with 
mobilities $(m_{12},m_{13},m_{23})$ equal to $(1,1,1)$, $(0,1,1)$, and  $(0,1,0)$
respectively. 
In each image, the phases $u_1$ and $u_2$ are depicted in blue and red, respectively.  \\

Our last example, shown in Figure \ref{fig_3D1},  concerns a more complex situation with 3 phases
where the initial geometry represents  a toy truck.
We compare evolutions obtained with different sets of mobilities,
and with surface tensions $\sigma_{i,j}$ all equal to $1$.

 \begin{figure}[htbp]
\centering
	\includegraphics[width=5.2cm]{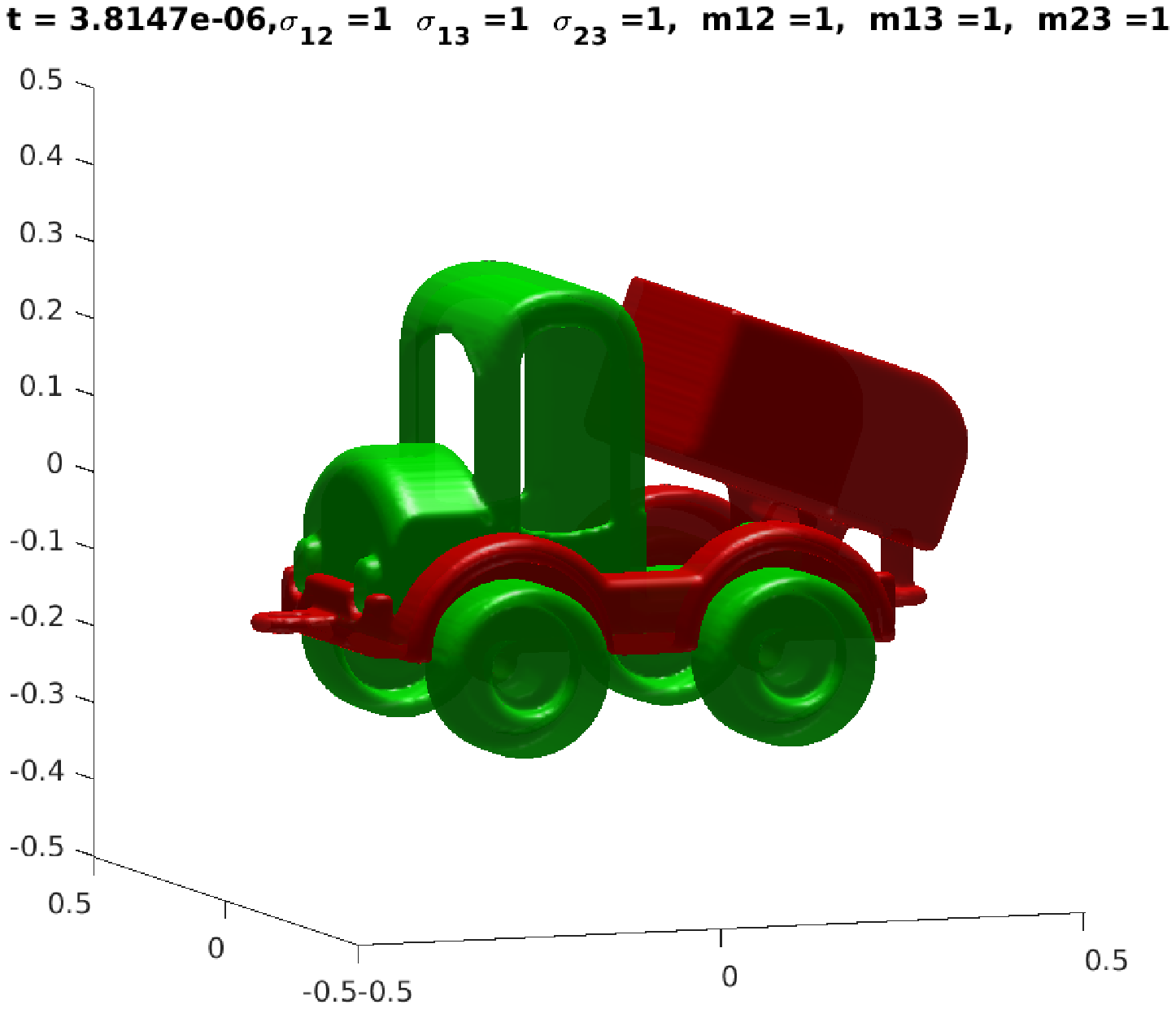}
	\includegraphics[width=5.2cm]{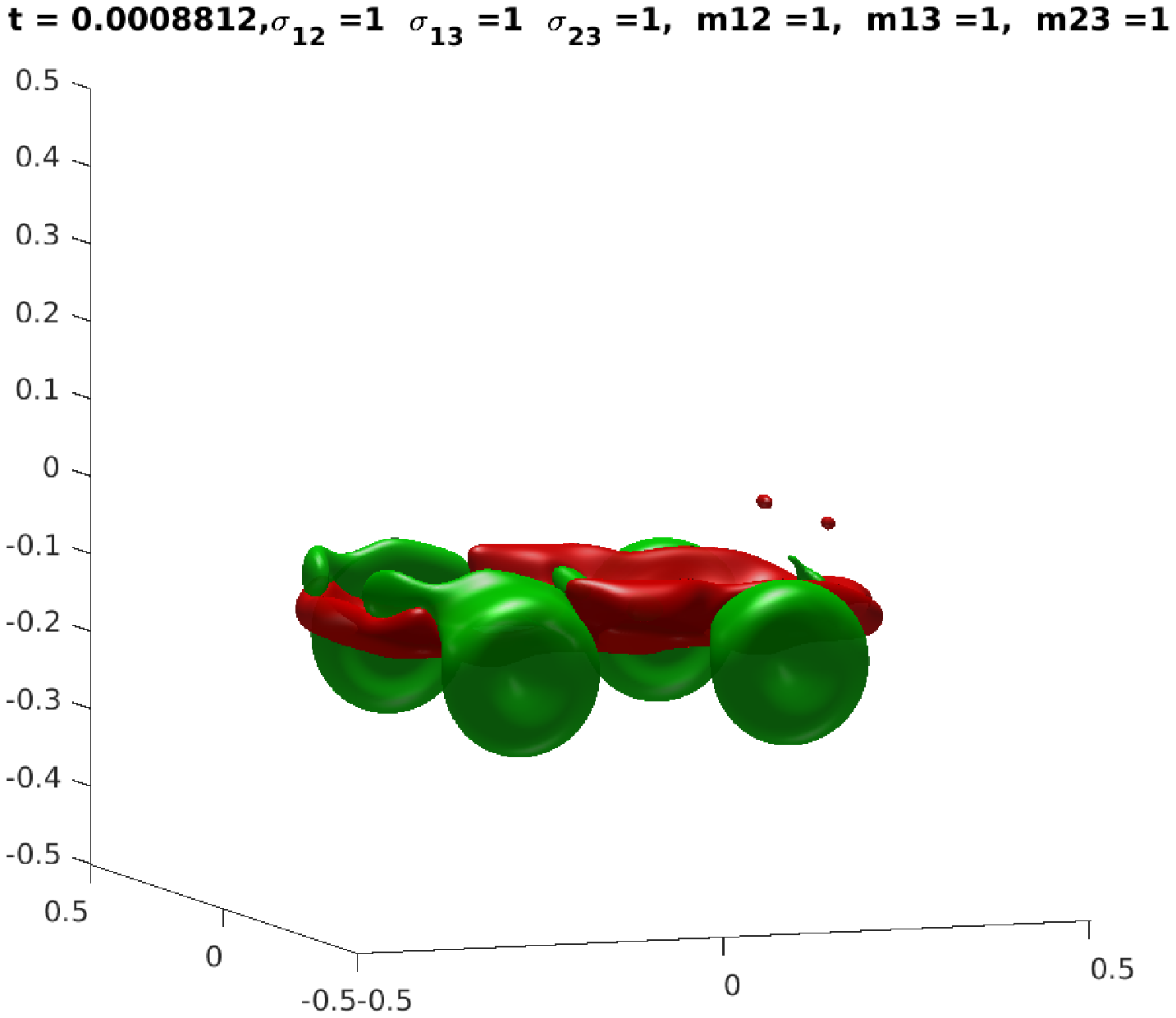}
	\includegraphics[width=5.2cm]{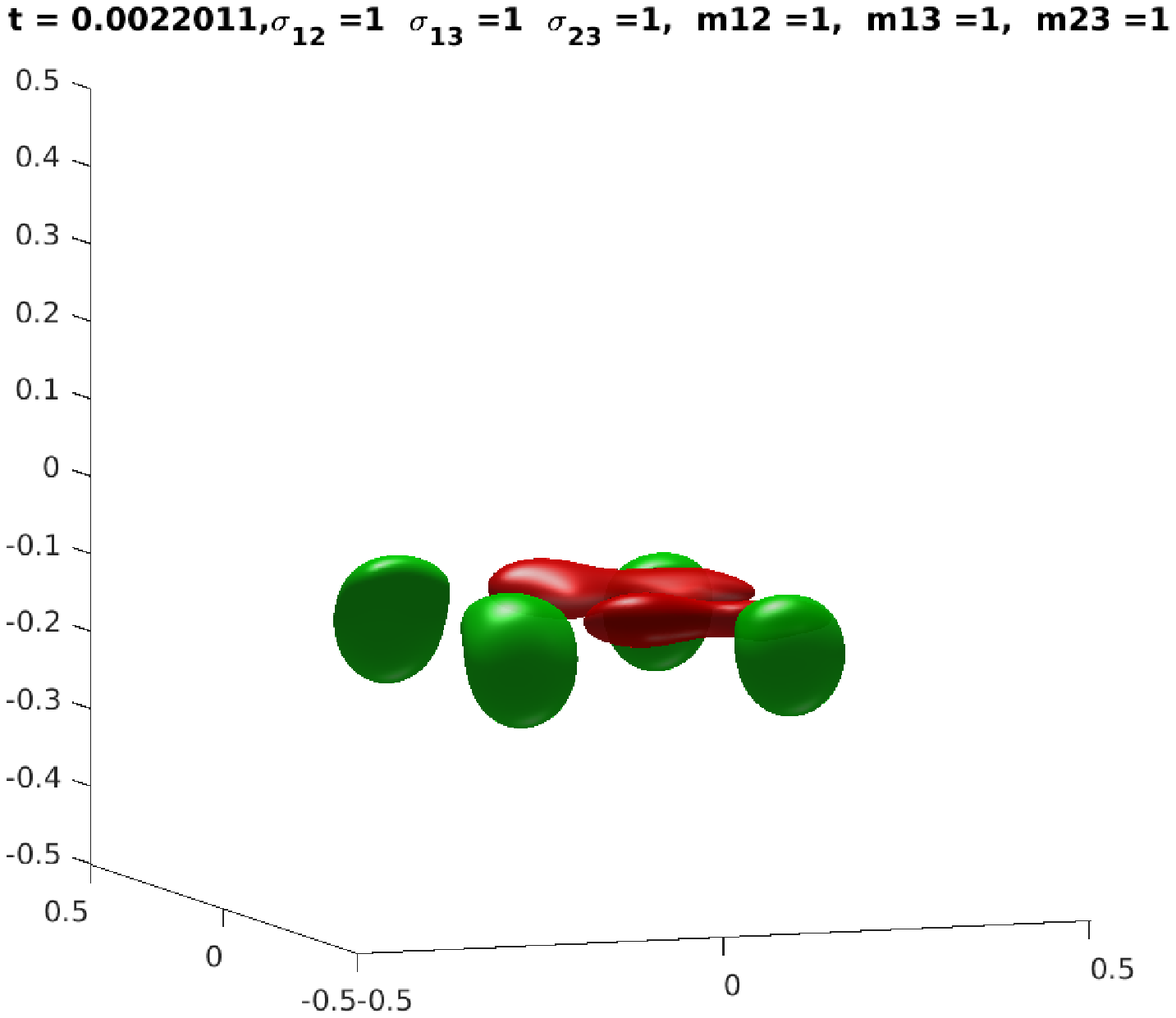} \\
	\includegraphics[width=5.2cm]{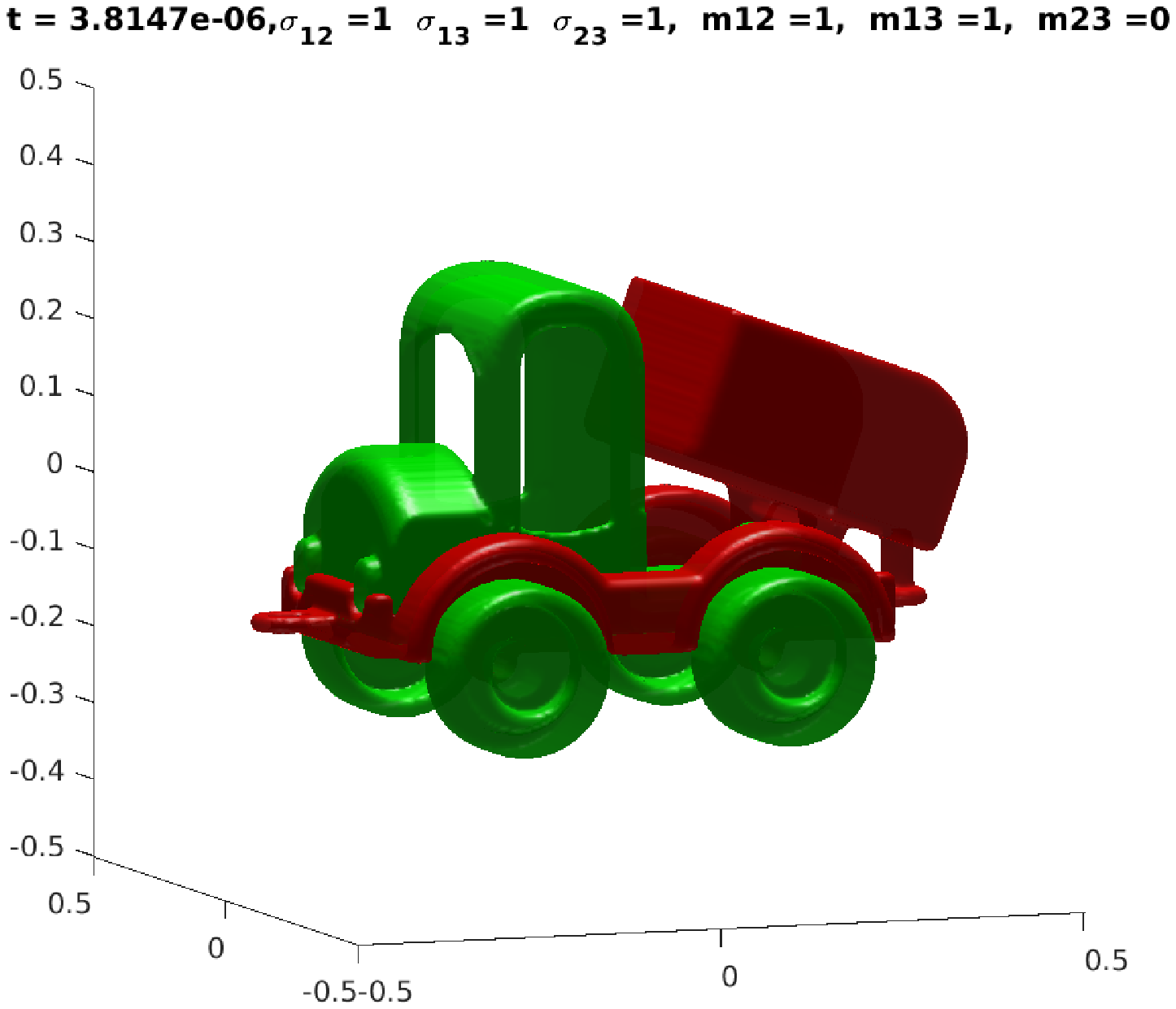}
	\includegraphics[width=5.2cm]{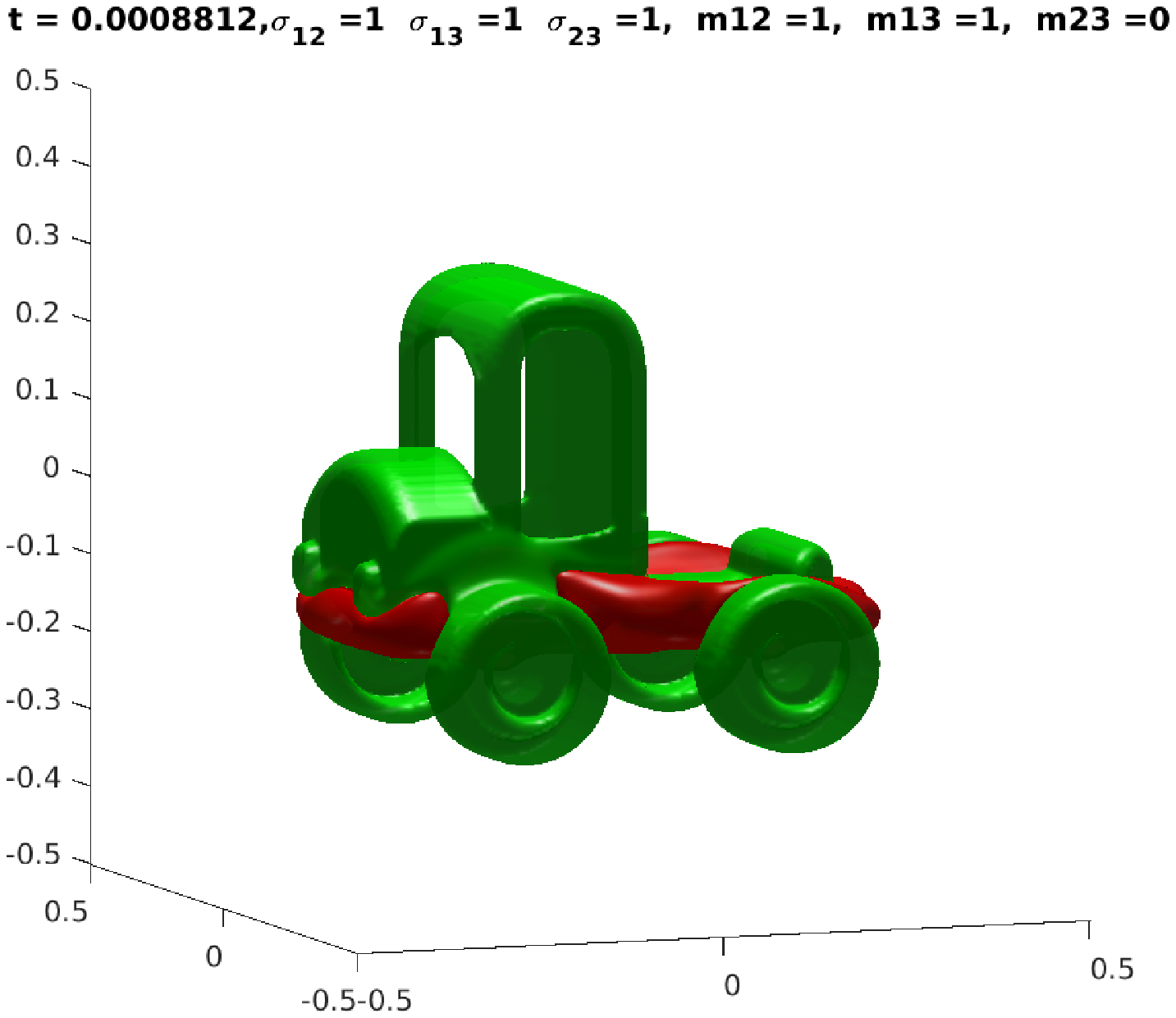}
	\includegraphics[width=5.2cm]{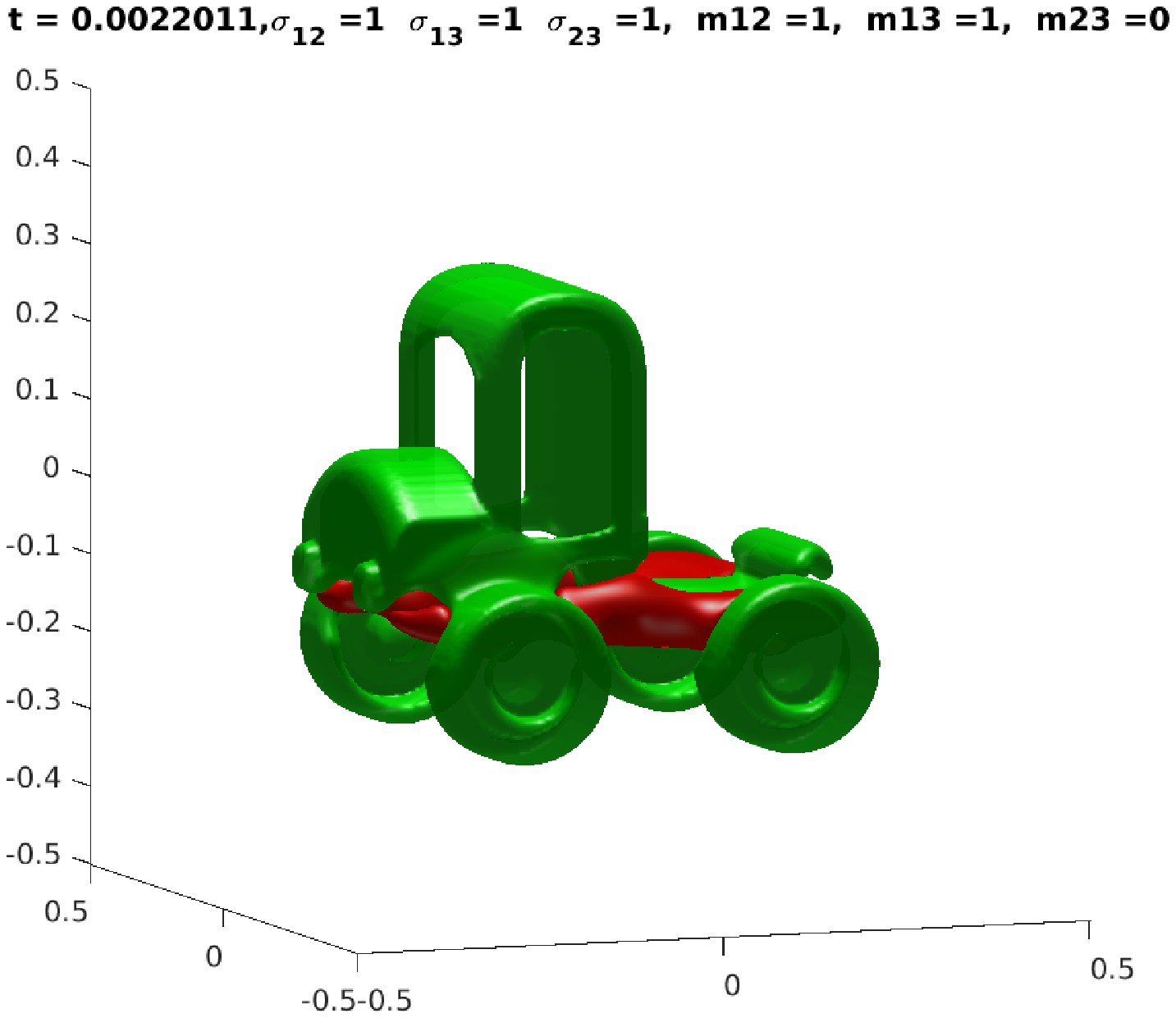} \\
	\includegraphics[width=5.2cm]{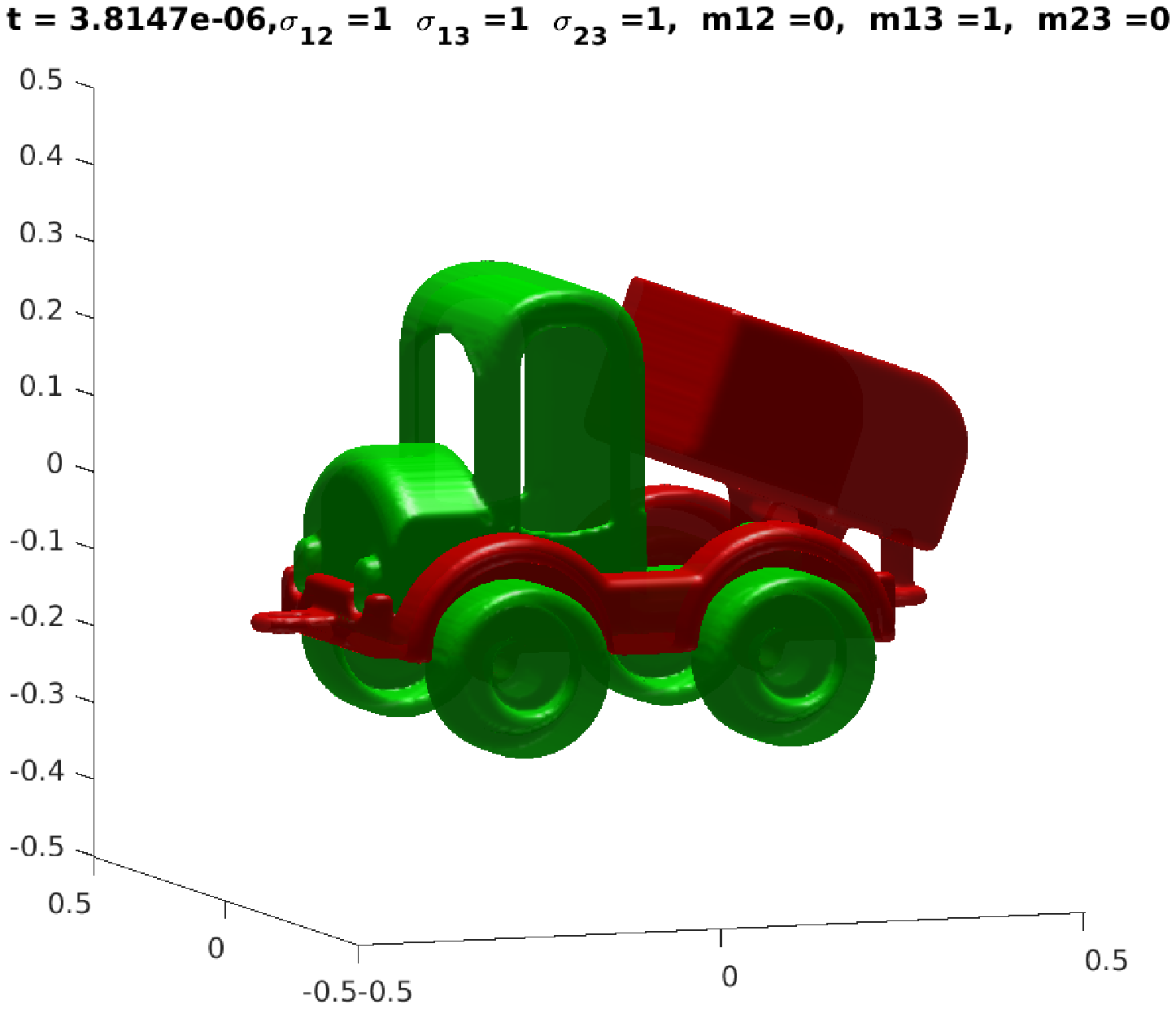}
	\includegraphics[width=5.2cm]{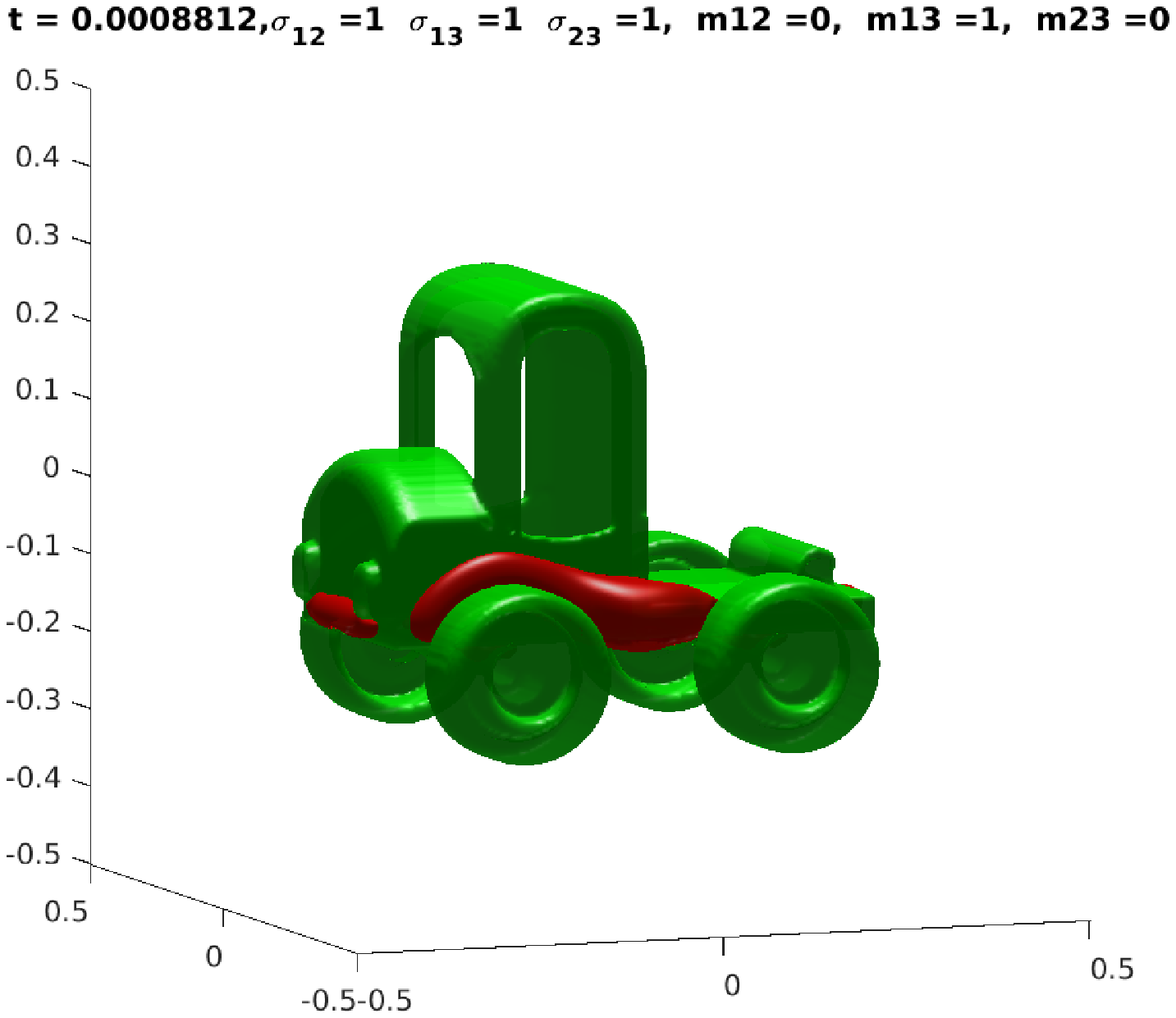}
	\includegraphics[width=5.2cm]{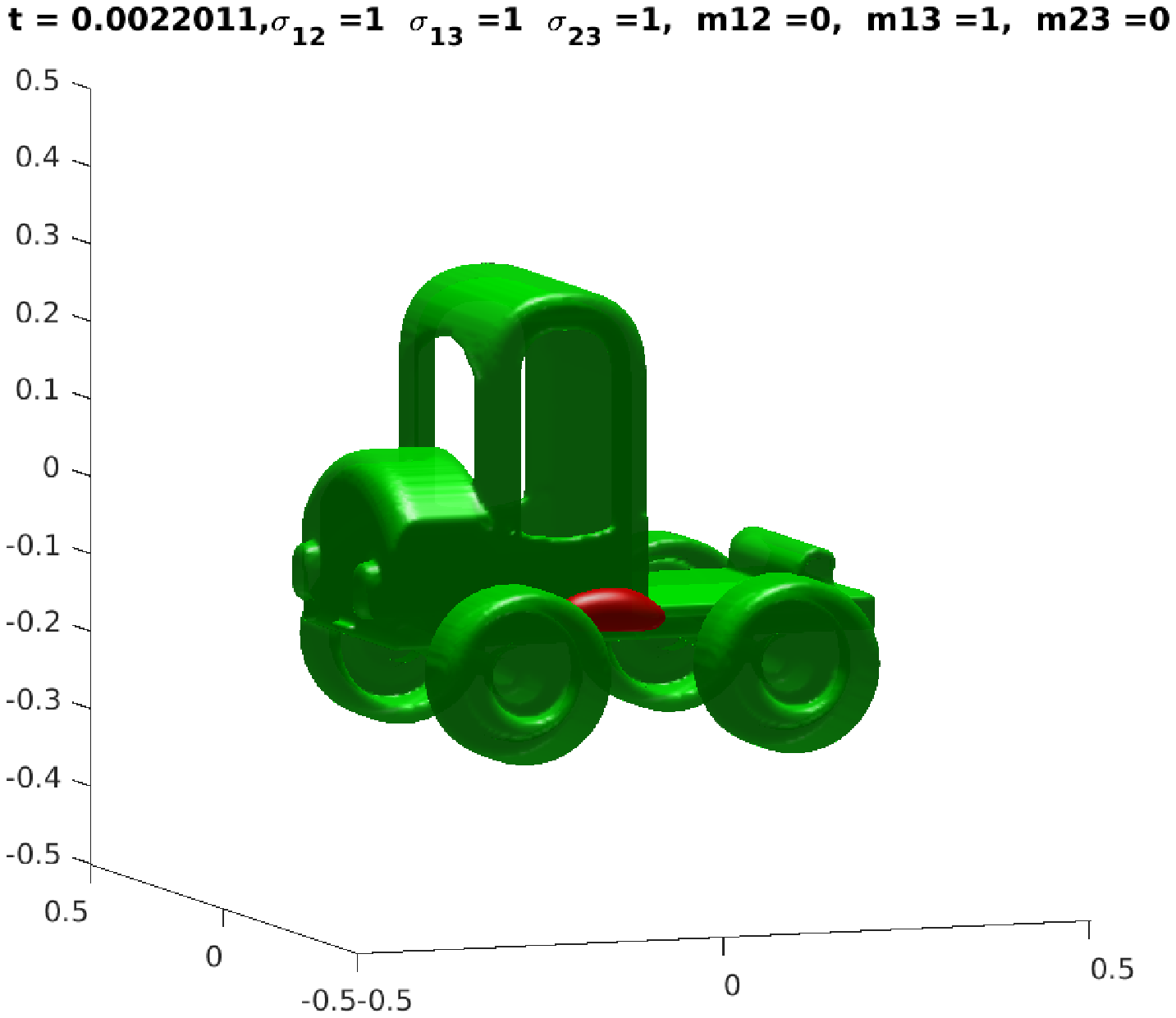} 
\caption{3D multiphase mean curvature flows with homogeneous surface tensions $\sigma_{12} = \sigma_{13} = \sigma_{23} = 1$.
The rows correspond to $(m_{12},m_{13},m_{23}) = (1,1,1)$, $(m_{12},m_{13},m_{23}) = (1,1,0)$, and
$(m_{12},m_{13},m_{23}) = (0,1,0)$, respectively.
{The images show the 1/2-level sets of $u_1$ (red) and $u_2$ (green) at different times}.
}
\label{fig_3D2}
\end{figure}

\section{Conclusion}

We introduced in this paper a numerical scheme for the approximation of multiphase mean curvature flow with 
additive surface tensions and general nonnegative mobilities. The scheme uses a decomposition of the set of mobilities
as sums of harmonically additive mobilities. We provided a formal asymptotic
expansion showing that smooth solutions of the associated Allen-Cahn
system approximate a sharp interface motion driven by
$V_{ij}= m_{i,j}\sigma_{ij}H_{ij}, 1 \leq i<j\leq N$,
up to order 2 in the order parameter $\varepsilon$.
The numerical tests we report are consistent with this expected accuracy.
In particular, when the contrast between mobilities is large, our
scheme provides approximate flows characterized by a width of the diffuse interface 
between phases that is not affected by the mobility contrast.

\section*{Acknowledgments}
The authors thank Roland Denis for fruitful discussions. They acknowledge support from the French National Research Agency (ANR) under grants ANR-18-CE05-0017 (project BEEP) and ANR-19-CE01-0009-01 (project MIMESIS-3D). Part of this work was also supported by the LABEX MILYON (ANR-10-LABX-0070) of Universit\'e de Lyon, within the program "Investissements d'Avenir" (ANR-11-IDEX- 0007) operated by the French National Research Agency (ANR), and by the European Union Horizon 2020 research and innovation programme under the Marie Sklodowska-Curie grant agreement No 777826 (NoMADS).

\end{document}